\documentclass{gtpart}     
\usepackage{amsrefs}
\usepackage{amsfonts}
\usepackage{amscd}
\usepackage{amsxtra} 
\title[Wall-Crossings in Toric Gromov--Witten Theory I: Crepant Examples]
{Wall-Crossings in Toric Gromov--Witten Theory I:\\ Crepant Examples}

\author{Tom Coates}
\givenname{Tom}
\surname{Coates}
\address{Department of Mathematics\\
Imperial College London\\ \newline
180 Queen's Gate\\
London SW7 2AZ \\ UK}
\email{t.coates@imperial.ac.uk}
\urladdr{}

\author{Hiroshi Iritani}
\givenname{Hiroshi}
\surname{Iritani}
\address{Faculty of Mathematics\\
Kyushu University\\ \newline
6-10-1, Hakozaki \\
Higashiku, Fukuoka, 812-8581 \\
Japan}
\email{iritani@math.kyushu-u.ac.jp}
\urladdr{}

\author{Hsian-Hua Tseng}
\givenname{Hsian-Hua}
\surname{Tseng}
\address{Department of Mathematics\\
University of Wisconsin--Madison\\ \newline
Van Vleck Hall, 480 Lincoln Drive\\
Madison, WI 53706-1388\\
USA}
\email{tseng@math.wisc.edu}


\arxivreference{math.AG/0611550}
\arxivpassword{827ah}

\volumenumber{}
\issuenumber{}
\publicationyear{}
\papernumber{}
\startpage{}
\endpage{}
\doi{}
\Zbl{}
\received{}
\revised{}
\accepted{}
\published{}
\publishedonline{}
\proposed{}
\seconded{}
\corresponding{}
\editor{}
\version{}

\theoremstyle{plain}

\newtheorem{thm}{Theorem}[section]
\newtheorem{pro}[thm]{Proposition}
\newtheorem{lem}[thm]{Lemma}

\newtheorem{cor}[thm]{Corollary}

\theoremstyle{definition}

\newtheorem{dfn}[thm]{Definition}
\newtheorem{nt}[thm]{Notation}

\newtheorem{assumptions}[thm]{Assumptions}
\newtheorem{rem}[thm]{Remark}
\newtheorem{exa}[thm]{Example}
\newtheorem*{exa*}{Example}

\newcommand{\PP}{\mathbb{P}}
\newcommand{\FF}{\mathbb{F}}

\newcommand{\CC}{\mathbb{C}}
\newcommand{\ZZ}{\mathbb{Z}}
\newcommand{\RR}{\mathbb{R}}
\newcommand{\QQ}{\mathbb{Q}}
\newcommand{\Cstar}{\CC^\times}
\newcommand{\cD}{\mathcal{D}}
\newcommand{\cF}{\mathcal{F}}
\newcommand{\inertia}{\mathcal{I}}

\newcommand{\cO}{\mathcal{O}}
\newcommand{\ev}{\mathrm{ev}}

\newcommand{\bt}{\mathbf{t}}
\newcommand{\bq}{\mathbf{q}}
\newcommand{\bp}{\mathbf{p}}

\renewcommand{\(}{\left(}
\renewcommand{\)}{\right)}

\DeclareMathOperator{\QC}{QC}

\DeclareMathOperator{\Res}{Res}
\newcommand{\fun}{\mathbf{1}}
\DeclareMathOperator{\age}{age}

\newcommand{\cX}{\mathcal{X}}
\newcommand{\cIX}{\mathcal{IX}}
\newcommand{\cLX}{\mathcal{L}_{\mathcal{X}}}
\newcommand{\cLY}{\mathcal{L}_Y}
\newcommand{\HorbX}{H^\bullet_{\text{\rm orb}}(\cX;\CC)}
\newcommand{\HorbXL}{H^\bullet_{\text{\rm orb}}(\cX;\Lambda)}
\newcommand{\HY}{H^\bullet(Y;\CC)}
\newcommand{\correlator}[1]{\left \langle #1 \right \rangle}

\newcommand{\Correlator}[2]{
  \big \langle \hspace{-1.1ex} \big \langle 
   {#1} 
  \big \rangle \hspace{-1.1ex} \big \rangle^{\cX}_{#2}
}
\newcommand{\BigCorrelator}[2]{
  \Bigg \langle \hspace{-1.3ex}\Bigg \langle 
  {#1} 
  \Bigg \rangle\hspace{-1.3ex} \Bigg \rangle^{\cX}_{#2}}
\DeclareMathOperator{\id}{id}
\DeclareMathOperator{\Eff}{Eff}
\newcommand{\EffX}{\Eff(\cX)}
\newcommand{\VHS}{VSHS\ }
\newcommand{\VHSns}{VSHS}
\newcommand{\U}{\mathbb{U}}
\newcommand{\Ubar}{\overline{\U}}
\newcommand{\tti}{\mathtt{i}}

\newcommand{\fr}[1]{\left\langle #1 \right\rangle}
\newcommand{\fl}[1]{\left\lfloor #1 \right\rfloor}

\newcommand{\CHpair}[2]{\int_{#2} #1}
\newcommand{\Ppair}[2]{\(#1\vphantom{\big|},#2\)}
\newcommand{\OPpair}[2]{\(#1\vphantom{\big|},#2\)_{\text{\rm orb}}}
\newcommand{\Abspair}[2]{\(#1\vphantom{\big|},#2\)_{z\cH_-/\cH_-}}

\DeclareMathOperator{\Aut}{Aut}

\DeclareMathOperator{\Gr}{Gr}

\DeclareMathOperator{\End}{End}

\DeclareMathOperator{\Ker}{Ker}

\theoremstyle{plain}
\newtheorem{conj}[thm]{Conjecture}

\newcommand{\E}{\mathbb{E}}
\newcommand{\JJ}{\mathbb{J}}
\newcommand{\Cz}{\CC \{z \}} 
\newcommand{\Czinv}{\cO\big({\PP^1\setminus \{ 0\} }\big)}
\newcommand{\Czzinv}{\CC \{ z, z^{-1} \} }

\newcommand{\cM}{\mathcal{M}}
\newcommand{\cJ}{\mathcal{J}}
\newcommand{\cI}{\mathcal{I}}
\newcommand{\cE}{\mathcal{E}}

\newcommand{\cH}{\mathcal{H}}
\newcommand{\cHX}{\cH_{\cX}}
\newcommand{\cHY}{\cH_Y}
\newcommand{\cL}{\mathcal{L}}
\newcommand{\cR}{\mathcal{R}}
\newcommand{\cT}{\mathcal{T}}
\newcommand{\ccM}{\overline{\mathcal{M}}}
\newcommand{\tcM}{\widetilde{\cM}}
\newcommand{\hcM}{\widehat{\cM}}
\newcommand{\cOzzinv}[1]{\mathcal{O}_{#1}\{z,z^{-1}\}}
\newcommand{\cOz}[1]{\mathcal{O}_{#1}\{z\}}

\newcommand{\atau}{\mathfrak{t}}
\newcommand{\ap}{\mathfrak{p}}
\newcommand{\aq}{\mathfrak{y}}

\newcommand{\maptau}{\hat{\tau}}
\newcommand{\bof}{\mathbf{f}}

\newcommand{\mapI}{\mathbb{I}}

\newcommand{\Hess}{\operatorname{Hess}}
\newcommand{\KS}{\operatorname{KS}}

\DeclareMathOperator{\mir}{mir}
\DeclareMathOperator{\Mir}{Mir}

\newcommand{\Sym}{\operatorname{Sym}}

\DeclareMathOperator{\Span}{Span}

\def\ov#1{\overline{#1}}
\def\pairW#1#2{\left(#1,#2\right)_W}

\def\parfrac#1#2{\frac{\partial #1}{\partial #2}}

\BibSpec{article}{%
+{}{\sc \PrintAuthors} {author}
+{,}{ \textrm} {title}
+{,}{ \textit} {journal}
+{,}{ } {volume}
+{}{ \IfEmptyBibField{journal}{}{\PrintDate{date}}} {transition}
+{,}{ } {pages}
+{,}{ } {status}
+{.}{} {transition}
+{}{ Preprint\IfEmptyBibField{date}{}{ \PrintDate{date}}, available at \tt} {eprint}
+{.}{} {transition}
}

\BibSpec{book}{%
+{}{\sc \PrintAuthors} {author}
+{,}{ \textit} {title}
+{.}{ } {series}
+{,}{ } {volume}
+{.}{ } {publisher}
+{,}{ } {place}
+{,}{ } {date}
+{.}{} {transition}
}

\BibSpec{collection.article}{
+{} {\sc \PrintAuthors} {author}
+{,} { \textrm } {title}
+{,} { in \it \PrintConference} {conference}
+{,} { pp.~} {pages}
+{.} {\PrintBook} {book}
+{,} { \PrintDateB} {date}
+{.} {} {transition}
}

\BibSpec{innerbook}{
+{.} { \emph} {title}
+{.} { } {part}
+{:} { \emph} {subtitle}
+{.} { } {series}
+{,} { } {volume}
+{.} { Edited by \PrintNameList} {editor}
+{.} { Translated by \PrintNameList}{translator}
+{.} { \PrintContributions} {contribution}
+{.} { } {publisher}
+{.} { } {organization}
+{,} { } {address}
+{,} { \PrintEdition} {edition}
+{,} { \PrintDateB} {date}
+{.} { } {note}
}

\begin{document}

\begin{abstract}
  Let $\cX$ be a Gorenstein orbifold with
  projective coarse moduli space $X$ and let $Y$ be a crepant
  resolution of $X$.  We state a conjecture relating the genus-zero
  Gromov--Witten invariants of $\cX$ to those of $Y$, which differs in
  general from the Crepant Resolution Conjectures of Ruan and
  Bryan--Graber, and prove our conjecture when $\cX = \PP(1,1,2)$ and
  $\cX = \PP(1,1,1,3)$.  As a consequence, we see that the original
  form of the Bryan--Graber Conjecture holds for $\PP(1,1,2)$ but is
  probably false for $\PP(1,1,1,3)$. Our methods are based on mirror
  symmetry for toric orbifolds.
\end{abstract}

\maketitle


\section{Introduction}

In this paper we use mirror symmetry to determine the relationship
between the quantum orbifold cohomology of an orbifold $\cX$ and the
quantum cohomology of a crepant resolution $Y$ of $\cX$ in the cases
$\cX = \PP(1,1,2)$ and $\cX = \PP(1,1,1,3)$.

\subsection*{A Picture From Physics}

Quantum cohomology and quantum orbifold cohomology occur in string
theory as a small part of a much larger picture.  There is supposed to
be a moduli space of physical theories --- the \emph{stringy K\"ahler
  moduli space} $\cM$ --- and a bundle of algebras over this moduli
space formed by the chiral rings of the theories.  Near certain limit
points of $\cM$, called \emph{large radius limit points} or
\emph{cusps}, the bundle of algebras is given by the quantum
cohomology or quantum orbifold cohomology of a target space $\cX$; at
a general point of $\cM$, however, there will be no such description.

More precisely, near each cusp there are distinguished co-ordinates on
$\cM$, called \emph{flat co-ordinates}, and a distinguished
trivialization of the bundle, called a \emph{flat trivialization},
such that when expressed in flat co-ordinates and with respect to the
flat trivialization the bundle of algebras is isomorphic to the
quantum orbifold cohomology algebra of $\cX$.  From this point of
view, the \emph{quantum parameters} $q_i$ occurring in the definition
of the quantum product\footnote{See Section~\ref{sec:smallQC} below
  for the definition.}  are exponentiated flat co-ordinates on a
neighbourhood of the corresponding cusp in $\cM$.

Different cusps can correspond to different target spaces, and in
particular the quantum orbifold cohomology $\QC(\cX)$ of $\cX$ and the
quantum cohomology $\QC(Y)$ of a crepant resolution $Y$ of $\cX$ are
expected to come from different cusps of the same moduli space $\cM$.
Since $\QC(\cX)$ and $\QC(Y)$ are supposed to be parts of the same
global family of algebras this motivates the conjecture, made in
various forms by various authors and discussed in detail below, that
$\QC(\cX)$ and $\QC(Y)$ coincide after analytic continuation in
quantum parameters.  We have already seen, however, the first hint
that this conjecture is probably too na\"\i ve in general: one should
also take into account whether or not the flat co-ordinates near the cusps
associated to $\cX$ and to $Y$ coincide after analytic continuation.

\subsection*{Overview of our Results}

In what follows we build on work of Givental and Barannikov to
construct a rigorous version of this physical picture in the cases
$\cX = \PP(1,1,2)$ and $\cX = \PP(1,1,1,3)$.  Our key tool is
(mathematical) mirror symmetry for toric orbifolds.  Rather than give
a global construction of the stringy K\"ahler moduli space $\cM$ ---
we do not know how to do this --- we instead construct the so-called
\emph{B-model moduli space} $\cM_{\rm B}$ and then identify subsets near
certain cusps in $\cM_{\rm B}$ with the subsets of $\cM$ on which the
quantum cohomology of $Y$ and the quantum orbifold cohomology of $\cX$
are defined.  The B-model moduli space $\cM_{\rm B}$ is expected to
coincide under (string theoretic) mirror symmetry with the stringy
K\"ahler moduli space $\cM$, but it has the advantage that we can give
it rigorous mathematical meaning.

We construct $\cM_{\rm B}$ from the toric data --- it is the toric orbifold
associated to the secondary fan for the crepant resolution $Y$ of
$\cX$.  Rather than constructing just a family of algebras over
$\cM_{\rm B}$ we construct a significantly finer structure called a
\emph{variation of semi-infinite Hodge structure} or \VHSns.  This
\VHS determines, as we will see in Section~\ref{sec:VHS} below, a
family of algebras over $\cM_{\rm B}$.  It also, together with some extra
data canonically associated to each cusp, determines flat co-ordinates
and a flat trivialization near each cusp, and allows us to compare the
flat structures associated to different cusps.  The \VHS here consists
of a vector bundle $V\to \cM_{\rm B}$ with flat connection and a family of
subspaces $E_y$,~$y \in \cM_{\rm B}$, in the fibers of $V$.  The vector
bundle $V$ is infinite-dimensional and the subspaces $E_y$ are in an
appropriate sense\footnote{We consider $E_y$ to be an element of the
  Segal--Wilson Grassmannian of $V_y$.}
semi-infinite.  The family of subspaces $\{E_y\}$ is an analog of a
variation of Hodge structure and it satisfies a version of Griffiths
transversality --- see Section~\ref{sec:VHS}.  The extra data at each
cusp consists of an \emph{opposite subspace} and a \emph{dilaton
  shift}.  Fix a cusp in $\cM_{\rm B}$ and choose $y$ near that cusp.  By
parallel transport one can, for each $x \in \cM_{\rm B}$, regard $E_{x}$ as
a subspace of the fiber $V_y$ and thus define a \emph{limiting Hodge
  structure} $E_{\text{lim}} \subset V_y$ associated\footnote{Here
  $E_{\text{lim}}$ is roughly speaking the limit of $E_x \subset V_y$
  as $x$ approaches the cusp.} to the cusp.  The opposite subspace
associated to the cusp is a subspace $V_-$ of $V_y$ such that
$E_{\text{lim}} \oplus V_- = V_y$, and the dilaton shift is a non-zero
element of $E_{\text{lim}}$.  The opposite subspace $V_-$ is uniquely
determined by monodromy properties --- it is required to be invariant
under the local monodromy around the cusp --- and a homogeneity
condition (Theorem~\ref{thm:uniqueness}); $V_-$ is the analog of the
\emph{weight filtration} on a limiting mixed Hodge structure.

We will see below that when $\cX = \PP(1,1,2)$, so its crepant
resolution $Y$ is the Hirzebruch surface $\FF_2$, the opposite
subspaces at the cusps of $\cM_{\rm B}$ associated to $\cX$ and to $Y$ agree
under parallel transport in $V$.  This implies that the flat
structures determined by $\cX$ and $Y$ agree: that not only the
families of algebras $\QC(\cX)$ and $\QC(Y)$, but also the flat
trivializations and flat co-ordinates associated to $\cX$ and $Y$, are
related by analytic continuation.  We deduce:

\begin{thm} \label{thm:P112}
  Let $\cX = \PP(1,1,2)$ and $Y = \FF_2$.  There is a linear
  isomorphism $\Theta:\HorbX \to \HY$ between the Chen--Ruan orbifold
  cohomology of $\cX$ and the cohomology of $Y$ such that the small
  quantum orbifold cohomology algebra of $\cX$ with quantum parameter
  $q$ and the small quantum cohomology algebra of $Y$ with quantum
  parameters $q_1$,~$q_2$ are isomorphic via $\Theta$, after analytic
  continuation in $(q_1,q_2)$ and the substitution
  \begin{align}
    \label{eq:P112subst}
    q_1 = -1 && q_2 = \tti \sqrt{q}.
  \end{align}
  An explicit formula for $\Theta$ is given as
  equation~\eqref{eq:mapforP112} below.  Furthermore, the map $\Theta$
  and the specialization \eqref{eq:P112subst} identify the quantum
  cohomology Frobenius manifolds associated to $\cX$ and $Y$.
\end{thm}

We will see further that when $\cX = \PP(1,1,1,3)$, so its crepant
resolution $Y$ is the scroll $\FF_3$, the opposite subspaces at the
cusps of $\cM_{\rm B}$ associated to $\cX$ and to $Y$ do not agree under
parallel transport in $V$.  This implies that the flat structures
determined by $\cX$ and $Y$ \emph{are different}: they do not agree
under analytic continuation.  We have:

\begin{thm} \label{thm:P1113}
  Let $\cX = \PP(1,1,1,3)$ and $Y = \FF_3$.  There is a linear
  isomorphism $\Theta(q)$, which depends non-trivially on $q$, between
  $\HorbX$ and $\HY$ such that the small quantum orbifold cohomology
  algebra of $\cX$ with quantum parameter $q$ and the small quantum
  cohomology algebra of $Y$ with quantum parameters $q_1$,~$q_2$ are
  isomorphic via $\Theta(q)$, after analytic continuation in $(q_1,q_2)$
  followed by the substitution
  \begin{align}
    \label{eq:P1113subst}
    q_1 = 1 && q_2 = \sqrt[3]{q}.
  \end{align}
  The isomorphism $\Theta(q)$ matches the Poincar\'e pairing on $\HY$
  with the orbifold Poincar\'e pairing on $\HorbX$.
\end{thm}

An explicit formula for $\Theta(q)$ can be found in
Section~\ref{sec:proofP1113}.  Note that the isomorphism $\Theta(q)$
cannot be induced by any isomorphism of the quantum cohomology
Frobenius manifolds associated to $\cX$ and $Y$, as it depends
non-trivially on $q$.

\subsection*{Mirror Symmetry}

Let us call the \VHS which we construct the \emph{B-model \VHSns}.
Mirror symmetry identifies the B-model \VHS with an object familiar
in Gromov--Witten theory.  Givental has observed that if one encodes
genus-zero Gromov--Witten invariants of $\cX$ in a
certain\footnote{This encoding is described in
  Section~\ref{sec:Giventalapproach} below; here the key point is that
  knowing $\cLX$ is equivalent to knowing all genus-zero
  Gromov--Witten invariants of $\cX$.} Lagrangian submanifold $\cLX$
of a symplectic vector space $\cHX$ then many seemingly-complicated
statements in Gromov--Witten theory are in fact simple geometric
assertions about $\cLX$ \cite{Givental:symplectic}.  Mirror symmetry
identifies the B-model \VHS with the so-called A-model \VHSns, which
consists of \emph{the family of tangent spaces to Givental's
  Lagrangian submanifold} $\cLX$.  The A-model \VHS is canonically
trivialized --- it is a family of subspaces of a fixed vector space
$\cHX$ --- and mirror symmetry here asserts that the B-model \VHSns,
expressed with respect to the flat trivialization near the cusp
corresponding to $\cX$, coincides with the family of tangent spaces to
$\cLX$.  In the cases at hand this follows from mirror theorems due
to Givental \cite{Givental:toric} and Coates--Corti--Lee--Tseng
\cite{CCLT}; it implies in particular that the family of algebra
structures over $\cM_{\rm B}$ determined by the B-model \VHS coincides
near the cusps of $\cM_{\rm B}$ with the quantum orbifold cohomology of
$\cX$ and $Y$.

\subsection*{The Crepant Resolution Conjecture}

The results we prove have the following consequence when $\cX =
\PP(1,1,2)$ and $\PP(1,1,1,3)$.  We conjecture that this holds in
general.

\begin{conj} \label{conj} Suppose that $\cX$ is an orbifold with
  projective coarse moduli space $X$, and that $Y$ is a crepant
  resolution of $X$.  Let $\cLX \subset \cHX$ be Givental's Lagrangian
  submanifold for $\cX$, and let $\cLY \subset \cHY$ be Givental's
  Lagrangian submanifold for $Y$.  Then there exists a linear
  symplectic isomorphism $\U:\cHX \to \cHY$, satisfying the conditions
  enumerated in \eqref{eq:CRCconditions} below, such that after
  analytic continuation of $\cLX$ and $\cLY$ we have $\U(\cLX) =
  \cLY$.
\end{conj}

We have not defined $\cLX$ and $\cLY$ at this point, so
Conjecture~\ref{conj} is necessarily slightly vague; we give a precise
statement as Conjecture~\ref{conj:preciseCRC} below.  As we will see
in Sections~\ref{sec:P1113} and \ref{sec:P112}, the symplectic
transformation $\U$ here records the effect of parallel transport in
the B-model \VHS (i.e. of parallel transport in the fibers of the
vector bundle $V \to \cM_{\rm B}$).

Conjecture~\ref{conj} is our version of the Crepant Resolution
Conjecture.  We now discuss its relationship to earlier versions of
the Crepant Resolution Conjecture formulated by Ruan and by
Bryan--Graber.  For the rest of this section, let $\cX$ be an orbifold
with projective coarse moduli space $X$ and let $Y$ be a crepant
resolution of $X$.

The first attempt to describe the relationship between the small
quantum cohomology algebra of $Y$ and the small quantum orbifold
cohomology algebra of $\cX$ is due to Ruan.  He conjectured that the
small quantum orbifold cohomology of $\cX$ is isomorphic to the
algebra obtained from the small quantum cohomology of $Y$ by analytic
continuation in quantum parameters followed by specializing some of
those parameters to roots of unity.  Theorems~\ref{thm:P112}
and~\ref{thm:P1113} prove the Ruan Conjecture for $\cX = \PP(1,1,2)$
and $\cX = \PP(1,1,1,3)$; the relationship between our Conjecture and
the Ruan Conjecture is discussed further in Section~\ref{sec:RuanConjecture}. 

A significant strengthening of the Ruan Conjecture has been proposed
by Bryan--Graber \cite{Bryan--Graber}.  They have modified their
conjecture in the light of the examples in this paper, but initially
they asserted that the quantum cohomology Frobenius manifolds
associated to $\cX$ and $Y$ become isomorphic after analytic
continuation in quantum parameters.  Thus Theorem~\ref{thm:P112}
proves the original form of the Bryan--Graber Conjecture for
$\PP(1,1,2)$, but Theorem~\ref{thm:P1113} \emph{does not} prove the
original form of the Bryan--Graber Conjecture for $\PP(1,1,1,3)$.  In
Theorem~\ref{thm:hardLefschetz} below we show that if the orbifold
cohomology of $\cX$ satisfies a Hard Lefschetz property --- this
property holds for $\cX = \PP(1,1,2)$ but not for $\cX = \PP(1,1,1,3)$
--- then our Conjecture implies the original form of the Bryan--Graber
Conjecture.  In general, however, our Conjecture \emph{does not} imply
the original Bryan--Graber Conjecture and we expect that the latter is
false.  The most recent version of the Bryan--Graber Conjecture
includes the Hard Lefschetz condition as a hypothesis.

We should emphasize that our results here \emph{do not} show that $\cX
= \PP(1,1,1,3)$ is a counterexample to the original form of the
Bryan--Graber Conjecture.  It is possible that there is a different
path of analytic continuation and a different choice of specialization
\eqref{eq:P1113subst} which produces an isomorphism of Frobenius
manifolds.  But we think that this is unlikely.  Conjecture~\ref{conj}
expresses the relationship between the quantum cohomology algebras of
$\cX=\PP(1,1,1,3)$ and $Y=\FF_3$ which is forced upon us by mirror
symmetry; in this sense it is the natural conjecture to make.
Furthermore the original form of the Bryan--Graber Conjecture ignores
some flexibility in parts of the structure --- the \emph{flat
  trivialization} and \emph{flat co-ordinates} --- which topological
string theory suggests should be background dependent rather than
fixed.  So we see no compelling reason for the original form of
the Bryan--Graber Conjecture to hold.  Conjecture~\ref{conj} has been
proved in a number of local toric Calabi--Yau examples by Coates
\cite{Coates:wallcrossings2}, and in forthcoming work Iritani will
prove it for general toric crepant birational transformations
\cite{Iritani:wallcrossings3}.

\subsection*{Singularity Theory}

Our results also have consequences in singularity theory. 
We construct the B-model \VHS from a so-called 
\emph{Landau--Ginzburg model}.  Singularity theorists have long known
how to construct the germ of 
a Frobenius manifold from a Landau--Ginzburg model:
for local singularities (germs of isolated
hypersurface singularities) 
this is due to Kyoji Saito \cite{saitoK} and Morihiko
Saito \cite{saitoM}; for global singularities (our case) 
this is due to Douai--Sabbah \cite{Douai--Sabbah}.  
It has long been known also that 
there are in general many possible germs of Frobenius structures
for a given singularity:
in our language, this is the statement that 
one can choose from many possible opposite subspaces. 
From this point of view, the content of this paper 
is that more global considerations --- monodromy and homogeneity
properties --- single out a canonical opposite subspace associated to
each cusp, and that the opposite subspaces 
associated to different cusps 
can be compared via analytic continuation. 

\subsection*{Plan of the Paper}

In Section~\ref{sec:formalism} we fix notation and develop our general
theory: we define variations of semi-infinite Hodge structure,
introduce Givental's symplectic formalism, and explain what we mean by
mirror symmetry.  In Section~\ref{sec:P1113} we analyze the case $\cX
= \PP(1,1,1,3)$, proving Theorem~\ref{thm:P1113} and
Conjecture~\ref{conj}.  The argument which proves
Theorem~\ref{thm:P112} and Conjecture~\ref{conj} for $\cX =
\PP(1,1,2)$ is very similar and we summarize it in
Section~\ref{sec:P112}.  In Section~\ref{sec:CRC} we describe a more
detailed version of Conjecture~\ref{conj} and prove that it implies
the most recent form of the Bryan--Graber Conjecture.  We conclude
with an Appendix A describing the Mellin--Barnes method for analytic
continuation of hypergeometric functions.

\subsection*{Acknowledgements}

This paper was begun whilst the authors held Postdoctoral Fellowships
at the Mathematical Sciences Research Institute as part of the ``New
Topological Structures in Physics'' program.  We are grateful to MSRI
for financial support and for providing such a stimulating research
environment, and to Yongbin~Ruan for many inspiring discussions.  This
paper forms part of a larger joint project with Alessio~Corti
\cite{CCLT}\cite{CCIT:computing}: we thank him for a
productive and enjoyable collaboration, and for many useful
conversations.  We also thank Mina~Aganagic, Vincent~Bouchard,
Alexander~Givental, and Andrei~Mustata for helpful conversations.
Tom~Coates was partially supported by the Clay Mathematics Institute,
the Royal Society, and by NSF grant DMS-0401275.  
Hiroshi~Iritani was partially supported by 
Grant-in-Aid for JSPS Fellows, 
Scientific Research 18-15108 and 
EPSRC EP/E022162/1. 
Hsian-Hua~Tseng was partially supported by a
postdoctoral fellowship at the Pacific Institute for the Mathematical
Sciences.

\section{Variations of Semi-Infinite Hodge Structure}

\label{sec:formalism}
In this section we fix notation for Gromov--Witten invariants, give an
introduction to Barannikov's theory of variations of semi-infinite
Hodge structure, indicate how this meshes with Givental's geometric
approach to Gromov--Witten theory, and explain what we mean by mirror
symmetry.  We assume that the reader is familiar with quantum
cohomology and quantum orbifold cohomology.  The quantum cohomology
and Gromov--Witten theory of algebraic varieties have been quite
widely studied: good introductions to the field include
Fulton--Pandharipande \cite{Fulton--Pandharipande}, Cox--Katz
\cite{Cox--Katz}, and Hori \emph{et al.} \cite{ClayVolume}.  The
quantum cohomology and Gromov--Witten theory of orbifolds were
introduced into mathematics by Chen and Ruan
\cite{Chen--Ruan:orbifold, Chen--Ruan:GW} in the setting of symplectic
geometry; an algebro-geometric version of the theory has been
developed by Abramovich, Graber, and Vistoli \cite{AGV:1, AGV:2}.  An
overview of this material, in compatible notation, can be found in
Section~2 of Coates--Corti--Lee--Tseng \cite{CCLT}.  Givental
introduced his formalism in \cite{Givental:quantization} and gave an
expository account of it in \cite{Givental:symplectic}.

\subsection{Notation and Conventions}

\label{sec:notation}

We work in the algebraic category and over $\CC$: by ``manifold'' we
mean ``smooth projective algebraic variety'' and we use the terms
``orbifold'' and ``smooth Deligne--Mumford stack'' interchangeably.
Introduce notation as follows.
\begin{center}
  \begin{tabular}{lp{0.75\textwidth}}
    $\cX$ & a compact orbifold \\
    $\cIX$ & the inertia stack of $\cX$.  A point of $\cIX$ is a pair
    $(x,g)$ with $x$ a point of $\cX$ and $g \in \Aut_{\cX}(x)$\\
    $I$ & the involution of $\cIX$ which sends $(x,g)$ to
    $(x,g^{-1})$ \\
    $\HorbX$ & the Chen--Ruan orbifold cohomology groups of $\cX$.  These are the
    cohomology groups $H^\bullet(\cIX;\CC)$ of the inertia stack  \\
    $\age$ & a rational number associated to each component $\cX_i$ 
    of the inertia stack. The grading on orbifold cohomology is shifted
    by the age: $\alpha\in H^p(\cX_i;\CC)$ has degree
    $\deg\alpha=p + 2 \age(\cX_i)$ \\
    $\OPpair{\alpha}{\beta}$ & the orbifold Poincar\'e pairing
    $\int_{\cIX} \alpha \cup I^\star \beta$ \\
    $\EffX$  & the set  of degrees of representable maps from possibly-stacky curves to
    $\cX$ (i.e. of degrees of effective curves in $\cX$) \\
    $\fr{r}$ &  the fractional part  $r - \fl{r}$ of a rational number $r$
  \end{tabular}
\end{center}

\begin{exa}
  Weighted projective space $\PP(w_0,w_1,\dots,w_n)$ is the stack
  quotient $\big[ \big( {\CC^{n+1} - \{0\}} \big) / \Cstar \big]$ where
  $\Cstar$ acts with weights $-w_0,\ldots,-w_n$.  Components of the
  inertia stack of $\PP(w_0,\dots,w_n)$ are indexed by
  \[
  F = \left\{{\textstyle {k \over w_i}} \, \Big\vert \, \, 0 \leq k <
    w_i, \; 0 \leq i \leq n\right\}
  \]
  via:
  \[
  \inertia \PP(w_0,\ldots,w_n) = \coprod_{f \in F} \, \PP\big(V^f\big)
  \]
  Here
  \[
    V^f = \left\{\(x_0,\ldots,x_n\) \in \CC^{n+1} \mid \text{$x_i = 0$
        unless $w_i f \in \ZZ$} \right\} 
    \]
    and
    \[
    \PP\big(V^f\big) = \left[ \( {V^f - \{0\}} \) / \Cstar \right]
    \]
    so that $\PP\big(V^f\big)$ consists of those points of
    $\PP(w_0,\ldots,w_n)$ with isotropy group containing $\exp(2 \pi i
    f) \in \Cstar$.  The locus $\PP\big(V^f\big)$ is itself a weighted
    projective space.  The involution $I$ maps the component
    $\PP\big(V^f\big)$ to the component $\PP\big(V^{\langle -f \rangle}\big)$.  The
    age of $\PP\big(V^f\big) \subset \inertia$ is $\langle - w_0 f \rangle +
    \cdots + \langle - w_n f \rangle$.
\end{exa}

\begin{exa}
  The orbifold cohomology of $\PP(1,1,2)$ is
  \[
  H^\bullet_{\text{\rm orb}}\big(\PP(1,1,2);\CC\big)  = H^\bullet\big(\PP\big(V^0\big);\CC\big) \oplus
  H^{\bullet-1}\big(\PP\big(V^{1 \over 2}\big);\CC\big)
  \]
  where:
  \begin{align*}
    &\PP\big(V^0\big) = \PP(1,1,2) &
    \age &= 0 \\
    &\PP\big(V^{1 \over 2}\big) = \PP(2) &
    \age &= 1 
  \end{align*}
  The involution $I$ is trivial.
\end{exa}

\begin{exa}
  The orbifold cohomology of $\PP(1,1,1,3)$ is
  \[
  H^\bullet_{\text{\rm orb}}\big(\PP(1,1,1,3);\CC\big)  = H^\bullet\big(\PP\big(V^0\big);\CC\big) \oplus
    H^{\bullet-2}\big(\PP\big(V^{1 \over 3}\big);\CC\big) \oplus
    H^{\bullet-1}\big(\PP\big(V^{2 \over 3}\big);\CC\big)
  \]
  where:
  \begin{align*}
    &\PP\big(V^0\big) = \PP(1,1,1,3) &
    \age &= 0 \\
    &\PP\big(V^{1 \over 3}\big) = \PP(3) &
    \age &= 2 \\   
    &\PP\big(V^{2 \over 3}\big) = \PP(3) &
    \age &= 1   
  \end{align*}
  The involution $I$ exchanges $\PP(V^{1 \over 3})$ and $\PP(V^{2 \over 3})$.
\end{exa}

\subsubsection{Generators and Bases for Homology and Orbifold
  Cohomology}
\label{sec:cohomology}

We now fix notation for the homology and orbifold cohomology of the
spaces which we will consider.  When discussing $\cX = \PP(1,1,2)$:
\begin{center}
  \begin{tabular}{lp{0.75\textwidth}}    
    $\fun_0$ & is the fundamental class of $\PP\big(V^0\big)$ \\
    $\fun_{1 \over 2}$ & is the fundamental class of $\PP\big(V^{1 \over 2}\big)$ \\ 
    $p$ & is the first Chern class $c_1(\cO(1)) \in H^2(\cX;\CC)$
  \end{tabular}
\end{center}
and:
\begin{align*}
   \phi_0 & = \fun_0 &
  \phi_1 & = p&
  \phi_2 & = p^2&
  \phi_3 & = \fun_{1 \over 2}
 \\
  \phi^0 & = 2 p^2& 
  \phi^1 & = 2 p&
  \phi^2 & = 2 \fun_0&
  \phi^3 &= 2 \fun_{1 \over 2}
\end{align*}
Note that, here and below, $\OPpair{\phi_i}{\phi^j} = {\delta_i}^j$.

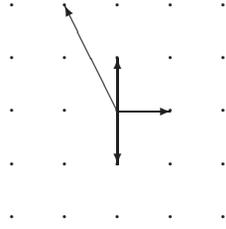
\begin{figure}[ht!]
  \begin{center}
    \begin{picture}(90,90)(-45,-45)
      \multiput(-40,-40)(20,0){5}{\makebox(0,0){$\cdot$}}
      \multiput(-40,-20)(20,0){5}{\makebox(0,0){$\cdot$}}
      \multiput(-40,-0)(20,0){5}{\makebox(0,0){$\cdot$}}
      \multiput(-40,20)(20,0){5}{\makebox(0,0){$\cdot$}}
      \multiput(-40,40)(20,0){5}{\makebox(0,0){$\cdot$}}
      \put(0,0){\vector(1,0){20}}
      \put(0,0){\vector(0,1){20}}
      \put(0,0){\vector(-1,2){20}}
      \put(0,0){\vector(0,-1){20}}
    \end{picture}
    \caption{The fan for the toric variety $\FF_2$.}
    \label{fig:fanF2}
  \end{center}
\end{figure}

When discussing the Hirzebruch surface $\FF_2$, which is the
projective bundle $\PP(\cO(-2) \oplus \cO)$ over $\PP^1$:
\begin{center}
  \begin{tabular}{lp{0.75\textwidth}}    
    $p_1$ & is the class in $H^2(\FF_2;\CC)$ Poincar\'e-dual to a fiber
    \\
    $p_2$ & is the class in $H^2(\FF_2;\CC)$ Poincar\'e-dual to the
    infinity section
  \end{tabular}
\end{center}
and:
\begin{align*}
   \phi_0 & = 1 &
  \phi_1 & = p_1&
  \phi_2 & = p_2&
  \phi_3 & = p_1 p_2
 \\
  \phi^0 & = p_1 p_2& 
  \phi^1 & = p_2&
  \phi^2 & = p_1&
  \phi^3 &= 1
\end{align*}
The surface $\FF_2$ is the toric variety corresponding to the fan in
Figure~\ref{fig:fanF2}.

It can therefore (see for example \cite[Chapter~VII]{Audin}) be
constructed as the quotient of $\(\CC^2-\{0\} \) \times \(\CC^2-\{0\}
\)$ by the action of $\big(\Cstar\big)^2$
\[
(s,t) \colon 
\begin{pmatrix}
  x \\ y \\ z \\ w
\end{pmatrix}
\longmapsto 
\begin{pmatrix}
  s x \\ s y \\ s^{-2} t z \\ t w
\end{pmatrix}
\]
and its cohomology ring is
\[
H^\bullet(\FF_2;\CC) = \CC[p_1,p_2] \slash \left \langle p_1^2, p_2^2 - 2 p_1 p_2
\right \rangle.
\]

When discussing $\cX = \PP(1,1,1,3)$:
\begin{center}
  \begin{tabular}{lp{0.75\textwidth}}    
    $\fun_0$ & is the fundamental class of $\PP\big(V^0\big)$ \\
    $\fun_{1 \over 3}$ & is the fundamental class of $\PP\big(V^{1 \over 3}\big)$ \\ 
    $\fun_{2 \over 3}$ & is the fundamental class of $\PP\big(V^{2 \over 3}\big)$ \\ 
    $p$ & is the first Chern class $c_1(\cO(1)) \in H^2(\cX;\CC)$
  \end{tabular}
\end{center}
and:
\begin{align*}
   \phi_0 & = \fun_0 &
  \phi_1 & = p&
  \phi_2 & = p^2&
  \phi_3 & = p^3&
  \phi_4 & = \fun_{1 \over 3} &
  \phi_5 & = \fun_{2 \over 3}
 \\
   \phi^0 & = 3 p^3 &
  \phi^1 & = 3 p^2&
  \phi^2 & = 3 p&
  \phi^3 & = 3 \fun_0&
  \phi^4 & = 3 \fun_{2 \over 3} &
  \phi^5 & = 3 \fun_{1 \over 3}
\end{align*}

When discussing the projective bundle $\FF_3 = \PP(\cO(-3) \oplus
\cO)$ over $\PP^2$:
\begin{center}
  \begin{tabular}{lp{0.75\textwidth}}    
    $p_1$ & is the class in $H^2(\FF_3;\CC)$ Poincar\'e-dual to the
    preimage in $\FF_3$ of a hyperplane in $\PP^2$
    \\
    $p_2$ & is the class in $H^2(\FF_3;\CC)$ Poincar\'e-dual to the
    infinity section
    \\
  \end{tabular}
\end{center}
and: 
\begin{align*}
  \phi_0 & = 1 & \phi^0 & = p_1^2p_2 \\
  \phi_1 & = \frac{p_2}{3}&   \phi^1 & = p_1 p_2 \\
  \phi_2 & = \frac{p_1p_2}{3}&   \phi^2 & = p_2 \\
  \phi_3 & = \frac{p_2-3p_1}{3}  &   \phi^3 & = -p_1(p_2-3p_1) \\
  \phi_4 & = -\frac{p_1(p_2-3p_1)}{3}&   \phi^4 & = p_2-3p_1 \\
  \phi_5 & = \frac{p_1^2 p_2}{3} & \phi^5 & = 3 
\end{align*}
The scroll $\FF_3$ is the toric variety corresponding to the fan with
rays
\[
e_1 = \begin{pmatrix}  1 \\ 0 \\ 0 \end{pmatrix},
e_2 = \begin{pmatrix}  0 \\ 1 \\ 0 \end{pmatrix},
e_3 = \begin{pmatrix}  -1 \\ -1 \\ 3 \end{pmatrix},
e_4 = \begin{pmatrix}  0 \\ 0 \\ 1 \end{pmatrix},
e_5 = \begin{pmatrix}  0 \\ 0 \\ -1 \end{pmatrix}
\]
and three-dimensional cones spanned by 
\[
\{e_1,e_2,e_4\},
\{e_1,e_3,e_4\},
\{e_2,e_3,e_4\},
\{e_1,e_2,e_5\},
\{e_1,e_3,e_5\},
\{e_2,e_3,e_5\}.
\]
It can be constructed as the quotient of $\(\CC^3-\{0\} \) \times
\(\CC^2-\{0\} \)$ by the action of $\(\Cstar\)^2$
\[
(s,t) \colon 
\begin{pmatrix}
  x \\ y \\ z \\ u \\ v
\end{pmatrix}
\longmapsto 
\begin{pmatrix}
  s x \\ s y \\ s z \\ s^{-3} t u \\ t v
\end{pmatrix}
\]
and its cohomology ring is
\[
H^\bullet(\FF_3;\CC) = \CC[p_1,p_2] \slash \left \langle p_1^3, p_2^2 - 3 p_1 p_2
\right \rangle.
\]

\subsubsection{Gromov--Witten Invariants and Quantum Cohomology}
\label{sec:GWdefs}
We denote Gromov--Witten invariants using correlators, writing, with
notation as in \cite{CCLT}:
  \begin{equation}
    \label{eq:GW}
    \correlator{\alpha_1 \psi^{k_1},\ldots,\alpha_n
      \psi^{k_n}}^{\cX}_{g,n,d}
    = \int_{\cX_{g,n,d}^{vir}} \prod_{i=1}^n \ev_i^\star \alpha_i \cdot
    \psi_i^{k_i}
  \end{equation}
  The integral here means cap product with the virtual fundamental
  class.  If any of the $k_i$ are non-zero then \eqref{eq:GW} is
  called a \emph{gravitational descendant}. 

  Double correlators denote generating functions for Gromov--Witten
  invariants:
\begin{equation}
 \label{eq:doublecorrelator}
  \Correlator{\alpha_1 \psi^{i_1},\ldots,\alpha_k
    \psi^{i_k}}{\tau} = 
  \sum_{d \in \EffX} \sum_{n \geq 0} {Q^d \over n!}
  \correlator{\alpha_1 \psi^{i_1},\ldots,\alpha_k
    \psi^{i_k},\tau,\tau,\ldots,\tau}^{\cX}_{0,n+k,d}
\end{equation}
where $\tau \in \HorbX$ and:
\[
Q^d =
\begin{cases}
  Q^{\CHpair{p}{d}} & \text{$\cX = \PP(1,1,2)$ or $\PP(1,1,1,3)$} \\
  Q_1^{\CHpair{p_1}{d}} Q_2^{\CHpair{p_2}{d}} & \text{$\cX = \FF_2$ or
    $\FF_3$}
\end{cases}
\]
These generating functions are formal series in the co-ordinates
$\tau^1, \tau^2,\ldots,\tau^N$ of $\tau = \tau^1 \phi_1 + \ldots +
\tau^N \phi_N$ and the variables $Q^{1/2}$ or $Q^{1/3}$ or $Q_1$, $Q_2$.
We make these latter variables, which are included to make the series
\eqref{eq:doublecorrelator} converge, into elements of our ground
ring $\Lambda$, setting:
\[
\Lambda =
\begin{cases}
  \CC[\![Q^{1/2}]\!] & \text{when $\cX = \PP(1,1,2)$} \\
  \CC[\![Q^{1/3}]\!] & \text{when $\cX = \PP(1,1,1,3)$} \\
  \CC[\![Q_1,Q_2]\!] & \text{when $\cX = \FF_2$ or $\cX = \FF_3$}
\end{cases}
\]
$\Lambda$ is called the \emph{Novikov ring}.  For later use, we define
the rings:
\begin{align*}
  \Lambda\{z\} & = 
\begin{cases}
  \CC[z][\![Q^{1/2}]\!] & \text{when $\cX = \PP(1,1,2)$} \\
  \CC[z][\![Q^{1/3}]\!] & \text{when $\cX = \PP(1,1,1,3)$} \\
  \CC[z][\![Q_1,Q_2]\!] & \text{when $\cX = \FF_2$ or $\cX = \FF_3$}
\end{cases} \\
  \Lambda\{z,z^{-1}\} & = 
\begin{cases}
  \CC[z,z^{-1}][\![Q^{1/2}]\!] & \text{when $\cX = \PP(1,1,2)$} \\
  \CC[z,z^{-1}][\![Q^{1/3}]\!] & \text{when $\cX = \PP(1,1,1,3)$} \\
  \CC[z,z^{-1}][\![Q_1,Q_2]\!] & \text{when $\cX = \FF_2$ or $\cX = \FF_3$}
\end{cases}
\end{align*}

The {\em genus-zero descendant potential of $\cX$} is
\[
\cF^0_{\cX}(t_0,t_1,\ldots) = 
\sum_{d \in \EffX} \sum_{n \geq 0}
{Q^d \over n!}  \correlator{\bt(\psi),\ldots,\bt(\psi)}^{\cX}_{0,n,d}
\]
where $t_0, t_1,\ldots$ are orbifold cohomology classes on $\cX$ and
$\bt(\psi) = t_0 + t_1 \psi + t_2 \psi^2 + \ldots$\phantom{.} This is a formal
power series in the co-ordinates $t_i^\alpha$ of $t_i = t_i^1 \phi_1 +
\ldots + t_i^N \phi_N$ with Taylor coefficients given by genus-zero
Gromov--Witten invariants:
\[
\cF^0_{\cX}(t_0,t_1,\ldots) = 
\sum_{\substack{d \in \EffX \\n \geq 0}}
\sum_{\substack{k_1,\ldots,k_n \\ \alpha_1,\ldots,\alpha_n}}
{Q^d t_{k_1}^{\alpha_1} \ldots t_{k_n}^{\alpha_n} \over n!}
\correlator{\phi_{\alpha_1} \psi^{k_1}, \phi_{\alpha_2}
  \psi^{k_2}, \ldots, \phi_{\alpha_n} \psi^{k_n}}^{\cX}_{0,n,d}
\]

The \emph{big quantum orbifold cohomology} of $\cX$ is the family of
$\Lambda$-algebra structures on $\HorbXL$ defined by
\begin{equation}
  \label{eq:bigQC}
  \OPpair{\phi_\alpha \bullet_\tau \phi_\beta}{\phi_\gamma} =
  \Correlator{\phi_\alpha,\phi_\beta,\phi_\gamma}{\tau}.
\end{equation}
This family of products $\bullet_{\tau}$ is parametrized by $\tau$ in
a formal neighbourhood of zero 
\footnote{This just means that the
  right-hand side of \eqref{eq:bigQC} is a formal power series in the
  co-ordinates $\tau^1, \tau^2,\ldots,\tau^N$ of $\tau $.} 
in $\HorbX$.  The \emph{small quantum orbifold cohomology} of $\cX$ is a
related family of algebra structures on $\HorbX$ which will be
described in detail in Section~\ref{sec:smallQC}.  It is defined,
roughly speaking, by restricting the parameter $\tau$ in
$\bullet_\tau$ to lie in $H^2(\cX;\CC) \subset \HorbX$.

\begin{rem}
  If $\cX$ is a manifold then orbifold cohomology, quantum orbifold
  cohomology, the orbifold Poincar\'e pairing, and orbifold
  Gromov--Witten invariants coincide respectively with usual
  cohomology, usual quantum cohomology, the usual Poincar\'e pairing,
  and usual Gromov--Witten invariants.
\end{rem}

\subsection{Variations of Semi-Infinite Hodge Structure} 
\label{sec:VHS}

The key notion in this paper is that of a \emph{variation of
  semi-infinite Hodge structure} or \VHSns.  This was introduced by
Barannikov \cite{Barannikov:periods} as part of his study of
higher-dimensional mirror symmetry.  \VHSns s occur both in the
mathematical version of the A-model (quantum cohomology and
Gromov--Witten theory) and in the mathematical version of the
B-model (singularity theory and Landau--Ginzburg models).  As we
will see, mirror symmetry in this context amounts to the assertion
that the \VHS associated with an A-model is isomorphic to the \VHS
associated with its B-model mirror.

A more traditional formulation of mirror symmetry is as the equality
of certain families of Frobenius algebras: small quantum cohomology on
the A-side and certain Jacobi rings on the B-side.  As mentioned in
the Introduction, one can obtain a family of Frobenius algebras from a
\VHS by choosing an \emph{opposite subspace} and a \emph{dilaton
  shift}. In good cases --- in the \emph{miniversal} situation --- a
\VHSns, an opposite subspace, and a dilaton shift together determine a
Frobenius manifold in the sense of Dubrovin \cite{Dubrovin}; this is
also known as a \emph{flat structure} in the sense of Kyoji Saito
\cite{saitoK}.  The method of constructing flat structures which we
describe here was originally developed by Morihiko Saito \cite{saitoM}
in the context of singularity theory, and was reformulated in terms of
\VHSns s and applied to mirror symmetry by Barannikov
\cite{Barannikov:projective}.

\begin{nt}
  Let $\Czzinv$ denote the ring of Laurent power series which converge
  on $\{z: 0<|z|<\epsilon\}$ for some $\epsilon>0$ which depends on
  the series under consideration.  Let $\Cz$ be the subring of
  $\Czzinv$ consisting of functions regular at $z=0$, and let $\Czinv$
  be the ring of holomorphic functions on $\PP^1\setminus \{0\}$.  We
  have 
  \[
  \Czzinv = \Cz \oplus z^{-1}\Czinv
  \]
  Let $(\cM,\cO_{\cM})$ be a smooth complex analytic space or its
  formal germ.  When $\cM$ is a complex analytic space, we define
  $\cOzzinv{\cM}$ to be the sheaf of relative Laurent series in $z$:
  for an open set $U\subset \cM$, $\cOzzinv{\cM}(U)$ is the set of
  functions which are holomorphic on $\{(q,z)\in U\times \CC :
  0<|z|<\epsilon(q)\}$ for some positive continuous map $\epsilon:U
  \to \RR$ which depends on the function under consideration.  Let
  $\cOz{\cM}$ be the subsheaf of $\cOzzinv{\cM}$ consisting of
  functions regular at $z=0$.  When $(\cM,\cO_{\cM})$ is a formal
  germ, we consider Laurent series in $z$ convergent in an adic
  topology: for a regular parameter system $t^1,\dots, t^l$ on $\cM$
  we set $\cOzzinv{\cM}:=\CC[z,z^{-1}][\![t^1,\dots,t^l]\!]$ and
  $\cOz{\cM}:=\CC[z][\![t^1,\dots,t^l]\!]$.
\end{nt}

\begin{dfn}
  \label{dfn:VHS}
  A \emph{variation of semi-infinite Hodge structure} (\VHSns) with
  base $\cM$ is a locally free $\cOz{\cM}$-module $\cE$ of finite rank
  equipped with a flat $z$-connection\footnote{A $z$-connection is a
    connection multiplied by $z$.}
  \[
  \nabla^z \colon \cE \longrightarrow \Omega^1_{\cM}
  \otimes_{\cO_{\cM}} \cE
  \]
  and a pairing 
  \[
  \Ppair{\cdot}{\cdot}_{\cE} \colon \cE \times \cE \rightarrow
  \cOz{\cM}
  \]
  which satisfy 
  \begin{align*}
    \nabla_X^z\big(f(q,z) s\big) &= \big(z X\, f(q,z)\big)s + f(q,z) \,
    \nabla^z_X s\\
    \big[\nabla_X^z, \nabla_Y^z\big] &= z \nabla_{[X,Y]}^z \\ 
    \Ppair{s_1}{s_2}_{\cE} &= \Ppair{s_2}{s_1}_{\cE}\big |_{z\mapsto -z} \\
    \Ppair{f(q,-z) s_1}{s_2}_{\cE} &= \Ppair{s_1}{f(q,z)s_2}_{\cE} =
    f(q,z) \Ppair{s_1}{s_2}_{\cE} \\
    zX\Ppair{s_1}{s_2}_{\cE} &= -\Ppair{\nabla^z_X s_1}{s_2}_{\cE} +
    \Ppair{s_1}{\nabla^z_X s_2}_{\cE}
  \end{align*}
  for all $f\in \cOz{\cM}$ and all vector fields $X,Y$ on $\cM$.  The
  pairing is assumed to be \emph{non-degenerate} in the sense that the
  induced pairing
  \[
  \big(\cE/z\cE\big) \otimes_{\cO_{\cM}} \big(\cE/z\cE\big) \rightarrow \cO_{\cM}
  \]
  is non-degenerate.  A \emph{grading} on this \VHS is a
  $\CC$-endomorphism $\Gr\colon \cE\rightarrow \cE$ such that there
  exists a vector field $E$ on $\cM$ and a constant $D \in
  \CC$ satisfying
  \begin{equation}
    \label{eq:Grconditions}
    \begin{aligned}
      \Gr\big(f(q,z) s\big) &= \big((2z\partial_z +2E)f(q,z)\big) s + f(q,z) \Gr (s)  \\
      \big[\Gr, \nabla^z_X\big] &= 2\nabla^z_X + \nabla^z_{[2E,X]} \\
      (2z\partial_z + 2E) \Ppair{s_1}{s_2}_{\cE}
      &=\Ppair{\Gr(s_1)}{s_2}_{\cE} + \Ppair{s_1}{\Gr(s_2)}_{\cE} - 2
      D\Ppair{s_1}{s_2}_{\cE}
    \end{aligned}
  \end{equation}
  for all $f \in \cOz{\cM}$ and all vector fields $X$ on $\cM$.  The
  vector field $E$, which is uniquely determined by $\Gr$, is called
  the \emph{Euler vector field}.
\end{dfn}

The analogy with a usual variation of Hodge structure comes from the
family of filtrations $\cdots \supset z^{-1}\cE \supset \cE \supset
z\cE \supset \cdots$ of $\cE\otimes_{\cOz{\cM}} \cOzzinv{\cM}$.  The
existence of the $z$-connection is \emph{Griffiths transversality} for
this family.

\begin{rem} \label{rem:convenience}
  When defining \VHSns s one can choose from many function rings in
  $z$: polynomial functions, entire functions, formal power series,
  $L^2(S^1,\CC)$, etc.  All the \VHSns s in our paper can in fact be
  defined over polynomial functions in $z$.  Also, the A-model \VHS is
  always defined over polynomial functions in $z$.  We chose the ring
  $\Cz$ for technical convenience: it lets us use the Segal--Wilson
  Grassmannian below.
\end{rem} 

Suppose  that the \VHS $\cE$ is generated by one section $s_0$ together
with its derivatives
\[
\nabla^z_{X_1} \nabla^z_{X_2} \cdots \nabla^z_{X_k} s_0 \in \cE.
\] 
In this situation $\cE$ gives rise to a \emph{family of Frobenius
  algebras} over $\cM$. Let $T^*\cM$ be the cotangent bundle of $\cM$
and set
\[
\cO_{T^*\cM}:=\bigoplus_{k=0}^\infty \Sym^k(\cT_\cM)
\]
where $\cT_{\cM}$ is the tangent sheaf of $\cM$.  Then $\cE/z\cE$
becomes an $\cO_{T^*\cM}$-module via the map 
\[
\Sym^k(\cT_{\cM}) \ni X_1X_2 \cdots X_k \longmapsto \big[\nabla^z_{X_1} \nabla^z_{X_2} \cdots
\nabla^z_{X_k} s_0\big] \in \cE/z\cE
\]
or in other words 
\[
\cT_{\cM} \ni X \longmapsto \big[\nabla^z_X \cdot\big]
\in \End(\cE/z\cE).
\]
Our assumption implies that there is an exact sequence
\[
\begin{CD}
  0 @>>> \cI @>>> \cO_{T^*\cM} @>>> \cE/z\cE @>>> 0 
\end{CD} 
\]
where $\cI$ is an ideal sheaf.  This identification of $\cE/z\cE$ with
$\cO_{T^*\cM}/\cI$ gives $\cE/z\cE$ an algebra structure.  The
pairing $\Ppair{\cdot}{\cdot}_{\cE}$ induces a pairing $\cE/z\cE
\otimes_{\cO_{\cM}} \cE/z\cE \rightarrow \cO_{\cM}$ and makes
$\cE/z\cE$ into a Frobenius algebra.  Note that the ideal sheaf $\cI$
is independent of the choice of generator $s_0$, so even though
different choices of $s_0$ produce different Frobenius algebra
structures on $\cE/z\cE$ they are all isomorphic as algebras.

\begin{dfn} \label{dfn:miniversal}
  A \VHS is said to be \emph{miniversal} if there is a section $s_0$ of
  $\cE$ such that the $\cO_{\cM}$-module map 
  \begin{equation}
    \label{eq:miniversality}
    \cT_{\cM} \ni X \longmapsto \big[\nabla^z_X s_0\big] \in \cE/z\cE
  \end{equation}
   is an isomorphism. This is equivalent to the natural composition
  \[
  \cT_{\cM} \hookrightarrow \cO_{T^*\cM} \rightarrow \cO_{T^*\cM}
  /\cI
  \]
  being an isomorphism.
\end{dfn}

In the miniversal case each tangent space $T_q\cM$ is naturally
equipped with a ring structure, independent of the choice of $s_0$
above.  As we will see below, if we can choose a ``good'' opposite
subspace then this product structure arises from a Frobenius manifold;
this implies that $\cM$ has the structure of an
$F$-manifold\footnote{See Hertling--Manin \cite{Hertling--Manin} or
  Manin \cite{Manin} for $F$-manifolds and Coates--Ruan
  \cite{Coates--Ruan}*{Section~6.2} or 
  Iritani \cite[Section~3.2]{Iritani:RIMS} for
  an expanded version of this remark.} without any choice.

\subsubsection{A Moving Subspace Realization}
\label{sec:moving}
As was indicated in the Introduction, the A-model \VHS arises as a
family of subspaces moving in a fixed symplectic vector space $\cHX$.
We now explain how to give such a ``moving subspace'' realization of
any \VHSns.

Consider the universal cover $\pi\colon \tcM\rightarrow \cM$ and let
$\cH$ denote the space of flat sections of
$\pi^\star\big(\cE\otimes_{\cOz{\cM}}\cOzzinv{\cM}\big)$:
\begin{equation*}
\cH = \left\{s\in \Gamma\(\tcM, \pi^\star\(\cE\otimes_{\cOz{\cM}}\cOzzinv{\cM}\)\) 
: \nabla^z s =0\right\} 
\end{equation*}
$\cH$ is a free $\Czzinv$-module of the same rank as $\cE$.  For each $q\in
\tcM$, we can embed the fiber $(\pi^\star\cE)_q$ into $\cH$ via the
map $\iota_q$, where
\[
\iota_q \colon (\pi^\star \cE)_q \ni v \longmapsto s \in \cH 
\text{ such that } s(q)= v.
\]
Let $\E_q\subset \cH$ denote the image of this embedding.  Because the
pairing $\Ppair{\cdot}{\cdot}_{\cE}$ on $\cE$ is $\nabla^z$-flat,
$\Ppair{s_1}{s_2}_{\cE}$ is a constant as an element of $\Czzinv$ for
any $s_1,s_2\in \cH$.  Define a symplectic form $\Omega$ on $\cH$ by
\[
\Omega(s_1,s_2) := \Res_{z=0} \Ppair{s_1}{s_2}_{\cE} dz. 
\]
In the graded case, the grading operator $\Gr$ induces a
$\CC$-endomorphism $\Gr$ of $\cH$ satisfying
\[
\Gr \big(f(z) s\big) = \big(2z\partial_z f(z)\big) s + f(z) \Gr(s) 
\]
for all $f\in \Czzinv$ and all $s\in \cH$.

We construct our moving subspace realization first in the case where
$\cM$ is a complex analytic space and the \VHS is analytic.  Take a
point $q_0 \in \cM$ and fix a basis for $\E_{q_0}/z\E_{q_0}$ over
$\CC$.  Pick a $\Cz$-basis $e_1,\dots e_N$ for $\E_{q_0}$ such that
$[e_1]$,\ldots,$[e_N]$ is our chosen basis for $\E_{q_0}/z\E_{q_0}$
and choose local sections $s_i$ of $\cE$ such that $s_i(q_0) =
e_i(q_0)$.  Then for $q$ near $q_0$ the subspace $\E_q$ is spanned
over $\Cz$ by vectors $\iota_q(s_i)$ and we can write $\iota_q(s_i) =
\sum_{j=1}^N s_{ij}(q,z) e_j$.  Let $S(q,z)$ be the matrix with
$(i,j)$ entry equal to $s_{ij}(q,z)$.  There exists $\epsilon>0$ such
that each matrix entry $s_{ij}(q,z)$ converges on $0<|z|<2\epsilon$,
and $S$ therefore defines a $C^\infty$ loop $S^1 \ni z \mapsto
S(q,\epsilon z)\in GL(N,\CC)$.  This loop depends on the choice of
local sections $s_1$, \ldots,$s_N$.  Removing this choice we obtain,
after fixing an isomorphism $\E_{q_0}/z\E_{q_0} \cong \CC^N$, a
holomorphic map from a neighbourhood of $q_0$ to the \emph{
  Segal--Wilson Grassmannian} $LGL(N,\CC)/L^+GL(N,\CC)$
\cite{Pressley--Segal}; the choice of isomorphism here is the choice
$[e_1]$,\ldots,$[e_N]$ of basis for $\E_{q_0}/z\E_{q_0}$.  At least
locally, therefore, the assignment $q\mapsto \E_q$ gives an analytic
family of elements of the Segal--Wilson Grassmannian.  When $\cM$ is a
formal germ and the \VHS is formal we proceed in the same way,
obtaining instead an $\infty$-jet in the Segal--Wilson Grassmannian.
With these interpretations in place, we have:

\begin{pro}
  A \VHS with base $\cM$ gives a family of subspaces $\E_q$ in $\cH$
  parameterized by $\cM$ such that
  \begin{itemize}
  \item $\E_q$ is a free $\Cz$-module of dimension equal to the rank
    of $\cE$
  \item $X \E_q \subset z^{-1} \E_q$ for any tangent vector $X\in T_q\cM$
  \item $\E_q$ is maximal isotropic with respect to the symplectic
    form $\Omega$
\end{itemize} 
In the graded case we have $(2E+\Gr)\E_q \subset \E_q$, so that in
particular $\Gr\E_q \subset z^{-1}\E_q$.
\end{pro} 
\begin{proof}
The first three properties are obvious from the definition. 
Consider the graded case and take a local section 
$s_i(q)=\sum_{j=1}^N s_{ij}(q,z)e_j(q)$
of $\cE$ as above, where $e_j(q)$ is $\nabla^z$-flat:
$e_j(q)=\iota_q^{-1}(e_j)$ for $e_j\in \cH$.  
Then:
\begin{align*}
\iota_q(\Gr(s_i)) &= \sum_{j=1}^N 
\iota_q\left(\big((2E+2z\partial_z)s_{ij}(q,z)\big) e_j(q) + s_{ij}(q,z)\Gr(e_j(q))\right) \\
&=\sum_{j=1}^N \left(2E (s_{ij}(q,z)e_j)+ \Gr(s_{ij}(q,z)e_j)\right)\\
&=(2E+\Gr)\iota_q(s_i)
\end{align*}
The left hand side here is an element of $\E_q$.
\end{proof}

\subsubsection{Opposite Subspaces and Frobenius Manifolds}
\begin{dfn}
  Given a \VHS with base $\cM$, an \emph{opposite subspace} at $q\in
  \tcM$ is a free $\Czinv$-submodule $\cH_-$ of $\cH$ such that the
  natural map
  \[
  \cH_- \oplus \E_q \rightarrow \cH
  \]
  is an isomorphism.  This implies that the projections
  \[
  z\cH_-/\cH_- \longleftarrow z\cH_- \cap \E_q \longrightarrow \E_q/z\E_q
  \]
  are isomorphisms.
\end{dfn}

Being opposite to $\cH_-$ is an open condition on $\cM$.  We can see
this using the geometry of the Segal--Wilson Grassmannian.  By
choosing an opposite subspace $\cH_-$ at $q_0$ and an isomorphism
$z\cH_-/\cH_-\cong \CC^N$ we can identify each subspace $\E_q$ with a
point in $LGL(N,\CC)/L^+GL(N,\CC)$.  Then $\cH_-$ is opposite at $q$
if and only if the point in the Grassmannian corresponding to $\E_q$
is in the \emph{big cell}, and the big cell is an open orbit of
$L^-GL(N,\CC)$.

An opposite subspace $\cH_-$ at $q_0\in \tcM$ also defines a
trivialization of $\pi^\star \cE$ near $q_0$:
\begin{equation}
\label{eq:trivialization} 
(\pi^\star \cE)_q \cong \E_q \cong (\E_q\cap z\cH_-)\otimes \Cz
\cong (z\cH_-/\cH_-) \otimes \Cz  
\end{equation}
for $q$ in some open neighbourhood of $q_0 \in \tcM$.  We call
this the \emph{flat trivialization} associated to $\cH_-$.

\begin{pro}
  \label{pro:normalformofconn}
  Let $\cH_-$ be an opposite subspace at $q_0\in \tcM$.  Under the 
  flat trivialization \eqref{eq:trivialization} associated to $\cH_-$,
  the flat $z$-connection $\nabla^z$ becomes
  \begin{equation}
    \label{eq:normalformofconn}
    \nabla^z_X = z \, d_X + A_X 
  \end{equation}
  where $d$ is the connection defined by the flat trivialization and
  $A$ is a $z$-independent $\End(z\cH_-/\cH_-)$-valued 1-form on
  $\tcM$.  If moreover $\cH_-$ is isotropic with respect to $\Omega$
  then the pairing $\Ppair{\cdot}{\cdot}_{\cE}$ induces and can be
  recovered from the symmetric $\CC$-bilinear pairing:
  \begin{equation}
    \label{eq:normalformofpairing}
    \begin{split}
      \Abspair{\cdot}{\cdot} \colon \big(z\cH_-/\cH_-\big) \otimes
      \big(z\cH_-/\cH_-\big) & \longrightarrow \CC \\
      [f_1] \otimes [f_2] & \longmapsto \Ppair{f_1}{f_2}_{\cE}\Big|_{z=\infty}
    \end{split}
  \end{equation}
  This pairing \eqref{eq:normalformofpairing} is non-degenerate and satisfies:
  \begin{equation}
    \label{eq:Frobeniusalgebra}
    \Abspair{A_X v}{w} = \Abspair{v}{A_X w}
  \end{equation} 
  If in addition the \VHS is graded and $\cH_-$ is preserved by the
  operator $\Gr$ then $\Gr$ takes the form
  \begin{equation}
    \Gr= \label{eq:normalformofgrading}
    2z\partial_z + 2d_{E} + \Gr_0 
  \end{equation}
  where $\Gr_0\in \End\big(z\cH_-/\cH_-\big)$ is a constant operator induced
  by the action of $\Gr$ on $z\cH_-\subset \cH$.  This $\Gr_0$ satisfies
  \begin{equation}
    \label{eq:integrationanddimension}
    2D\Abspair{v}{w}= \Abspair{\Gr_0(v)}{w} + \Abspair{v}{\Gr_0(w)}
  \end{equation}
  where $D$ is the constant from Definition~\ref{dfn:VHS}.
\end{pro}
\begin{proof}
  Take a section $s$ of $\cE$ which corresponds to a constant vector
  in $z\cH_-/\cH_-$ under the trivialization
  \eqref{eq:trivialization}.  Then $\iota_q(s(q)) = v_0 + w(q)$ for
  some constant vector $v_0\in z\cH_-$ and some $w(q)\in \cH_-$.
  Thus $\iota_q(\nabla_X^z s(q))= zX \iota_q(s(q)) = zXw(q)$ is in
  $\E_q\cap z\cH_-$.  This proves \eqref{eq:normalformofconn}.

  If $\cH_-$ is isotropic with respect to $\Omega$ then it is clear
  that the restriction of $\Ppair{\cdot}{\cdot}_{\cE}$ to $z\cH_-\cap
  \E_q$ takes values in $\CC$.  Under the identification $z\cH_-\cap
  \E_q \cong z\cH_-/\cH_-$, this coincides with the pairing
  \eqref{eq:normalformofpairing}.  Non-degeneracy follows from the
  non-degeneracy of the pairing $\Ppair{\cdot}{\cdot}_{\cE}$, and
  equation \eqref{eq:Frobeniusalgebra} follows from the
  $\nabla^z$-flatness of $\Ppair{\cdot}{\cdot}_{\cE}$.

  If $\cH_-$ is preserved by $\Gr$ then so is $z\cH_-$, and thus $\Gr$
  induces a constant operator $\Gr_0\in \End(z\cH_-/\cH_-)$; equations
  \eqref{eq:normalformofgrading} and
  \eqref{eq:integrationanddimension} follow immediately from
  \eqref{eq:Grconditions}.
\end{proof}

In the miniversal and graded case, the structures in Proposition
\ref{pro:normalformofconn} assemble to give a Frobenius manifold with
base an open subset of $\tcM$.  Take an $\Omega$-isotropic,
$\Gr$-invariant opposite subspace $\cH_-$.  Assume that there exists
an eigenvector $[v_0]\in z\cH_-/\cH_-$ of $\Gr_0$ such that the
corresponding section $s_0\in \cE$ under the trivialization
\eqref{eq:trivialization} makes \eqref{eq:miniversality} into an
isomorphism.  The eigenvector $[v_0]$ is called the \emph{dilaton shift}.
The affine subspace $v_0+\cH_-$ meets each $\E_q$ in a unique point
$\iota_q(s_0)$.  Barannikov's \emph{period map} is 
\begin{equation}
  \label{eq:periodmap}
  \begin{split}
    \Psi\colon \tcM &\longrightarrow \cH  \\
    q &\longmapsto \iota_q\big(s_0(q)\big).  
  \end{split}
\end{equation}
The miniversality condition implies that the map 
\begin{equation}
\label{eq:flatcoordinate}
\begin{split}
  \tcM & \longrightarrow z\cH_-/\cH_- \\ 
  q & \longmapsto \big[z(\Psi(q)-v_0)\big]
\end{split}
\end{equation}
is a local isomorphism.  The linear co-ordinates on $\tcM$ given by
this map are called \emph{flat co-ordinates}.  By differentiating
\eqref{eq:flatcoordinate} we obtain the \emph{Kodaira--Spencer map}
\[
\KS\colon \cT_{\tcM}\ni X\longmapsto  A_X [v_0] \in z\cH_-/\cH_-
\]
Pulling back the metric $\Abspair{\cdot}{\cdot}$ on $z\cH_-/\cH_-$
along $\KS$ gives a non-degenerate symmetric bilinear pairing:
\[
g\colon \cT_{\tcM}\otimes \cT_{\tcM} \rightarrow \cO_{\tcM}
\]
This is constant with respect to the flat co-ordinates, so the ``flat
co-ordinates'' really are flat co-ordinates for the
metric\footnote{Note that the metric $g$ is a $\CC$-bilinear form not a
  Hermitian form.} $g$.  The product $\bullet$ on tangent spaces is
defined by
\[
A_{X\bullet Y} [v_0] = A_X A_Y [v_0].
\]
The identity vector field $e$ is given by 
\[
A_e [v_0] = [v_0].
\]
It is easy to check that these data satisfy all the axioms for a 
Frobenius manifold:
\begin{pro}
  \label{pro:VHStoFrobeniusmanifold}
  Take an $\Omega$-isotropic, $\Gr$-invariant opposite subspace
  $\cH_-$ and a dilaton shift $[v_0]\in z\cH_-/\cH_-$ such that
  the corresponding section $s_0\in \cE $ under the trivialization
  \eqref{eq:trivialization} makes \eqref{eq:miniversality} into an
  isomorphism.  Then the data $(\bullet, e, g, E)$ defined above
  determine a Frobenius manifold with base an open subset of $\tcM$.
  In other words:
  \begin{enumerate}
  \item the Levi--Civita connection $\nabla^{\rm LC}$ of $g$ is flat
  \item $\big(T_q \cM, \bullet, g\big)$ is an associative, commutative
    Frobenius algebra
  \item the pencil of connections $\nabla^{\lambda}_X = \nabla^{\rm
      LC}_X +\lambda X\bullet$ is flat
  \item the identity vector field $e$ is flat
  \item the Euler vector field satisfies $\big(\nabla^{\rm LC}\big)^2 E=0$ and
    \begin{align*}
      Eg(X,Y) &= g([E,X],Y) + g(X, [E,Y]) + (\alpha+ 2 - D) g(X,Y), \\
      [E, X\bullet Y] & = [E, X] \bullet Y + X \bullet [E, Y] + [X,Y] 
    \end{align*} 
    where $\alpha$ is the eigenvalue of $[v_0]$ with respect to $\Gr_0$
  \end{enumerate}
  In these terms, $\Gr_0$ is given by 
  \[
  \Gr_0\big(\KS(X)\big) = \KS\big((2+\alpha)X - 2 \nabla^{\rm LC}_X E\big).
  \]
\end{pro} 

\begin{rem}
  Even in the non-miniversal case, if the map
  \eqref{eq:flatcoordinate} is injective with image an affine subspace
  of $z\cH_-/\cH_-$ then we still refer to the linear co-ordinates on 
  $\tcM$ given by \eqref{eq:flatcoordinate} as \emph{flat
    co-ordinates}.  In this case the pullback of the pairing on
  $z\cH_-/\cH_-$ via the Kodaira--Spencer map will not in general give
  a metric on the base $\tcM$.  But these ``flat co-ordinates'' and the
  constant pairing on $z\cH_-/\cH_-$ are what would remain from the
  flat structure on the Frobenius manifold if our non-miniversal \VHS
  arose as a ``slice'' of a miniversal \VHSns.  This is exactly the
  relationship between small and big quantum cohomology.
\end{rem}

\subsection{The Big A-Model \VHS and Givental's Symplectic Formalism}
\label{sec:Giventalapproach}

We now define a \VHS which gives rise, through an appropriate choice
of opposite subspace and dilaton shift, to the Frobenius manifold
structure on big quantum orbifold cohomology.  This is the \emph{big
  A-model \VHSns} described in the Introduction.  As we will see
below, the moving subspace realization of this \VHS gives the family
of tangent spaces to Givental's Lagrangian cone.

\subsubsection{The Big A-Model \VHSns} 
\label{sec:bigAVHS}

Let $H$ be a formal neighbourhood of zero in $H^\bullet_{\rm
  orb}(\cX,\Lambda)$, where $\Lambda$ and associated rings are defined
in Section~\ref{sec:GWdefs}.  Recall that, for a linear co-ordinate
system $\tau^1,\ldots,\tau^N$ on $H^\bullet_{\rm orb}(\cX,\CC)$,
we have:
\begin{align*}
  &\cO_H = \Lambda[\![\tau^1, \ldots,\tau^N]\!] \\
  &\cOz{H} = \Lambda\{z\}[\![\tau^1, \ldots,\tau^N]\!] \\
  &\cOzzinv{H} = \Lambda\{z,z^{-1}\}[\![\tau^1, \ldots,\tau^N]\!]
\end{align*}
The big A-model \VHS has base $H$ and is given by
\[
\cE^{\rm big}_{\rm A} := H^\bullet_{\rm orb}(\cX,\CC) \otimes \cOz{H}
\]
The flat $z$-connection is the \emph{Dubrovin connection}
\[
\nabla^z = z\, d + \sum_{\alpha=1}^N \(\phi_\alpha \bullet_\tau\) d \tau^\alpha
\] 
the pairing is
\begin{align*}
\cE^{\rm big}_{\rm A} \times \cE^{\rm big}_{\rm A} 
&\longrightarrow \cOz{H}  \\
f(z)\times g(z) &\longmapsto \OPpair{f(-z)}{g(z)}
\end{align*}
and the grading operator $\Gr\colon \cE_{\rm A}^{\rm big}\rightarrow
\cE_{\rm A}^{\rm big}$ and Euler vector field $E$ are 
\begin{align*}
\Gr &= 2z\partial_z+  2 d_E + \Gr_0 \\ 
E &=\sum_{\alpha} \(1- \frac{\deg \phi_\alpha}{2}\) 
\tau^\alpha \parfrac{}{\tau^\alpha} +  
\sum_{\alpha : \deg \phi_\alpha=2} \rho_\alpha \parfrac{}{\tau^\alpha} 
\end{align*}
where $\Gr_0$ is the usual grading operator on orbifold cohomology
\[
\Gr_0(\phi_\alpha)= \deg (\phi_\alpha) \, \phi_\alpha
\]
and $c_1(\cX)= \sum_{\alpha:\deg \phi_{\alpha}=2} \rho_\alpha \phi_\alpha$.
These data satisfy the axioms for a graded \VHS (Definition
\ref{dfn:VHS}) with $D=\dim_{\CC}\cX$, except that the ground ring
$\CC$ there is replaced here by $\Lambda$.

\subsubsection{Givental's Symplectic Formalism} 

Following Givental \cite{Givental:symplectic}, we now describe the
genus-zero Gromov--Witten theory of $\cX$ in terms of a Lagrangian
submanifold $\cLX$ of the symplectic vector space
\begin{align*}
  \cHX &:= \HorbX \otimes \Lambda\{z,z^{-1}\} && 
  \text{(the vector space)}
  \\
  \Omega(f,g) &:= \Res_{z=0} \OPpair{f(-z)}{g(z)}  dz && \text{(the
    symplectic form)}.
\end{align*}
Relations between genus-zero Gromov--Witten invariants --- the String
Equation, Dilaton Equation, and Topological Recursion Relations ---
translate into very strong constraints on the geometry of $\cLX$: see
Coates--Givental \cite{Coates--Givental:QRRLS}*{Proposition~1},
Givental \cite{Givental:symplectic}, and Coates--Corti--Iritani--Tseng
\cite{CCIT:computing}*{Appendix~B}.  These constraints can be
rephrased as the statement that the tangent spaces to $\cLX$ form the
moving subspace realization of a \VHSns; we will see in the next
section that this is the big A-model \VHSns.

The space $\cHX$ is the direct sum of Lagrangian subspaces
\begin{align*}
  \cHX^+ = \HorbX \otimes \Lambda\{z\}, && 
  \cHX^- = z^{-1} \HorbX \otimes \Lambda\{z^{-1}\}
\end{align*}
and this polarization identifies the symplectic space $\(\cHX,\Omega\)$
with the cotangent bundle $T^\star \cHX^+$.  We will regard the
genus-zero descendant potential $\cF^0_{\cX}$ as the formal germ of
a function on $\cHX^+$.  A general point of $\cHX^+$ takes the form
\[
\bq(z) = q_0 + q_1 z + q_2 z^2 + \ldots
\]
where $q_0, q_1,\ldots$ are orbifold cohomology classes on $\cX$.
Setting
\begin{equation}
  \label{eq:dilatonshift}
  q_k =
  \begin{cases}
    t_k & k \ne 1 \\
    t_1 - 1 & k=1
  \end{cases}
\end{equation}
makes $\cF_{\cX}^0$ into a function on a formal
neighbourhood\footnote{This just means that $\cF^0_{\cX}$ is a formal
power series in the variables $t^\alpha_k$.} 
of the point $-z$ in $\cHX^+$. 
This change of variables $\bq(z) = \bt(z) - z$ is
called the {\em dilaton shift}.

The Lagrangian submanifold $\cLX$ is the graph of the differential of
$\cF^0_{\cX}$.  Since $\cF^0_{\cX}$ is defined only in a formal
neighbourhood  of $-z \in \cHX^+$,
$\cLX$ is a formal germ of a Lagrangian submanifold of $T^\star
\cHX^+$.  The polarization $\cHX = \cHX^+ \oplus \cHX^-$ 
identifies $T^\star \cHX^+$ with $\cHX$, and we regard 
$\cLX$ as a formal germ of a
Lagrangian submanifold of $\cHX$ via this identification.  
$\cLX$ has a more concrete description as follows.  
A general point of $\cHX$ has the form 
\[
\sum_{k=0}^{\infty} \sum_{\mu = 1}^{N} q^\mu_k \, \phi_\mu z^k +
\sum_{l=0}^\infty \sum_{\nu=1}^N p_{l,\nu} \, \phi^\nu (-z)^{-1-l}
\]
where $\phi^1,\ldots,\phi^N$ is the basis for $\HorbX$ such that
$\OPpair{\phi_i}{\phi^j} = \delta_i^{\phantom{i}j}$, and this defines
Darboux co-ordinates $\{q^\mu_k, p_{l,\nu}\}$ on $\cHX$.  In these
co-ordinates $\cLX$ is given by
\[
p_{l,\nu} = {\partial \cF^0_{\cX} \over \partial q^\nu_l}
\]
so a general point on $\cLX$ takes the form:
\begin{equation}
  \label{eq:pointonL}
-z + \bt(z) + 
\sum_{\substack{d \in \EffX \\n \geq 0}} \sum_{\alpha = 1}^N
{Q^d \over n!}  \correlator{\bt(\psi),\ldots,\bt(\psi), {\phi_{\alpha}
  \over -z-\psi}}^{\cX}_{0,n+1,d} \phi^\alpha
\end{equation}
The expression ${1 \over -z-\psi}$ here should be expanded as a power
series in $z^{-1}$.  Note that $\cLX$ encodes all genus-zero
Gromov--Witten invariants of $\cX$.

Turning off gravitational descendants, by setting
\[
t_k =
\begin{cases}
  \tau & k=0 \\
  0 & k \ne 0
\end{cases}
\]
gives a family $\tau \mapsto J(\tau,-z)$ of elements of $\cLX$
parametrized by $H$.  We call this the {\em big $J$-function} of
$\cX$.   From \eqref{eq:pointonL} we have:
\begin{equation}
  \label{eq:Jfunction}
  J_{\cX}(\tau,-z) = -z + \tau + \sum_{\alpha = 1}^N 
  \BigCorrelator{\phi^\alpha \over
    -z-\psi}{\tau} \, \phi_\alpha
\end{equation}

As discussed, $\cLX$ has a very special geometric structure.  
Precise geometric properties we need 
are given in Proposition \ref{pro:TisLtau} below. 
These imply that $\cLX$ is the germ of a Lagrangian 
cone with vertex at the origin such that 
each tangent space $T$ is tangent to the cone 
exactly along $zT$. 
Define a submersion $\maptau\colon \cLX \to H$ by  
\[
\maptau(\bq,\bp) = 
\sum_{d\in \EffX}
\sum_{n\ge 0}
\sum_{\alpha=1}^N 
\frac{Q^d}{n!}
\correlator{1,\bt(\psi),\bt(\psi),\dots, \bt(\psi),\phi^\alpha}_{0,n+2,d}^\cX 
\phi_\alpha.  
\] 
Here we used the dilaton shift $\bq(z) = \bt(z) -z $ as before.  The
String Equation shows that $\maptau\big(J_\cX(\tau,-z)\big) = \tau$,
i.e. that the $J$-function is a section of the map
$\maptau\colon \cLX \to H$.  Denote by $T_\tau$ the tangent space to $\cL_\cX$ at
$J_{\cX}(\tau,-z)$.

\begin{pro} \label{pro:TisLtau} \ 
  \begin{itemize} 
  \item[(a)] The tangent space to $\cLX$ 
    at $(\bq,\bp)$ coincides with $T_\tau$
    for $\tau = \maptau(\bq,\bp)$. 
  \item[(b)] The tangent space $T_\tau$ is closed under 
    multiplication by $z$ and has the structure of a 
    $\Lambda\{z\}$-module. Moreover, it is freely 
    generated over $\Lambda\{z\}$ by the partial derivatives 
    \[
    {\partial \over \partial \tau^1} J_{\cX}(\tau,-z), \ldots,
    {\partial \over \partial \tau^N} J_{\cX}(\tau,-z). 
    \]
  \item[(c)] 
  The fiber at $\tau\in H$ 
  of the map $\maptau\colon \cLX\to H$ 
  is given by 
  \[
  z T_\tau \cap \big(\cHX, -z\big) 
  \]
  where $\big(\cHX,-z\big)$ is the formal neighborhood of $-z$ in
  $\cHX$; see \cite{CCIT:computing}*{Appendix B} for the definition of
  $\big(\cHX,-z\big)$ as a formal scheme over $\Lambda$.  The
  intersection here should be interpreted as the set of $R$-valued
  points for an arbitrary $\Lambda$-algebra $R$ when $\tau$ is an
  $R$-valued point of $H$.
  \end{itemize}
\end{pro}

\begin{proof}
  Part (a) is \cite{CCIT:computing}*{Lemma~B6} and part (b) is
  \cite{CCIT:computing}*{Proposition~B3}.  We sketch the proof of (c)
  following Givental \cite{Givental:symplectic}.  Take $\bof \in
  \maptau^{-1}(\tau)$.  The String Equation says that $z^{-1} \bof \in
  T_\bof \cLX$.  Thus $\bof \in z T_\bof \cLX = z T_\tau$ by (a), and
  so $\maptau^{-1}(\tau) \subset zT_\tau$.  In particular $\Ker
  d_\bof\maptau \subset zT_\tau$.  Since $\maptau$ is submersive, both
  $\Ker d_\bof \maptau$ and $zT_\tau$ are codimension-$(\dim H)$
  subspaces of $T_\tau$, so $\Ker d_\bof \maptau = z T_\tau$.  Part
  (c) follows.
\end{proof}

Part (c) of this proposition shows that 
$\cLX$ can be reconstructed from its tangent spaces as 
\begin{equation}
\label{eq:T_reconstructs_L}
\cLX = \bigcup_{\tau\in H} 
\Big(z T_\tau \cap \big(\cHX, -z\big)\Big). 
\end{equation} 

\subsubsection{The Big A-Model \VHS and Givental's Cone}
\label{sec:bigQCandcone}

Consider the map:
\begin{align*}
\JJ\colon \cE_{\rm A}^{\rm big} & 
\longrightarrow H^{\bullet}_{\rm orb}(\cX,\Lambda)  
\otimes \cOzzinv{H} \\  
\phi_\alpha & \longmapsto \parfrac{}{\tau^\alpha} J_{\cX}(\tau,z) 
\end{align*}
Using the well-known relationship between the $J$-function and the big
quantum product
\[
z\parfrac{}{\tau^\alpha} \parfrac{}{\tau^\beta} J_{\cX}(\tau,z) 
= \sum_{\gamma=1}^N \Correlator{\phi_\alpha,\phi_\beta,\phi^\gamma}{\tau} 
\parfrac{}{\tau^\gamma} J_{\cX}(\tau,z)
\]
we see that $\JJ$ satisfies
\begin{equation}
\label{eq:diffeq_JJ}
\JJ\circ \nabla^z_X = z \,d_X\circ  \JJ 
\end{equation} 
for any vector field $X$ on $H$.  Thus $\JJ$ sends any
$\nabla^z$-parallel section to a constant section: it identifies the
space of $\nabla^z$-flat sections of the big A-model \VHS --- which
we denoted by $\cH$ in Section~\ref{sec:moving} --- with Givental's
symplectic space $\cHX$.  Proposition~\ref{pro:TisLtau}(b) implies
that the image of each fiber $\cE_{\rm A, \tau}^{\rm big}$ under $\JJ$ 
coincides, if we flip the sign of $z$, 
with the tangent space $T_{\tau}$ to
$\cL_{\cX}$: 
\begin{align}
\label{eq:EisT}
\E_\tau(\Lambda) := \JJ(\cE_{\rm A, \tau}^{\rm big}) = T_\tau |_{z\mapsto -z}. 
\end{align} 
Thus the moving subspace realization $\E_\tau(\Lambda)$ 
of the big A-model \VHS determines
the Lagrangian submanifold $\cLX$ via \eqref{eq:T_reconstructs_L}. 
Since $\JJ \circ \Gr = \big(2z\partial_z + 2 d_E + \Gr_0 - 2c_1(\cX)/z\big)
\circ \JJ$, we see that the grading operator acts on $\cHX$ by:
\begin{equation}
\label{eq:gradingonH}
\Gr|_{\cHX} = 2 z\partial_z + \Gr_0 - 2c_1(\cX)/z.
\end{equation}

The standard opposite subspace for the big A-model \VHS is $\cHX^-$.
This is clearly isotropic with respect to $\Omega$; it is also
preserved by $\Gr$.  The period map \eqref{eq:periodmap} associated to
this opposite subspace and the choice $[v_0] = 1 \in z\cHX^-/\cHX^-$
(the dilaton shift) is $\tau \mapsto z^{-1} J_\cX(\tau,z)$:
\begin{equation}
\label{eq:J_intersection}
\E_\tau(\Lambda)\cap (1 + \cHX^-) = \left\{z^{-1} J_\cX(\tau, z)\right\}. 
\end{equation} 
Since $J(\tau, z) = z + \tau + O(z^{-1})$, the flat co-ordinates here
are the usual linear co-ordinates $\tau^i$ on $H$.  The big A-model
\VHS is miniversal, and these choices of opposite subspace $\cHX^-$
and dilaton shift $[v_0]$ produce the usual big quantum orbifold
cohomology Frobenius manifold defined by Chen--Ruan
\cite{Chen--Ruan:GW}.

\subsubsection{The Analytic Big A-Model \VHS and Givental's Cone 
Over $\C$}
\label{subsubsec:analyticbig} 
Suppose that the big quantum product $\bullet_\tau$ 
in \eqref{eq:bigQC} is convergent as a power series in $\tau$ and $Q$. 
Then we can consider the big A-model \VHS 
and the Lagrangian submanifold over $\C$ 
by specializing all the Novikov variables $Q_i$ to 1 as follows. 
By the Divisor Equation (see e.g.
Abramovich--Graber--Vistoli \cite[Theorem~8.3.1]{AGV:2}) 
and equation \eqref{eq:bigQC}, 
the big quantum product becomes 
\begin{equation}
\label{eq:bigproduct_divisor}
\big(\phi_\alpha\bullet_\tau \phi_\beta,\phi_\gamma\big) 
= \sum_{d\in \EffX} 
\sum_{n\ge 0} 
\frac{Q^d e^{\CHpair{\sigma}{d}}}{n!} 
\correlator{\phi_\alpha,\phi_\beta,\phi_\gamma,
\tau',\dots,\tau'}_{0,d,n+3}^\cX
\end{equation} 
where $\tau = \sigma + \tau'$ is the decomposition of $\tau$ into the
non-twisted second cohomology class $\sigma$ and the sum $\tau'$ of
other components.  Our convergence assumption therefore implies that
the specialization $\bullet_\tau|_{Q_i=1}$ is analytic on a domain
$U \subset \HorbX$ of the form:
\[
U = \Big\{\tau = (\sigma, \tau') \;:\; \text{ $\textstyle\Re(\int_d \sigma) < -M$
  for all $d \in \EffX\setminus\{0\}$, $\|\tau'\| < \epsilon$}\Big\}  
\]
for a sufficiently big $M>0$ and a sufficiently small $\epsilon>0$.
Note that $\bullet_\tau$ for sufficiently small values of $Q$ equals
$\bullet_{\hat\tau}|_{Q_i=1}$ with $\hat\tau = \tau + \sum_i p_i \log
Q_i$.  In particular, the origin $\tau=Q=0$ of $H$ corresponds to the
limit direction $\Re(\CHpair{\sigma}{d})\to -\infty$,~$\tau'\to 0$ in
$U$.  This is called the \emph{large radius limit}.  The
specialization $\cE_{\rm A}^{\rm big}|_{Q_i=1}$ is defined as an
analytic \VHS on the base $U$ in the same way as in Section
\ref{sec:bigAVHS}.  Because the map $\JJ$ is a solution to the
differential equation \eqref{eq:diffeq_JJ}, the specialization
\[
\JJ|_{Q_i=1} \colon \cE_{\rm A}^{\rm big}|_{Q_i=1}
\longrightarrow \HorbX \otimes \cOzzinv{U} 
\]
is well-defined on $U$ and gives an analytic family of moving
subspaces:
\begin{equation}
\label{eq:moving_A}
\E_\tau := \JJ|_{Q_i=1}\big(\cE_{{\rm A},\tau}^{\rm big}|_{Q_i=1}\big) 
\subset \cHX|_{Q_i=1} := \HorbX\otimes \Czzinv
\end{equation} 
In view of the relations \eqref{eq:T_reconstructs_L} and
\eqref{eq:EisT}, we \emph{define} the analytic version of Givental's
cone to be the set 
\begin{equation}
\label{eq:analyticLag}
\cLX^{\rm an} := \bigcup_{\tau \in U} 
z \E_\tau|_{z\mapsto -z}. 
\end{equation} 
Note that the analytic version is no longer a germ at some point.  In
fact if we work with an $L^2$ version of the symplectic formalism,
replacing $\cHX|_{Q_i=1}$ with $\HorbX\otimes L^2(S^1,\C)$, then in a
neighbourhood of the section $\tau \mapsto J_\cX(\tau,-z)|_{Q_i=1}$,
$\cLX^{\rm an}$ has the structure of a Hilbert submanifold of
$\HorbX\otimes L^2(S^1,\C)$ and this submanifold is the graph of the
differential of the analytic function $\cF_\cX^0|_{Q_i=1}$.  When the
quantum product $\bullet_\tau|_{Q_i=1}$ is analytically continued to a
bigger domain, the analytic Givental's cone $\cLX^{\rm an}$ can be
enlarged to a bigger submanifold using the construction above.  This
is the \emph{analytic continuation of Givental's cone} appearing in
Conjecture \ref{conj}.

\subsection{Small Quantum Cohomology}
\label{sec:smallQC}

Small quantum orbifold cohomology is a family of algebra structures on
$\HorbX$ defined, roughly speaking, by restricting the parameter
$\tau$ of the big quantum product $\bullet_\tau$ to lie in
$H^2(\cX;\CC) \subset \HorbX$.  In this section we make this precise,
and also explain how small quantum orbifold cohomology arises from a
non-miniversal \VHS --- the \emph{small A-model \VHSns}.

\subsubsection{Small Quantum Orbifold Cohomology}
\label{sec:smallQCsub}

From the equation \eqref{eq:bigproduct_divisor}, 
the big quantum product $\bullet_{\tau}$ at a
non-twisted second cohomology class $\tau \in H^2(\cX;\CC)$ becomes:
\[
\OPpair{\alpha\bullet_{\tau} \beta}{\gamma} = \sum_{d\in \EffX} Q^d 
e^{\CHpair{\tau}{d}} \correlator{\alpha,\beta,\gamma}_{0,3,d}^{\cX}
\]
This shows that the Novikov parameters keep track of the modes of
Fourier expansion in $\tau\in H^2(\cX,\CC)$ and that the product
$\bullet_{\tau}$ depends on $\tau$ only through the \emph{exponentiated
  K\"{a}hler parameters} $e^{\tau^i}$.  The \emph{small quantum
  orbifold cohomology} of $\cX$ is the possibly multi-valued family of
algebras $(H^\bullet_{\rm orb}(\cX;\CC), \circ_q)$ defined by
\begin{equation}
  \label{eq:smallQC}
  \OPpair{\alpha\circ_{q} \beta}{\gamma} = \sum_{d\in \EffX} q^d 
  \correlator{\alpha,\beta,\gamma}_{0,3,d}^{\cX} 
\end{equation}
where the parameter $q$ lies on the torus
\[
\Upsilon = H^2(\cX;\CC)/ 2\pi\tti H^2(\cX;\ZZ) 
\]
and $q^d$ denotes the following possibly multi-valued function on
$\Upsilon$:
\begin{align*}
  q^d \colon \Upsilon \ni [\tau] \longmapsto
  \exp\big(\textstyle\CHpair{\tau}{d}\big)\in \Cstar && \tau\in H^2(\cX;\CC)
\end{align*}
The cohomology groups here denote sheaf cohomology \emph{of the stack
  $\cX$} and not of the coarse moduli space.

In order to make geometric sense of \eqref{eq:smallQC}, we introduce
co-ordinates on the torus $\Upsilon$ and consider the associated
partial compactification.  Let $p_1,\dots, p_r$ be an integral basis
of the free part of $H^2(\cX,\ZZ)$ such that $p_i$ evaluates
non-negatively on $\EffX$.  We define $\Cstar$-valued co-ordinates
$q_i$ on $\Upsilon$ by
\[
q_i\colon \Upsilon\ni 
\big[\tau^1 p_1+ \cdots + \tau^r p_r\big] \longmapsto e^{\tau^i}\in \Cstar
\]
and for $d\in \EffX$ we write:
\[
q^d = q_1^{\CHpair{p_1}{d}} \cdots q_r^{\CHpair{p_r}{d}}
\] 
Each exponent $\CHpair{p_i}{d}$ here is a non-negative rational
number.  The co-ordinates $(q_1,\dots,q_r)$ give a partial
compactification of $\Upsilon$, $\Upsilon \hookrightarrow \CC^r$, and
the small quantum product $\circ_q$ defines a possibly multi-valued
family of associative algebras in a formal neighbourhood of the origin
in $\CC^r$.  The origin of $\CC^r$ is called the \emph{large radius
  limit point} of $\Upsilon$.

\begin{rem} 
\label{rem:orbifoldpoint}
When pulled back to a suitable finite cover $\widetilde{\CC^r}$ of $\CC^r$, 
the product $\circ_q$ becomes single-valued. 
The multi-valuedness of the product $\circ_q$ happens only 
for orbifolds and introduces an orbifold singularity at 
the large radius limit point. 
\end{rem} 

In our examples we already chose suitable integral bases $\{p_i\}$ for
$H^2(\cX,\ZZ)$ in Section \ref{sec:cohomology}.  This gives
co-ordinates on $\Upsilon$, which we denote by $q_1,q_2$ if $\cX=\FF_2$
or $\FF_3$ and by $q$ if $\cX=\PP(1,1,2)$ or $\PP(1,1,1,3)$. In many
cases, including the examples in our paper, the small quantum product
$\circ_q$ is known to be convergent in a neighbourhood of the large
radius limit point.  In what follows we will assume this, writing
$U_{\rm A}\subset \CC^r$ for the domain of convergence of $\circ_q$
and $\cM_{\rm A}=U_{\rm A}\cap \Upsilon$.  $\cM_{\rm A}$ is called the
\emph{K\"{a}hler moduli space} or \emph{A-model moduli space}.

\begin{rem}
  \label{rem:qandQ}
  The Novikov variables $Q$ and $Q_1$, $Q_2$ are not the same as the
  parameters $q$ and $q_1$, $q_2$ for small quantum cohomology.  But
  the restriction of the big quantum product $\bullet_\tau$ to the
  locus $\tau \in H^2(\cX;\CC)$ can be recovered from
  the small quantum product $\circ_q$ by setting
   \begin{align*}
     &q = Q e^t &\text{where } \tau & = t \, p \\
     &          &\text{and }  \cX  & = \text{$\PP(1,1,2)$ or
       $\PP(1,1,1,3)$} \\
     \intertext{or}
     &q_1  = Q_1 e^{\tau^1}, q_2 = Q_2 e^{\tau^2}, &\text{where } \tau & =
     \tau^1 p_1 + \tau^2 p_2 \\
     &          &\text{and }  \cX & = \text{$\FF_2$ or $\FF_3$}.
   \end{align*}
\end{rem}

\subsubsection{The Small A-Model \VHSns} 
\label{sec:smallAVHS}

The \emph{small A-model \VHSns} has base $\cM_{\rm A}$.  It is the
free $\cOz{\cM_{\rm A}}$-module
\[
\cE_{\rm A}= H^\bullet_{\rm orb}(\cX;\CC)\otimes \cOz{\cM_{\rm A}}
\] 
with flat $z$-connection given by
\[
\nabla^z = z \, d + \sum_{i=1}^r (p_i\circ_q)  \frac{d q_i}{q_i}
\]
pairing given by
\begin{align*}
\cE_{\rm A} \times \cE_{\rm A} & 
\longrightarrow \cOz{\cM_{\rm A}} \\ 
f(z)\times g(z) & \longmapsto \OPpair{f(-z)}{g(z)} 
\end{align*}
and grading operator $\Gr\colon \cE_{\rm A}\rightarrow \cE_{\rm A}$
and Euler vector field $E$ given by
\begin{align*}
\Gr &= 2 z\partial_z + 2 d_E + \Gr_0\\
E &= \sum_{i=1}^r \rho_i q_i \parfrac{}{q_i}  
\end{align*}
where $\Gr_0$ is the usual grading operator on orbifold cohomology and
$c_1(X) = \sum_{i} \rho_i p_i$.  The Dubrovin connection $\nabla^z$
here is independent of our choice of co-ordinates on $\cM_{\rm A}$; it
extends to a connection on $U_{\rm A}$ with a logarithmic singularity
along the normal crossing divisor $q_1q_2\cdots q_r=0$ (Deligne's
extension).  

The flat $z$-connection $\nabla^z$ makes $\cE_{\rm A}$ into a $D$-module in
the sense of Givental \cite{Givental:homological}.  When we want to
emphasize this structure, we will refer to $\cE_{\rm A}$ as the
\emph{quantum $D$-module}.  See Guest \cite{Guest} and Iritani
\cite{Iritani:gen} for more on this.

One obtains the small A-model \VHS from the big A-model \VHS by
restricting $\tau$ to lie in $H^2(\cX,\CC)$ and
specializing Novikov variables $Q_i$ to $1$.  In the
moving subspace realization, the small \VHS therefore corresponds a
subfamily of tangent spaces to Givental's Lagrangian cone 
$\cL_{\cX}^{\rm an}$,
and to the following subcone of $\cL_{\cX}^{\rm an}$: 
\begin{equation}
\label{eq:subcone}
\bigcup_{\tau\in H^2(\cX;\CC)} z \E^\cX_\tau|_{z\mapsto -z}.   
\end{equation} 
This ``small subcone'' has a standard slice, the \emph{small $J$-function},
which is obtained from the big $J$-function by restricting $\tau$ to
lie in $H^2(\cX;\CC)$ and then setting the Novikov variables to $1$.
In our examples we find, by applying the Divisor Equation to
\eqref{eq:Jfunction}, that the small $J$-functions of $\PP(1,1,1,3)$,
$\FF_3$, $\PP(1,1,2)$ and $\FF_2$ are:
\begin{align}
 \label{eq:smallJP112}
 &J_{\PP(1,1,2)}(q,z)  = z \, q^{p/z} \( \fun_0 + \sum_{\substack{d \colon 2
      d \in \ZZ \\ d > 0}} \sum_{\alpha = 1}^N q^d \correlator{\phi^\alpha\over z(z -
    \psi)}^{\PP(1,1,2)}_{0,1,d} \,
  \phi_\alpha\) \displaybreak[0] \\
  \label{eq:smallJF2}
  &J_{\FF_2}(q_1,q_2,z) = z \, q_1^{\: p_1/z} q_2^{\: p_2/z} \( 1 +
  \sum_{k,l \geq 0} \sum_{\alpha = 1}^N q_1^{\: k} q_2^{\: l}
  \correlator{\phi^\alpha\over z(z - \psi)}^{\FF_2}_{0,1,(k,l)} \,
  \phi_\alpha\)  \displaybreak[0]\\
  \label{eq:smallJP1113}
  &J_{\PP(1,1,1,3)}(q,z)  = z \, q^{p/z} \( \fun_0 + \sum_{\substack{d \colon 3
      d \in \ZZ \\ d > 0}}  \sum_{\alpha = 1}^N q^d \correlator{\phi^\alpha\over z(z -
    \psi)}^{\PP(1,1,1,3)}_{0,1,d} \,
  \phi_\alpha\) \displaybreak[0]\\
  \label{eq:smallJF3}
  &J_{\FF_3}(q_1,q_2,z) = z \, q_1^{\: p_1/z} q_2^{\: p_2/z} \( 1 +
  \sum_{k,l \geq 0} \sum_{\alpha = 1}^N q_1^{\: k} q_2^{\: l}
  \correlator{\phi^\alpha\over z(z - \psi)}^{\FF_3}_{0,1,(k,l)} \,
  \phi_\alpha\)
\end{align}
These are multi-valued analytic functions $\cM_{\rm A} \to \HorbX \otimes
\Czzinv$, for the appropriate choice of target space $\cX$.

The small $J$-function corresponds to the unit section of $\cE_{\rm A}$, 
i.e. $z^{-1} J_{\cX}(q,z) = \iota_q (\fun)$. 
Since in the cases at hand the small
quantum orbifold cohomology algebra of $\cX$ is generated by
$H^2(\cX;\CC)$, the small A-model \VHS is generated by this section
together with its derivatives.  As discussed above 
Definition~\ref{dfn:miniversal}, this makes $\cE_{\rm A}/z\cE_{\rm A}$ into a
Frobenius algebra: the algebra structure here is the small quantum
orbifold cohomology of $\cX$ and the pairing is the orbifold
Poincar\'e pairing.

\subsection{Mirror Symmetry}
\label{sec:mirrorsym}

We now define the \emph{B-model \VHS} discussed in the Introduction,
and explain what we mean by mirror symmetry.  

\subsubsection{The B-Model \VHSns} 
A \emph{Landau--Ginzburg model} in this context is a holomorphic family
$\pi\colon Z\rightarrow {\cM}_{\rm B}$ of affine Calabi--Yau manifolds ---
for us they will be algebraic tori --- together with a function
$W\colon Z\rightarrow \CC$ called the \emph{superpotential} and a
section $\omega$ of the relative canonical sheaf $K_{Z/\cM_{\rm B}}$ which
gives a holomorphic volume form $\omega_y$ on each fiber
$Z_y=\pi^{-1}(y)$. The base space $\cM_{\rm B}$ of the family is called the
\emph{B-model moduli space}. Landau--Ginzburg models which
correspond under mirror symmetry to the quantum cohomology of toric
varieties have been constructed by Givental \cite{Givental:ICM,
  Givental:toric} and Hori--Vafa \cite{Hori--Vafa}. In this section we
explain how to obtain a \VHS --- the B-model \VHS --- from a
Landau--Ginzburg model. 

\begin{assumptions}
  \label{assumptions}
  We can assume (by deleting any points at which this condition fails
  to hold) that for each $y_0 \in \cM_{\rm B}$ there exists a neighbourhood
  $U$ of $y_0$ and a constant $M>0$ such that for all $y \in U$, all
  the critical points of $W_y$ are contained in the set $\{x\in Z_y :
  |W_y(x)|<M\}$.  We further assume that:
  \begin{itemize}
  \item[(a)] the family of pairs
    \begin{align*}
      \big(Z_y,N_{y,\theta} \big) && \text{where} 
      &&&  N_{y,\theta} = \left\{x : \Re\big(W_y(x) e^{\tti\theta}\big)\le -M\right\} \\
      &&&&&\big(y,e^{\tti\theta}\big)\in U \times S^1
    \end{align*}
    is topologically locally trivial on $U\times S^1$
  \item[(b)] all the data $\pi:Z\rightarrow \cM_{\rm B}$, $W$, $\omega$
    are algebraic
  \item[(c)] there is a complete K\"{a}hler metric on $Z_y$ such that
    the set 
    \[
    \big\{ x\in Z_y: \| \mathop{\mathrm{grad}} W_y (x) \| <C \big\}
    \]
    is compact for some $C>0$
  \item[(d)] the critical points of $W_y$ are isolated and generically
    non-degenerate
  \end{itemize}
  These assumptions are satisfied by the Landau--Ginzburg mirrors to
  $\PP(1,1,2)$, $\FF_2$, $\PP(1,1,1,3)$, and $\FF_3$; furthermore 
  (a--c) here have been proved for the mirror to a general compact
  toric orbifold by Iritani \cite{Iritani:integral}*{Section~4.2}.
  We will use condition (a) when constructing a local system of
  relative homology groups.  Condition (b) is much stronger than we
  need: we use it only to ensure the convergence of certain integrals,
  and this certainly follows from a polynomial-growth condition on the
  integrand \eqref{eq:diffform}.  The remaining conditions allow us to
  use Morse theory.  Assumption (c) implies that we can choose a
  metric without introducing critical points ``at infinity'': it holds
  for the mirrors to toric varieties.  In the examples at hand, the
  critical points of $W_y$ are always distinct and non-degenerate: we
  denote them by $\sigma_1,\dots,\sigma_N$.
\end{assumptions}

Under our assumptions, a Landau--Ginzburg model determines a local
system $R^\vee$ on $\cM_{\rm B}\times \Cstar$ with fiber over $(y,z)$
equal to the relative homology group
\[
R^\vee_{(y,z)}= H_n\big(Z_y, \{x \in Z_y : \Re (W_y(x)/z)\ll 0\}\big).  
\]
Let $\cO_{\cM_{\rm B}\times \Cstar}$ denote the analytic structure sheaf.
The associated locally free sheaf $\cR^\vee=R^\vee\otimes
\cO_{\cM_{\rm B}\times \Cstar}$ has a Gauss--Manin connection, which is
flat.  We construct flat sections of $\cR^\vee$ using Morse theory,
defining the cycle $\Gamma_k(y,z)$, $k\in \{1,\ldots,N\}$, to be the
closure of the union of downward gradient flowlines for the function
$x\mapsto \Re\big(W_y(x)/z\big)$ from the critical point $\sigma_k$ of
$W_y$.  If the imaginary parts of the critical values of $W_y/z$ are
all distinct then the image of $\Gamma_k(y,z)$ under $W_y/z$ becomes a
negative half-line from $W_y(\sigma_k)/z$ parallel to the real axis.
The cycles $\Gamma_k(y,z)$ form a basis for the relative homology
group $R_{(y,z)}^\vee$.  Note that these sections become multi-valued
under analytic continuation: $\Gamma_k(y,z)$ has monodromy in both $y$
and $z$.

The dual bundle $\cR$ on $\cM_{\rm B}\times \Cstar$ has fiber over
$(y,z)$ equal to the relative cohomology group $H^n\big(Z_y, \{x \in Z_y :
\Re(W_y(x)/z)\ll 0\}\big)$.  This bundle also has a flat Gauss--Manin
connection, as well as a distinguished section
\begin{equation}
\label{eq:generatorofEB}
(y,z) \longmapsto \exp\big(W_y/z\big) \, \omega_y .
\end{equation}
The co-ordinates of this section with respect to the flat frame dual
to $\Gamma_1(y,z)$, \ldots,$\Gamma_N(y,z)$ are \emph{oscillating
  integrals}:
\begin{equation}
\label{eq:oscint}
\cI_{k}(y,z)=
\int_{\Gamma_k(y,z)} \exp\big(W_y/z\big) \, \omega_y
\end{equation}
We will consider only sections of $\cR$ represented by differential
forms
\begin{align}
\label{eq:diffform}
f(x,z) \exp\big(W(x)/z\big) \, \omega_y, && f(x,z) \in 
(\pi\times \id)_\star\cO_{Z\times \Cstar} 
\end{align}
such that $x \mapsto f(x,z)$ is algebraic on each fiber $Z_y$; note that
integrals of such forms over cycles $\Gamma_k(y,z)$ are convergent.

\begin{dfn}
\label{dfn:mirrorDmodule}
Given a Landau--Ginzburg model $\big(\pi\co Z \to \cM_{\rm B},W,\omega\big)$, we
define the associated \emph{B-model \VHS} as follows.  The base of
the B-model \VHS is $\cM_{\rm B}$.  Let $\cE_{\rm B}$ be the
$\cOz{\cM_{\rm B}}$-module consisting of sections of $\cR$ of the form
\eqref{eq:diffform} with $f$ regular in a small neighbourhood of
$z=0$.  Let $\nabla^z$ be the flat $z$-connection on $\cE_{\rm B}$ given by
\[
\nabla^z_X = z\nabla^{\rm GM}_X 
\]
where $X$ on the left-hand side is a vector field on $\cM_{\rm B}$,
$\nabla^{\rm GM}$ is the Gauss--Manin connection on $\cR$, and $X$ on
the right-hand side denotes the standard lift to a vector field on
$\cM_{\rm B} \times \Cstar$.  The pairing on $\cE_{\rm B}$ is defined as the dual
to the intersection pairing on relative homology groups:
\[
R_{(y,-z)}^\vee \otimes R_{(y,z)}^\vee \rightarrow \CC
\]
\end{dfn}

As in Section~\ref{sec:smallAVHS}, the flat $z$-connection $\nabla^z$
makes $\cE_{\rm B}$ into a $D$-module in the sense of Givental
\cite{Givental:homological}.  When we want to emphasize this
structure, we will refer to $\cE_{\rm B}$ as the \emph{mirror $D$-module}.

The pairing on $\cE_{\rm B}$ is given, at generic $y$, by
\[
\Ppair{[s_1]}{[s_2]}_{\cE_{\rm B}} = 
\frac{1}{\big(2\pi\tti z\big)^{n}}
\sum_{k=1}^N 
\int_{\Gamma_k(y,-z)} s_1(-z) \cdot \int_{\Gamma_k(y,z)} s_2(z)
\]
because $\{\Gamma_k(y,-z)\}_k$ and $\{\Gamma_k(y,z)\}_k$ are
mutually-dual bases for the relative homology groups.  We now check
that the data in Definition~\ref{dfn:mirrorDmodule} satisfy the axioms
for a \VHSns.

\begin{lem}
The pairing $\Ppair{\cdot}{\cdot}_{\cE_{\rm B}}$ takes values in 
$\cOz{\cM_{\rm B}}$. 
Also, $\cE_{\rm B}$ is a free $\cOz{\cM_{\rm B}}$-module of 
the same rank as $\cR$. 
\end{lem}
\begin{proof}
  The method of stationary phase gives the asymptotic expansion of the
  oscillating integrals:
  \begin{equation*}
    \int_{\Gamma_i(y,z)} f(x,z) \exp\big(W(x)/z\big) \, \omega_y 
    \sim \frac{\big({-2\pi z}\big)^{n/2}}
    {\sqrt{\Hess W_y(\sigma_i)}}
    e^{W_y(\sigma_i)/z} \Big(f(\sigma_i,0) + O(z)\Big)  
  \end{equation*}
  for $f$ regular at $z=0$.  Here $z$ goes to zero in an angular
  sector where the order of $\Im(W_y(\sigma_1)/z),\dots,
  \Im(W_y(\sigma_N)/z)$ is unchanged.  The Hessian of $W_y$ at
  $\sigma_i$ is calculated in terms of local co-ordinates
  $x^1,\dots,x^n$ on $Z_y$ such that $\omega_y = dx^1\wedge \cdots
  \wedge dx^n$ near $\sigma_i$.  At generic $y$, therefore, we have
  \[
  \Ppair{\big[f(x,z)e^{W/z}\,\omega_y\big]}{\big[g(x,z)e^{W/z}\,\omega_y\big]}_{\cE_{\rm B}}
  \sim \sum_{k=1}^N 
  \frac{f(\sigma_k,0) g(\sigma_k,0)}{\Hess W_y(\sigma_k)} + O(z) 
  \]
  The first term on the right hand side here is the \emph{residue
    pairing} of the elements $f(x,0)$ and $g(x,0)$ of the Jacobi ring
  of $W_y$.  The left hand side is holomorphic on $0<|z|<\epsilon$ for
  some $\epsilon$; the above asymptotics imply that it is actually
  regular at $z=0$.

  In a neighbourhood of each point $y\in \cM_{\rm B}$, we can find
  fiberwise-algebraic functions $\phi_k(x)\in \pi_\star \cO_Z$,
  $k=1,\dots, N$, such that $[\phi_k(x)]$ forms a basis of the Jacobi
  ring $J(W_y)$.  If $y$ is generic then we can choose $\phi_k$ such
  that $\phi_i(\sigma_j) = \delta_{ij}$.  Let
  $s_i=\big[\phi_i(x)e^{W/z}\,\omega_y\big]$ be the corresponding
  section of $\cE_{\rm B}$.  Then the Gram matrix
  $\Ppair{s_i}{s_j}_{\cE_{\rm B}}$ is of the form
  $\pairW{\phi_i}{\phi_j}+O(z)$ with the first term given by the
  residue pairing.  The non-degeneracy of the matrix
  $\pairW{\phi_i}{\phi_j}$ implies that $s_1,\dots,s_N$ form an
  $\cOz{\cM_{\rm B}}$-basis of $\cE_{\rm B}$.
\end{proof} 

In our examples it turns out that $\cE_{\rm B}$ is generated by the single
section \eqref{eq:generatorofEB} together with its derivatives.  In
other words
\begin{equation}
\label{eq:Dmodulegeneratedbyoscint}
\cE_{\rm B} \cong \cOz{\cM_{\rm B}}\langle z \cT_{\cM_{\rm B}}\rangle 
/\{P(y,z\partial,z):
 \text{$P \cI_{k}(y,z)=0$ for all $k$}\}
\end{equation}
and so $\cE_{\rm B}$ is generated as a $D$-module by oscillating integrals.
The family of Frobenius algebras $\cE_{\rm B}/z\cE_{\rm B}$ determined by $\cE_{\rm B}$
and the section \eqref{eq:generatorofEB} consists of the family of
Jacobi rings $\bigcup_{y \in \cM_{\rm B}} J(W_y)$ equipped with the residue pairing.

\begin{rem}
  In the Introduction we described the B-model \VHS as a family of
  subspaces $\E_q$, $q \in \cM_{\rm B}$, in the fibers of a vector bundle $V
  \to \cM_{\rm B}$ with flat connection.  The vector bundle $V$ is $\cE_{\rm
    B} \otimes_{\cOz{\cM_{\rm B}}} \cOzzinv{\cM_{\rm B}}$, the family of subspaces is the
  subbundle $\cE_{\rm B}$ of $V$, and the connection on $V$ is $\nabla^{\rm GM}$.
\end{rem}
\subsubsection{The Mirror Conjecture and Mirror Theorems} 
\label{sec:VHSmirrorsym}
We now formulate a mathematical version of mirror symmetry, in the
spirit of Givental's ICM lecture \cite{Givental:ICM}.

\begin{conj}
\label{conj:mirrorsym}
Let $\cE_{\rm A}$ be the small A-model \VHS of $\cX$; recall that
this has base the A-model moduli space $\cM_{\rm A}$.  Let $\cE_{\rm B}$
be the B-model \VHS associated to the Landau--Ginzburg mirror to
$\cX$; this has base the B-model moduli space $\cM_{\rm B}$.  Let
$\Delta_\epsilon^*\subset \CC$ be a punctured disc of radius
$\epsilon$.  There is an open set $U_{\rm B}\subset \cM_{\rm B}$ with co-ordinates
$(y_1,\dots,y_r)\colon U_{\rm B} \cong \({\Delta_\epsilon^*}\)^r$ and a
 map $\mir \colon U_{\rm B}  \to \cM_{\rm A}$ of the form
\begin{align*}
   (y_1,\dots, y_r)&\longmapsto (q_1,\dots,q_r)&&\text{where} &&
  q_i=y_i\exp\big(f_i(y_1,\ldots,y_r)\big) \\
   &&&&&
  f_i(0,\ldots,0)=0
\end{align*}
such that there is an isomorphism of \VHSns s:
\begin{equation}
  \label{eq:mirroriso}
  \bigg(\cE_{\rm B},\nabla^z,\Ppair{\cdot}{\cdot}_{\cE_{\rm B}}\bigg) \cong 
  \mir^\star\bigg(\cE_{\rm A}, \nabla^z, \Ppair{\cdot}{\cdot}_{\cE_{\rm A}}\bigg)
\end{equation}
In the graded case this isomorphism also preserves the operator $\Gr$.
\end{conj}

The map $\mir$ here is called the mirror map.  We have seen that the
\VHSns s $\cE_{\rm A}$ and $\cE_{\rm B}$ give rise to families of Frobenius
algebras: $\cE_{\rm A}/z \cE_{\rm A}$ gives the small quantum orbifold cohomology
algebra of $\cX$ equipped with the orbifold Poincar\'e pairing, and
$\cE_{\rm B}/z\cE_{\rm B}$ gives the family of Jacobi rings $\bigcup_{y \in \cM_{\rm B}}
J(W_y)$ equipped with the residue pairing.  Thus
Conjecture~\ref{conj:mirrorsym} implies that there is a
grading-preserving linear isomorphism
\[
\Mir_y:J(W_y) \to \HorbX
\]
which matches the product on the Jacobi
ring $J(W_y)$ with the small quantum product $\circ_{\mir(y)}$ and
matches the orbifold Poincar\'{e} pairing with the residue pairing.

\begin{proof}[How to Prove Conjecture~\ref{conj:mirrorsym}]
  In the examples that we consider below,
  Conjecture~\ref{conj:mirrorsym} is simply a reformulation of mirror
  theorems proved by Givental \cite{Givental:toric} (or equivalently
  by Lian--Liu--Yau \cite{Lian--Liu--Yau:toric}) and by
  Coates--Corti--Lee--Tseng \cite{CCLT}. But this reformulation is
  essential to our argument in Sections~\ref{sec:P1113} and
  \ref{sec:P112}: it allows us to give a systematic
  construction of flat structures near the cusps of the B-model
  moduli space $\cM_{\rm B}$, and to compare the flat structures associated
  to different cusps.  We proceed as follows.

  The oscillating integrals \eqref{eq:oscint} form a basis of
  solutions\footnote{Recall that a \emph{solution} to a $D$-module is
    a solution to the system of differential equations defining that
    $D$-module.}  to the mirror $D$-module $\cE_{\rm B}$.  These solutions
  are multi-valued in both $y$ and $z$.  The system of differential
  equations in $y$ satisfied by the oscillating integrals is called
  the \emph{Picard--Fuchs system} --- see equations
  \eqref{eq:PFofF_3}, \eqref{eq:PFofP1113}, \eqref{eq:PFF2},
  \eqref{eq:PFP112}.  Another basis of solutions to the Picard--Fuchs
  system is given by the components of a cohomology-valued
  hypergeometric function called the $I$-function $I_\cX(y,z)$ --- see
  equations \eqref{eq:I-functionofF3}, \eqref{eq:I-functionofP1113},
  \eqref{eq:IF2}, \eqref{eq:IP112}.  Thus the mirror
  $D$-module $\cE_{\rm B}$ is isomorphic to the $D$-module
  generated by the $I$-function. The $I$-functions are
  multi-valued in $y$ but single-valued in $z$. 

  On the other hand, in our examples the small quantum orbifold
  cohomology algebra is generated by $H^2(\cX;\CC)$ and so the
  A-model \VHS $\cE_{\rm A}$ is generated by the small $J$-function
  $J_\cX(q,z)$ together with its derivatives. This implies that $\cE_{\rm A}$
  is generated as a $D$-module by the small $J$-function. The small
  $J$-functions are multi-valued in $q$ but single-valued in $z$ ---
  see equations \eqref{eq:smallJP112}, \eqref{eq:smallJF2},
  \eqref{eq:smallJP1113}, \eqref{eq:smallJF3}.  

  A Givental-style mirror theorem states that \emph{the $I$-function
    and the small $J$-function coincide after a suitable change of
    variables $y\mapsto q(y)$}:
  \[
  J_{\cX}(q(y),z) = I_{\cX}(y,z)
  \]
  The change of variables $y\mapsto q(y)$ here gives the mirror map in
  Conjecture \ref{conj:mirrorsym}.  Such a mirror theorem implies
  Conjecture \ref{conj:mirrorsym}, as we can then define the
  isomorphism \eqref{eq:mirroriso} to be the $D$-module isomorphism
  which maps:
  \begin{align*}
    (\cE_{\rm B})_y \ni \big[\exp(W_y/z)\,\omega_y\big] && \text{to} &&
    z^{-1} I_{\cX}(y,z)=z^{-1} J_{\cX}(q,z)
    \in \iota_q((\cE_{\rm A})_q)
  \end{align*}
  The matching of gradings and pairings under this isomorphism will be
  explained in Propositions~\ref{pro:matchingofgrading} and
  \ref{pro:matchingofpairing} below. 
\end{proof}

\section{Example: $\FF_3$ and $\PP(1,1,1,3)$}
\label{sec:P1113}

We now apply our general theory to the cases $\cX = \PP(1,1,1,3)$ and
$Y=\FF_3$.  Following the prescriptions of Givental
\cite{Givental:toric} and Hori--Vafa \cite{Hori--Vafa}, we write down
Landau--Ginzburg models which correspond under mirror symmetry to
$\FF_3$ and to $\PP(1,1,1,3)$.  Let $\cM_{\FF_3}$ denote the base of
the Landau--Ginzburg mirror to $\FF_3$ and let $\cM_{\PP(1,1,1,3)}$
denote the base of the Landau--Ginzburg mirror to $\PP(1,1,1,3)$.  We
construct the B-model moduli space $\cM_{\rm B}$, described in the
Introduction, as a partial compactification of $\cM_{\FF_3}$.  The
space $\cM_{\rm B}$ consists of a copy of $\cM_{\FF_3}$ together with a copy
of $\cM_{\PP(1,1,1,3)}$ as a ``divisor at infinity''.  We form a
Landau--Ginzburg model with base $\cM_{\rm B}$ by patching together the
mirrors to $\FF_3$ and $\PP(1,1,1,3)$, and define the B-model
\VHS to be the \VHS associated to this Landau--Ginzburg model.

The space $\cM_{\rm B}$ has two cusps, one of which corresponds to $\FF_3$
and the other to $\PP(1,1,1,3)$.  We show that there is an opposite
subspace canonically associated to each cusp
(Theorem~\ref{thm:uniqueness}) and that these opposite subspaces give
rise, in the sense described in the Introduction, to the small quantum
cohomology of $\FF_3$ and the small quantum orbifold cohomology of
$\PP(1,1,1,3)$.  We show that these two opposite subspaces \emph{do
  not agree} with each other under parallel transport in $\cE_{\rm B}$ ---
this amounts to computing the analytic continuation of the
$I$-function $I_Y(y,z)$ mentioned in Section~\ref{sec:VHSmirrorsym}
--- and from this we prove Theorem~\ref{thm:P1113} and Conjecture~\ref{conj}.

\subsection{The Landau--Ginzburg Mirror to $\FF_3$}

The Landau--Ginzburg mirror of $\FF_3$ is a family of algebraic tori
$\pi\colon Z\to \cM_{\FF_3}$ together with a superpotential
$W\colon Z\to \CC$ and a holomorphic volume form on each fiber
of $\pi$.  Recall that ${\FF_3}$ is defined as a GIT quotient of
$\CC^5$ by $(\Cstar)^2$ where $(\Cstar)^2$ acts via the inclusion:
\begin{align*}
(\Cstar)^2 \hookrightarrow (\Cstar)^5 &&
(s,t) \mapsto (s,s,s,s^{-3}t,t)
\end{align*}
The mirror family $\pi:Z \to \cM_{\FF_3}$ is given by restricting the
dual of this inclusion
\begin{align*}
\pi\colon (\Cstar)^5 & \longrightarrow (\Cstar)^2 \\
(w_1,\dots,w_5) &\longmapsto (w_1w_2w_3w_4^{-3},w_4w_5) 
\end{align*}
to the open subset $\cM_{\FF_3} \subset (\Cstar)^2$ defined by:
\[
\cM_{\FF_3} = \Big\{(y_1,y_2) \in (\Cstar)^2 : \textstyle y_1 \neq -{1
  \over 27} \Big\}
\]  
The superpotential $W$ is
\[
W=w_1+w_2+w_3+w_4+w_5
\]
and the holomorphic volume form $\omega_y$ on the fiber
$Z_y=\pi^{-1}(y_1,y_2)$ is:
\[
\omega_y =\frac{d\log w_1\wedge \cdots \wedge d\log w_5}
{d\log y_1\wedge d\log y_2}
\] 
We deleted the locus $y_1 = {-{1 \over 27}}$ from $\cM_{\FF_3}$ to
ensure that Assumptions~\ref{assumptions} hold.  It is
straightforward to show that the oscillating integrals
\eqref{eq:oscint} satisfy the Picard--Fuchs equations:
\begin{equation}
  \label{eq:PFofF_3}
  \begin{split}
    &D_2 (D_2 - 3 D_1) \,f = y_2 \, f\\
    &D_1^3 D_2 (D_2 - z) (D_2 - 2 z) \, f= y_1 y_2^3 \, f \\
    &D_1^3 D_2 (D_2 - z) \, f = y_1 y_2^2 (D_2 - 3 D_1) \, f \\
    &D_1^3 D_2 \, f = y_1 y_2 (D_2 - 3 D_1) (D_2 - 3 D_1 - z) \, f \\
    &D_1^3 \, f = y_1 (D_2 - 3 D_1)(D_2 - 3 D_1 - z)(D_2 - 3 D_1 - 2z) \, f
  \end{split}
\end{equation}
where $D_1 = z y_1 {\partial \over \partial y_1}$ and $D_2 = z y_2
{\partial \over \partial y_2}$.  

\subsection{Mirror Symmetry for $\FF_3$}
\label{subsec:MS_F3}

We apply Givental's mirror theorem for toric varieties
\cite{Givental:toric}*{Theorem~0.1}.  The $I$\nobreakdash-function
\begin{equation}
\label{eq:I-functionofF3}
I_{\FF_3}(y_1,y_2,z) = z  \sum_{k,l\ge 0} 
\frac{y_1^{k + p_1/z}y_2^{l+p_2/z}\prod_{m=-\infty}^{0} (p_2-3p_1+m z)}
{\prod_{m=1}^{k}(p_1+m z)^3 \prod_{m=1}^{l}(p_2+m z)
\prod_{m=-\infty}^{l-3k}(p_2-3p_1+m z)}
\end{equation}
where $p_1$, $p_2$ is the basis of $H^2(\FF_3)$ defined in Section
\ref{sec:cohomology}, also satisfies the Picard--Fuchs system
\eqref{eq:PFofF_3}.  It coincides with the small $J$-function
\eqref{eq:smallJF3} after a change of variables $(y_1,y_2) \mapsto
(q_1,q_2)$:
\[
J_{\FF_3}(q_1,q_2,z) = I_{\FF_3}(y_1,y_2,z)
\]
As we explained in Section~\ref{sec:VHSmirrorsym}, this proves
Conjecture~\ref{conj:mirrorsym} for $\FF_3$.

Since $J_{\FF_3}(q_1,q_2,z) = z + p_1 \log q_1 + p_2 \log q_2 +
O(z^{-1})$, we can can read off the mirror map $(y_1,y_2) \mapsto
(q_1,q_2)$ by expanding the $I$-function as a Laurent series in
$z^{-1}$.  This gives:
\begin{equation*}
  \begin{split}
    q_1 &= y_1 \exp\(3 \sum_{k\ge 1} (-1)^k {(3k-1)! \over \(k!\)^3} y_1^k\) \\
    q_2 &= y_2 \exp\(\sum_{k\ge 1} (-1)^{k+1} {(3k-1)! \over \(k!\)^3} y_1^k\) 
  \end{split}
\end{equation*}
and hence:
\begin{equation*}
  \begin{split}
    y_1 &= q_1 + 6 q_1^2 + 9 q_1^3 + 56 q_1^4 - 300 q_1^5 + \ldots \\
    y_2 &= q_2 \(1 - 2 q_1 + 5 q_1^2 - 32 q_1^3 + 286 q_1^4 - 3038 q_1^5-\ldots\)
  \end{split}
\end{equation*}
The mirror map identifies a suitable neighbourhood\footnote{This
  neighbourhood is $U_{\rm B}$ from Conjecture~\ref{conj:mirrorsym}.} of
$y_1=y_2=0$ in $\cM_{\FF_3}$ with the K\"{a}hler moduli space $\cM_{\rm A}$
of $\FF_3$; here $(q_1,q_2)$ are co-ordinates on $\cM_{\rm A}$.  This
identification matches up the Jacobi ring $J(W_y)$ of $W_y$ with the
small quantum cohomology algebra of $\FF_3$ at $(q_1,q_2)$ and the
residue pairing with the Poincar\'e pairing.

\subsection{The Landau--Ginzburg Mirror to $\PP(1,1,1,3)$}

The mirror family $\pi$ is
\begin{align*}
\pi\colon Z=(\Cstar)^4 & \longrightarrow \cM_{\PP(1,1,1,3)}=(\Cstar) \\
(w_1,w_2,w_3,w_5) &\longmapsto w_1w_2w_3 w_5^3
\end{align*}
the superpotential $W$ is
\[
W=w_1+w_2+w_3+w_5
\]
and the holomorphic volume form $\omega_y$ on the fiber
$Z_y=\pi^{-1}(y)$ is:
\[
\omega_y =\frac{d\log w_1\wedge d\log w_2\wedge d\log w_3 \wedge d\log w_5}
{d\log y}
\] 
The non-standard numbering of the co-ordinates here will be convenient
later.  The oscillating integrals \eqref{eq:oscint} satisfy the
Picard--Fuchs equation
\begin{equation}
  \label{eq:PFofP1113}
  D^3 (3D) (3D - z)(3D-2z)\, f = y\, f
\end{equation}
where $D = z y {\partial \over \partial y}$.

\subsection{Mirror Symmetry for $\PP(1,1,1,3)$}
\label{subsec:MS_P1113}

The mirror theorem we need here was proved by
Coates--Corti--Lee--Tseng \cite{CCLT}.  The $I$-function
\begin{equation}
  \label{eq:I-functionofP1113}
  I_{\PP(1,1,1,3)}(y,z) = z \, y^{p/z} \sum_{\substack{d \colon 3 d \in \ZZ
      \\ d \geq 0}} { y^d \over \prod_{\substack{b \colon \fr{b} = \fr{d} \\ 0 < b \leq d}} (p + bz)^3
    \prod_{\substack{b \colon \fr{b} = 0 \\ 0 < b \leq 3d}} (3p + bz)}
  \fun_{\fr{d}}
\end{equation}
also satisfies the Picard--Fuchs equation \eqref{eq:PFofP1113}.  It
coincides with the small $J$-function \eqref{eq:smallJP1113} after
the (trivial) change of variables $q=y$:
\[
J_{\PP(1,1,1,3)}(q,z) = I_{\PP(1,1,1,3)}(y,z)
\]
As discussed in Section~\ref{sec:VHSmirrorsym}, this proves
Conjecture~\ref{conj:mirrorsym} for $\PP(1,1,1,3)$.  The moduli space
$\cM_{\PP(1,1,1,3)}$ with co-ordinate $y$ is identified via the map
$q=y$ with the K\"ahler moduli space $\cM_{\rm A}$ of $\PP(1,1,1,3)$; here
$q$ is once again a co-ordinate on $\cM_{\rm A}$.  This identification
matches the Jacobi ring $J(W_q)$ with the small quantum orbifold
cohomology algebra of $\PP(1,1,1,3)$ at $q$, and the residue pairing
with the Poincar\'e pairing.

\subsection{Constructing the B-model \VHSns}
\label{sec:BVSHSP1113}
\begin{figure}[ht!]
  \begin{center}
    \begin{picture}(220,50)
      \put(150,5){\vector(1,0){40}}
      \put(150,5){\vector(-3,1){120}}
      \put(150,5){\vector(0,1){40}}
      \multiput(30,5)(40,0){5}{\makebox(0,0){$\cdot$}}
      \multiput(30,45)(40,0){5}{\makebox(0,0){$\cdot$}}
      \put(195,0){$p_1$}
      \put(18,40){$\ap_1$}
      \put(150,50){\makebox(0,0){$\ap_2=p_2$}}
    \end{picture}
    \caption{The secondary fan for ${\FF_3}$.}
    \label{fig:secondaryF_3}
  \end{center}
\end{figure}
We now extend the Landau--Ginzburg mirror of $\FF_3$ to a
Landau--Ginzburg model with a larger base, defined in terms of the
secondary fan for $\FF_3$ (Figure \ref{fig:secondaryF_3}).  Take
$w_1,w_2,w_5$ as co-ordinates on the fiber $Z_y$, so that:
\begin{equation}
  \label{eq:LGF3patch1}
  W_y=w_1+w_2+\frac{y_1y_2^3}{w_1w_2w_5^3}+\frac{y_2}{w_5}+w_5
\end{equation} 
The toric orbifold $\ccM$ associated to the secondary fan for $\FF_3$
gives a compactification of $\cM_{\FF_3}$.  One co-ordinate patch on
$\ccM$ comes from the K\"{a}hler cone of $\FF_3$, which is the cone in
the secondary fan spanned by $p_1=(1,0)$ and $p_2=(0,1)$.  The vectors
$p_1, p_2$ are dual to the co-ordinates $y_1,y_2$ on $\cM_{\FF_3}$.
The adjacent cone, spanned by $\ap_1=(-3,1)$ and $\ap_2=(0,1)$,
defines another co-ordinate patch on $\ccM$: let $\aq_1,\aq_2$ be the
co-ordinates dual to $\ap_1,\ap_2$.  The two co-ordinate systems are
related by:
\begin{align*}
\aq_1=y_1^{-1/3}&& \aq_2= y_1^{1/3} y_2
\end{align*}
Note that $\aq_1,\aq_2$ are multi-valued and so are not honest
co-ordinates on $\ccM$.  One should think of $\ccM$ as an orbifold and
of $\aq_1,\aq_2$ as a uniformizing system\footnote{This is the mirror
  partner of Remark~\ref{rem:orbifoldpoint}. }  near a $\ZZ/3\ZZ$
quotient singularity at $(\aq_1,\aq_2)=0$.  In the co-ordinates
$(\aq_1,\aq_2)$, we have:
\begin{equation}
  \label{eq:LGF3patch2}
  W_y=w_1+w_2+\frac{\aq_2^3}{w_1w_2w_5^3}+\frac{\aq_1\aq_2}{w_5}+w_5
\end{equation}
We can therefore extend the family of tori $\pi$ and the
superpotential $W_y$ across the locus $\{\aq_1=0\}$, where we see
the Landau--Ginzburg mirror of $\PP(1,1,1,3)$:
\begin{align*}
W=w_1+w_2+w_3+w_5, && w_1w_2w_3w_5^3=\aq_2^3 .
\end{align*}
The locus $\aq_1=0$, $\aq_2\neq 0$ in $\ccM$ is identified with the base
$\cM_{\PP(1,1,1,3)}$ of the Landau--Ginzburg mirror of $\PP(1,1,1,3)$
via the map $y = \aq_2^3$.

\begin{figure}[ht!]
  \begin{center}
    \begin{picture}(265,90)(0,-10)
      \put(40,23){\line(1,0){200}}
      \put(90,5){\line(0,1){70}}
      \put(150,5){\line(0,1){60}}
      \multiput(220,5)(0,5){10}{\line(0,1){2.5}}
      \put(73,-7){$y_1=0$}
      \put(132,-7){$y_1=-\frac{1}{27}$}
      \put(202,-7){$y_1=\infty$}
      \put(202,58){$\aq_1=0$}
      \put(0,20){$y_2=0$}
      \put(250,20){$\aq_2=0$} 
    \end{picture}
    \caption{The B-model moduli space $\cM_{\rm B}$.}
    \label{fig:kaehlermoduliofF_3}
  \end{center}
\end{figure}

The base of our extended Landau--Ginzburg model, which we call the
\emph{B-model moduli space} $\cM_{\rm B}$, is obtained from $\ccM$ by
deleting the closures of the loci $ \{y_1y_2=0\}$ and $\{y_1=-1/27\}$.
Equations \eqref{eq:LGF3patch1} and \eqref{eq:LGF3patch2} define a
Landau--Ginzburg model over $\cM_{\rm B}$ which contains the
Landau--Ginzburg mirrors for $\FF_3$ and for $\PP(1,1,1,3)$ as
subsets.  The limit points $y_1=y_2=0$ and $\aq_1=\aq_2=0$ of $\cM_{\rm B}$
are called the \emph{large radius limit points} or \emph{cusps}
corresponding respectively to $\FF_3$ and to $\PP(1,1,1,3)$.

Let $\cE$ denote the B-model \VHS with base $\cM_{\rm B}$ defined by the
Landau--Ginzburg model just described.  (See Definition
\ref{dfn:mirrorDmodule} for the B-model \VHSns.)\phantom{.}  We equip $\cE$
with the grading operator $\Gr \colon \cE\rightarrow \cE$ defined by
\begin{align*}
\Gr \big[f(x,z)e^{W_y/z}\, \omega_y\big] = 
\left[\(2 z\partial_z + 2\sum_{i=1}^5 w_i \partial_{w_i}\) f(x,z)
e^{W_y/z} \, \omega_y\right]
\end{align*}
This satisfies the axioms for a graded \VHS with Euler field and dimension
\begin{align*}
E=2y_2\parfrac{}{y_2}= 2\aq_2\parfrac{}{\aq_2}, && 
D=\dim \FF_3 =3.
\end{align*}

\begin{rem}
  The superpotentials \eqref{eq:LGF3patch1} and
  \eqref{eq:LGF3patch2} have isolated non-degenerate critical points,
  and so the small quantum cohomology algebra of $\FF_3$ and the
  small quantum orbifold cohomology algebra of $\PP(1,1,1,3)$ are
  semisimple.
\end{rem}

\begin{rem}
  The mirror $D$-module develops a singularity along the lines
  $\{y_1=0\}$ and $\{y_2=0\} = \{\aq_2=0\}$.  These are the solid
  lines in Figure~\ref{fig:kaehlermoduliofF_3}.  It is non-singular
  along the (dashed) line $\{\aq_1=0\}$.
\end{rem}
\subsection{An Opposite Subspace At Each Cusp}

We now characterize those opposite subspaces for the B-model \VHS
which give rise, via mirror symmetry, to the big quantum cohomology
Frobenius manifolds for $\FF_3$ and for $\PP(1,1,1,3)$.  As we will
see in the next section, these opposite subspaces are \emph{not}
mapped into each other under parallel transport from cusp to cusp.

The $I$-functions \eqref{eq:I-functionofF3},
\eqref{eq:I-functionofP1113} define $D$-module homomorphisms
\begin{equation}
  \label{eq:Imaps}
  \begin{split}
    & \mapI_{\FF_3} \colon \cE \longrightarrow H^\bullet(\FF_3)\otimes 
    \cOzzinv{\tcM_{\rm B}}, \\
    & \mapI_{\PP(1,1,1,3)} \colon \cE|_{\cM_{\PP(1,1,1,3)}}
    \longrightarrow H^\bullet_{\text{\rm orb}}\big(\PP(1,1,1,3)\big)\otimes 
    \cOzzinv{\tcM_{\PP(1,1,1,3)}}
  \end{split}
\end{equation}
where we give the right-hand sides the trivial $D$-module structure.
These homomorphisms are defined using the isomorphism
\eqref{eq:Dmodulegeneratedbyoscint}, by sending a representative
$P(y,z\partial,z)\in \cOz{\cM_{\rm B}}\langle z\cT_{\cM_{\rm
    B}}\rangle$ to $z^{-1} P I_{\FF_3}$ or to $z^{-1} P
I_{\PP(1,1,1,3)}$.  This does not depend on our choice of
representative $P$ since the $I$-functions satisfy the Picard--Fuchs
equations \eqref{eq:PFofF_3}, \eqref{eq:PFofP1113}.  $\mapI_{\FF_3}$
is \emph{a priori} defined only in a small neighbourhood of the cusp
for $\FF_3$, but it can be extended to the whole of $\tcM_{\rm B}$ by
analytic continuation (or, which amounts to the same thing, by solving
the Picard--Fuchs equations on this larger region).  The maps
\eqref{eq:Imaps} send $\nabla^z$-parallel sections of $\cE$ to
constant sections, and thus identify the space $\cH$ of flat sections
of $\cE$ with Givental's symplectic vector space (with $Q=1$):
\begin{equation}
  \label{eq:identify}
  \begin{aligned}
    & \mapI_{\FF_3} \colon \cH \cong 
    H^\bullet(\FF_3)\otimes \Czzinv =\cH_{\FF_3}|_{Q_i=1} \\
    & \mapI_{\PP(1,1,1,3)} \colon \cH \cong H^\bullet_{\rm orb}(\PP(1,1,1,3)) 
    \otimes \Czzinv=\cH_{\PP(1,1,1,3)}|_{Q=1}
  \end{aligned}
\end{equation}
As discussed in Section~\ref{sec:VHSmirrorsym}, the mirror isomorphism
\eqref{eq:mirroriso} sends the generator $[\exp(W_y/z)\,\omega_y]$ of
$\cE$ to the family of vectors $y\mapsto I(y,z)=J(q,z)$ lying on
Givental's Lagrangian submanifold $\cL$.  The identifications
\eqref{eq:identify} are exactly those induced by \eqref{eq:mirroriso}.

An easy calculation using the explicit forms of the $I$-functions
yields:
\begin{pro}
\label{pro:matchingofgrading}
Under the identifications \eqref{eq:identify}, the B-model grading
operator corresponds to the A-model grading operator
\eqref{eq:gradingonH}. \hfill \qedsymbol
\end{pro}

The Hodge structure $\E_y$ near cusps behaves as follows.  As
$(y_1,y_2)\to 0$, we have:
\begin{equation}
\label{eq:asymptoticofHodgeF3}
\mapI_{\FF_3}(\E_y)  \sim e^{(p_1\log y_1 + p_2 \log y_2)/z}  
\Big(H^\bullet(\FF_3) \otimes \Cz + O(y_1,y_2)\Big)
\end{equation}
and as $\aq_2\to 0$ with $\aq_1=0$, we have:
\begin{equation}
\label{eq:asymptoticofHodgeP1113}
\mapI_{\PP(1,1,1,3)}(\E_y) \sim 
e^{3p\log \aq_2/z} \Big(H^\bullet_{\rm orb}\big(\PP(1,1,1,3)\big)\otimes \Cz 
+O(\aq_2)\Big)
\end{equation}
These are semi-infinite analogs of Schmid's Nilpotent Orbit Theorem
\cite{Schmid} --- the singularity of the Hodge structure near a cusp
is asymptotically given by the exponential of nilpotent operators.
This corresponds to the fact that $\cE$ has quasi-unipotent monodromy
at each cusp, and so has a regular singular extension (Deligne's
extension) on a finite cover of a neighbourhood of each cusp.  

\begin{dfn}
  The \emph{limiting Hodge structure} at a cusp $c$ is defined to be a
  subspace $\E^{\rm lim}_c$ of $\cH$ satisfying
  \[
  \E_{(x_1,x_2)} \sim \exp\left(\sum_{i=1}^2 N_i \log x_i\right) 
  \Big(\E^{\rm lim}_c+O(x_1,x_2)\Big)
  \] 
  where $N_i$ are nilpotent operators and $x_1, x_2$ are local
  co-ordinates centered at $c$ such that $\cE$ has a logarithmic
  singularity along the $x_i$-axes.  The limiting Hodge structure
  actually depends on the choice of such co-ordinates\footnote{A
    co-ordinate change of the form $\log x_i'=\log x_i + f_i(x_1,x_2)$
    with $f_i(0,0)=0$ does not change the limiting Hodge structure, so
    $\E_c^{\rm lim}$ depends only on the choice of ``origin'' of $\log
    x_i$.}, but in our examples we take $(x_1,x_2)=(y_1, y_2)$ near
  the cusp $c_1$ for $\FF_3$ and $(x_1,x_2)=(\aq_1,\aq_2)$ near the
  cusp $c_2$ of $\PP(1,1,1,3)$.
\end{dfn}

The above calculation shows that:
\begin{align*}
\mapI_{\FF_3}(\E_{c_1}^{\rm lim})&= 
H^\bullet(\FF_3)\otimes \Cz=\cH_{\FF_3}^+|_{Q_i=1} \\
\mapI_{\PP(1,1,1,3)}(\E_{c_2}^{\rm lim})& 
=H^\bullet_{\rm orb}(\PP(1,1,1,3))\otimes \Cz
=\cH_{\PP(1,1,1,3)}^+|_{Q=1}  
\end{align*} 

We now construct an opposite subspace $\cH_-$ for each cusp $c$.  We
postulate that $\cH_-$ should satisfy:
\begin{itemize}
\item[(a)] $\cH_-$ is opposite to the limiting Hodge structure
  $\E_c^{\rm lim}$
\item[(b)] $\cH_-$ is preserved by the grading operator $\Gr$
\item[(c)]$\cH_-$ is invariant under local monodromy, and moreover the
  monodromy action $M$ satisfies $M^N=\id$ on $z\cH_-/\cH_-$ where $N$
  is the order of the local isotropy group at the cusp $c$
\end{itemize} 
For $\FF_3$, the local monodromy means the monodromy around the axes
$y_1=0$ and $y_2=0$.  The corresponding monodromy actions on
$\cH_{\FF_3}$ are given by
\begin{align*}
M_1 = \exp\big(2\pi\tti p_1/z\big), && M_2 = \exp\big(2\pi\tti p_2/z\big).
\end{align*}
For $\PP(1,1,1,3)$, the local monodromy means the monodromy
coming from an orbifold loop $[0,1]\ni t\mapsto
(\aq_1,\aq_2)=\big(0,e^{2\pi\tti t/3}\big)$.  The corresponding action on
$\cH_{\PP(1,1,1,3)}$ is given by
\[
M = M_0 \exp\big(2\pi\tti p/z\big)
\]
where
\begin{align*}
M_0\big|_{H^\bullet(\PP(1,1,1,3))}=\id, && M_0\big(\fun_{1 \over 3}\big)=\alpha \fun_{1 \over 3}, 
&& M_0\big(\fun_{2 \over 3}\big)=\alpha^2 \fun_{2 \over 3} 
\end{align*}
and $\alpha=\exp\big(2\pi\tti/3\big)$.  We have $N=1$ for $\FF_3$ and
$N=3$ for $\PP(1,1,1,3)$.  Condition (c) above implies that the
connection 1-form $A$ in Proposition~\ref{pro:normalformofconn} is
well-defined on an $N$-fold cover (uniformizing system) of a
neighbourhood of each cusp.
 
\begin{thm} \label{thm:uniqueness} Let $c_1$ and $c_2$ be the cusps of
  $\cM_{\rm B}$ corresponding to $\FF_3$ and to $\PP(1,1,1,3)$ respectively.
  For each $c_i$ there exists a unique opposite subspace $\cH^-_{c_i}$
  satisfying conditions (a--c) above.  Moreover, $\cH^-_{c_i}$
  corresponds under \eqref{eq:identify} to the standard A-model
  opposite subspace:
  \begin{align*}
    \mapI_{\FF_3}\big(\cH^-_{c_1}\big) &= z^{-1}H^\bullet(\FF_3)\otimes \Czinv
    =\cH_{\FF_3}^-|_{Q_i=1} \\
    \mapI_{\PP(1,1,1,3)}\big(\cH^-_{c_2}\big) &= z^{-1}H^\bullet_{\text{\rm
      orb}}\big(\PP(1,1,1,3)\big)\otimes\Czinv =\cH_{\PP(1,1,1,3)}^-|_{Q=1}
  \end{align*}
\end{thm}
\begin{proof}
  We give a proof only for $\PP(1,1,1,3)$.  The $\FF_3$ case is
  similar and easier.  Throughout the proof we identify $\cH$ with
  $\cH_{\PP(1,1,1,3)}|_{Q=1}$ via the map $\mapI_{\PP(1,1,1,3)}$ and
  write $\cH_-$ for $\mapI_{\PP(1,1,1,3)} \cH_{c_2}^-$.

  Proposition \ref{pro:matchingofgrading} implies that the A-model
  grading operator \eqref{eq:gradingonH} preserves $\cH_-$:
  \begin{align*}
    \Gr = 2 z \partial_z + \Gr_0 - 2 c_1\big(\PP(1,1,1,3)\big)/z, &&
    c_1\big(\PP(1,1,1,3)\big)=6p.
  \end{align*}
  On the other hand, the logarithm $6\pi\tti p/z = \log (M^3)$ of
  the cube of the monodromy preserves $\cH_-$, and so the ``usual''
  grading operator $2z\partial_z + \Gr_0$ also preserves $\cH_-$.
  This means that $\cH_-$ is a homogeneous subspace of $\cH$.  Because
  $\cH_-$ is opposite to $\E_{c_2}^{\rm lim}$, there is a unique
  $\CC$-basis $\{\psi_0,\dots,\psi_3,\phi_1,\phi_2\}$ of $z\cH_-\cap
  \E_{c_2}^{\rm lim}$ such that:
  \begin{align*}
    \psi_i = p^i + O(z), && \phi_1=\fun_{1 \over 3}+O(z), && 
    \phi_2 = \fun_{2 \over 3}+ O(z)  
  \end{align*}
  These elements must be homogeneous.  Since both $M$ and
  $\exp\big(2\pi\tti p/z\big)$ preserve $\cH_-$, $M_0$ must also preserve
  $\cH_-$.  It is clear that $M_0$ preserves $\E^{\rm lim}_{c_2}$, so
  it acts on $z\cH_-\cap \E^{\rm lim}_{c_2}$.  Thus $z\cH_-\cap
  \E^{\rm lim}_{c_2}$ decomposes into eigenspaces for $M_0$; it
  follows that $\psi_i\in H^\bullet\big(\PP(1,1,1,3)\big)\{z\}$ and $\phi_i\in
  H^\bullet\big(\PP(V^{i/3})\big)\{z\}$.  Homogeneity now implies that
  $\phi_i=\fun_{i/3}$.

  Since $M^3$ acts trivially on $z\cH_-/\cH_-$ it follows that
  $\frac{1}{3}\log M^3$ sends $z\cH_-$ to $\cH_-$.  Thus
  $\frac{1}{3}z\log M^3=2\pi\tti p$ preserves $\E^{\rm lim}_{c_2}$
  and $z\cH_-$ simultaneously, and therefore acts on $z\cH_-\cap
  \E^{\rm lim}_{c_2}$.  By homogeneity again, we can write
  \[
  \psi_i = p^i + \sum_{j=0}^{i-1} c_{ij} z^{i-j} p^j
  \]
  for some $c_{ij}\in\CC$.  As $p^{4-i} \psi_i\in z\cH_-\cap \E^{\rm
    lim}_{c_2}$ and $p^{4-i} \psi_i$ is divisible by $z$, we know that
  $p^{4-i}\psi_i$ should be zero.  This shows that $c_{ij}=0$ and that
  $\psi_i=p^i$.  Since $z\cH_-$ is spanned over $\Czinv$ by
  $\psi_0,\dots,\psi_3$, $\phi_1$, $\phi_2$, it follows that $\cH_- =
  \cH_{\PP(1,1,1,3)}^-|_{Q=1}$.
\end{proof}

Monodromy properties also force, as we now show, the A-model and
B-model pairings to coincide up to an overall scalar factor.  In the
next section we will see that the composition $\mapI_{\FF_3}\circ
\mapI_{\PP(1,1,1,3)}^{-1}$ exactly preserves Givental's symplectic
form.

\begin{pro}
  \label{pro:matchingofpairing}
  Under the identifications \eqref{eq:identify}, the B-model
  symplectic form corresponds to a scalar multiple of Givental's
  symplectic form.  In particular, $\cH_{c_i}^-$ is isotropic with
  respect to the B-model symplectic form.
\end{pro} 
\begin{proof}
  We give a proof only for $\PP(1,1,1,3)$; the $\FF_3$ case is
  similar.  Let $\Ppair{\cdot}{\cdot}_{\rm B}$ be the $\Czzinv$-valued
  pairing on $\cH_{\PP(1,1,1,3)}|_{Q=1}$ induced by the B-model
  pairing.  The definition of the B-model pairing shows that this is
  monodromy-invariant:
  \[
  \Ppair{M\alpha}{M\beta}_{\rm B} = \Ppair{\alpha}{\beta}_{\rm B}
  \]
  Using $6\pi\tti p/z=\log M^3$, we have
  \begin{equation}
    \label{eq:pisselfadjoint}
    \Ppair{p\alpha}{\beta}_{\rm B} = \Ppair{\alpha}{p \beta}_{\rm B}
  \end{equation}
  and so $M_0 = M \exp\big({-2}\pi\tti p/z\big)$ also preserves the pairing
  $\Ppair{\cdot}{\cdot}_{\rm B}$.  This implies that
  \begin{align}
    \label{eq:orthogonality}
    \Ppair{H^\bullet\big(\PP(1,1,1,3)\big)}{\fun_{i/3}}_{\rm B} = 0 && \text{and} &&
    \Ppair{\fun_{i/3}}{\fun_{i/3}}_{\rm B} = 0  
  \end{align}
  for $i=1,2$.  

  From the asymptotics \eqref{eq:asymptoticofHodgeP1113} we know that
  for each $\alpha\in \cH_{\PP(1,1,1,3)}^+|_{Q=1}$, there exists a
  family of elements $\{\alpha_{\aq_2}\}$ in
  $\cH_{\PP(1,1,1,3)}|_{Q=1}$ such that
  \[
  \alpha_{\aq_2} = \exp\big(3p\log \aq_2/z\big) \Big(\alpha+ O(\aq_2)\Big) \in
  \E_{\aq_1=0, \aq_2}.
  \] 
  For any $\alpha,\beta\in \cH_{\PP(1,1,1,3)}^+$, the B-model
  pairing of $\alpha_{\aq_2}$ and $\beta_{\aq_2}$ takes values in
  $\Cz$:
  \[
  \Ppair{\alpha_{\aq_2}}{\beta_{\aq_2}}_{\rm B} =
  \Ppair{\alpha+O(\aq_2)}{\beta+O(\aq_2)}_{\rm B} \in \Cz.
  \]
  Taking the limit $\aq_2\to 0$, we see that $\Ppair{\alpha}{\beta}_{\rm B}$
  is in $\Cz$.

  The compatibility of grading and the pairing gives that for
  homogeneous elements $\alpha,\beta\in H^\bullet_{\rm
    orb}(\PP(1,1,1,3))$ we have:
  \begin{align*}
    2 z\partial_z \Ppair{\alpha}{\beta}_{\rm B}
    &= (\deg \alpha + \deg \beta - 6) \Ppair{\alpha}{\beta}_{\rm B} \\
    &\qquad -\Ppair{(c_1(\cX)/z)\cup \alpha}{\beta}_{\rm B}
    -\Ppair{\alpha}{(c_1(\cX)/z)\cup \beta}_{\rm B}
  \end{align*} 
  The second line vanishes by \eqref{eq:pisselfadjoint}.  This means
  that $\Ppair{\alpha}{\beta}_{\rm B}$ is homogeneous of degree $\deg
  \alpha+\deg \beta-6$.  From this homogeneity, equation
  \eqref{eq:pisselfadjoint}, the orthogonality
  \eqref{eq:orthogonality}, and the fact that
  $\Ppair{\alpha}{\beta}_{\rm B}\in \Cz$, it follows that
  the only non-vanishing pairings among basis elements are
  \[
  \Ppair{p^i}{p^j}_{\rm B} = \Ppair{\fun}{p^{i+j}}_{\rm B}\in \CC, \quad
  \Ppair{\fun_{1 \over 3}}{\fun_{2 \over 3}}_{\rm B} =
  \Ppair{\fun_{2 \over 3}}{\fun_{1 \over 3}}_{\rm B} \in \CC
  \] 
  with $i+j=3$.  This shows directly that $\cH^-_{\PP(1,1,1,3)}$
  is isotropic with respect to the B-model symplectic form.  But the
  general theory of \VHS in Section \ref{sec:VHS} implies that the
  B-model pairing also satisfies (see equation
  \eqref{eq:Frobeniusalgebra})
  \begin{align*}
    \Ppair{\alpha\circ_q \beta}{\gamma}_{\rm B} = 
    \Ppair{\alpha}{\beta\circ_q \gamma}_{\rm B} &&
    \text{for all $\alpha,\beta,\gamma \in H_{\text{\rm
        orb}}^\bullet\big(\PP(1,1,1,3)\big)$} 
  \end{align*} 
  because the small quantum orbifold cohomology algebra is generated by
  $p$.  Thus $\Ppair{\cdot}{\cdot}_{\rm B}$ is completely determined by the
  value $\Ppair{\fun}{p^3}_{\rm B}$ and is proportional to the orbifold
  Poincar\'{e} pairing.
\end{proof} 

\begin{rem}
  In proving the uniqueness in general of opposite subspaces and
  pairings which behave well under monodromy, the hard Lefschetz
  property of the \emph{usual} cohomology of a projective orbifold
  will play an important role.  
  This will be explained in Iritani \cite{Iritani:wallcrossings3}. 
  See also Iritani \cite[Theorem 3.13]{Iritani:RIMS} 
  for the uniqueness of opposite subspaces in the A-model. 
  In the proofs above, we implicitly
  used the hard Lefschetz property of $H^\bullet(\PP(1,1,1,3))$.  A
  hard Lefschetz property for \emph{orbifold} cohomology is discussed
  in Theorem \ref{thm:hardLefschetz} below.    
\end{rem} 

\begin{dfn} \emph{A polarization of $\cH$} at a cusp $c$ is a
  decomposition
  \[
  \cH =  \E^{\rm lim}_{c}\oplus \cH_c^-
  \]
  where $\E^{\rm lim}_{c}$ is the limiting Hodge structure and
  $\cH_c^-$ is an opposite subspace.
\end{dfn}

The polarization at a cusp will be mapped by $\mapI_{\FF_3}$ or
$\mapI_{\PP(1,1,1,3)}$ to the standard polarization:
\begin{align*}
\cH_{\FF_3}^+\oplus \cH_{\FF_3}^-\Big|_{Q_i=1} &&  
\text{or} &&
\cH_{\PP(1,1,1,3)}^+\oplus \cH_{\PP(1,1,1,3)}^-\Big|_{Q=1}
\end{align*}

\subsection{The Polarizations Are Different}
\label{sec:differenceofpolarization_F3}
We now compare the polarizations at the cusps of $\cM_{\rm B}$ corresponding
to $\FF_3$ and to $\PP(1,1,1,3)$.  
Let $\Ubar \colon \cH_{\PP(1,1,1,3)}|_{Q=1}
\rightarrow \cH_{\FF_3}|_{Q_i=1}$ be the linear transformation 
defined by the composition 
\[
\begin{CD}
\cH_{\PP(1,1,1,3)}|_{Q=1} 
 @>{(\mapI_{\PP(1,1,1,3)})^{-1}}>> \cH @>{\mapI_{\FF_3}}>> 
\cH_{\FF_3} |_{Q_i=1}  
\end{CD}
\]
and let $\U: \cH_{\PP(1,1,1,3)}|_{Q=1}
\rightarrow \cH_{\FF_3}|_{Q_i=1}$ be $\Ubar$
followed by changing the sign of $z$.  The transformation $\Ubar$ (or
equivalently $\U$) measures the difference between the polarizations
at the two cusps.  As we will see, the sign flip in the definition of
$\U$ comes from the sign flip which relates the A-model \VHS to the
tangent spaces to Givental's Lagrangian submanifold (see
Section~\ref{sec:bigQCandcone}).

\begin{pro}
\label{pro:symplectictransformationP1113}
The matrix of $\U$ with respect to the bases for $H^\bullet_{\text{\rm
  orb}}\big(\PP(1,1,1,3)\big)$ and $H^\bullet(\FF_3)$ defined in
Section~\ref{sec:cohomology} is:
\begin{equation}
  \label{eq:matrixofU}
  \begin{pmatrix} 
    1 & 0 &0 & 0&0 &0  \\
    0 & 1 &0 &0 &0 &0  \\
    0 & 0 &1 &0 & 0 &0  \\
    0 & 0 &0 &0 & \frac{2\sqrt{3}\pi}{3\Gamma(\frac{1}{3})^3}z & 
    \frac{2\sqrt{3}\pi}{3\Gamma(\frac{2}{3})^3}  \\
    -\frac{\pi^2}{3 z^2} & 0 &0 &0 &\frac{2\pi^2}{3\Gamma(\frac{1}{3})^3} &
    -\frac{2\pi^2}{3\Gamma(\frac{2}{3})^3 z} \\
    \frac{8\zeta(3)}{z^3} &0 &0 & 1 &
    \frac{2\sqrt{3}\pi^3}{9\Gamma(\frac{1}{3})^3z} &
    \frac{2\sqrt{3}\pi^3}{9\Gamma(\frac{2}{3})^3z^2} 
  \end{pmatrix}
\end{equation}
Here $\zeta$ is the Riemann zeta function.  The linear transformation
$\U$ preserves the grading and the symplectic forms but does not
preserve the standard opposite subspaces.
\end{pro} 
\begin{proof}
  The $I$-functions $I_{\FF_3}$ and $I_{\PP(1,1,1,3)}$ are the images
  of the generator $[\exp(W_y/z) \, \omega_y]\in \cE$ under the maps
  $\mapI_{\FF_3}$ and $\mapI_{\PP(1,1,1,3)}$.  It follows that
  \[
  \Ubar \big(I_{\PP(1,1,1,3)}\big) = I_{\FF_3}|_{\aq_1=0}
  \]
  where we regard $I_{\PP(1,1,1,3)}$ as a function of $\aq_2$ via the
  map $y = \aq_2^3$ discussed above
  Figure~\ref{fig:kaehlermoduliofF_3}.  We calculate $\Ubar$ (and
  hence $\U$) by analytically continuing $I_{\FF_3}$ to a
  neighbourhood of the large radius limit point for $\PP(1,1,1,3)$ and
  then comparing it with $I_{\PP(1,1,1,3)}$.  Using the Barnes method
  (see the Appendix A), one finds:
  \begin{multline}
    \label{eq:IF3continued}
    I_{\FF_3}(\aq_1,\aq_2,z) = \textstyle z \,
    \Gamma\big(1+\frac{p_1}{z}\big)^3 \,
    \Gamma\big(1+\frac{p_2}{z}\big)
    \Gamma\big(1+\frac{p_2 - 3 p_1}{z}\big)   \times \\
    \sum_{k,l\ge 0} 
    \frac{(-1)^{k+l}\sin\big(\frac{p_2-3p_1}{z}\pi\big)}
    {3\sin\big(\frac{p_2 - 3 p_1}{3z}\pi + \frac{l-k}{3}\pi\big)}
    \frac{\aq_1^k\aq_2^{l+p_2/z}} {k! \, z^{2l} \,
      \Gamma\big(1+\frac{p_2}{3z}+\frac{l-k}{3}\big)^3 \,
      \Gamma\big(1+\frac{p_2}{z}+l\big) }
  \end{multline}
  We compare this with:
  \begin{align*}
    I_{\PP(1,1,1,3)}(\aq_2,z) = z\sum_{m\ge 0}
    \frac{\Gamma\big(1-\langle{-\frac{m}{3}}\rangle+\frac{p}{z}\big)^3 \,
      \Gamma\big(1+\frac{3p}{z}\big)} 
    {\Gamma\big(1+\frac{p}{z}+\frac{m}{3}\big)^3 \,
      \Gamma\big(1+\frac{3p}{z}+m\big)} 
    \frac{\aq_2^{m+3p/z}}{z^{2m}}
    \frac{\fun_{\langle \frac{m}{3}\rangle }} {z^{3\langle
        {-\frac{m}{3}}\rangle }}.
  \end{align*}
  Since $\Ubar$ is equivariant with respect to the monodromy action
  around the axis $\{y_2=0\}=\{\aq_2=0\}$, we have $\Ubar e^{2\pi
    \tti 3p/z}= e^{2\pi \tti p_2/z}\Ubar$ and so $\Ubar 3p =
  p_2 \Ubar$.  Thus:
  \begin{align*}
    \Ubar\big(\fun_0\big) & = 
    \frac{
      \Gamma\big(1+\frac{p_1}{z}\big)^3 \,
      \Gamma\big(1+\frac{\ap_1}{z}\big)
    } 
    {
      \Gamma\big(1+\frac{p_2}{3z}\big)^3
    }
    \frac{
      \sin\big(\frac{\ap_1}{z}\pi\big)
    }
    {
      3\sin\big(\frac{\ap_1}{3z}\pi\big)
    } \displaybreak[0]\\
    z^{-1} \Ubar\big(\fun_{2 \over 3}\big) & =
    \frac{
      \Gamma\big(1+\frac{p_1}{z}\big)^3 \,
      \Gamma\big(1+\frac{\ap_1}{z}\big)
    }
    {
      \Gamma\big(\frac{2}{3}+\frac{p_2}{3z}\big)^3
    }
    \frac{
      \sin\big(\frac{\ap_1}{z}\pi\big)
    }
    {
      3\sin\big(\frac{\ap_1}{3z}\pi+\frac{2}{3}\pi\big)
    } \displaybreak[0]\\
    z^{-2} \Ubar\big(\fun_{1 \over 3}\big) & =
    -\frac{
      \Gamma\big(1+\frac{p_1}{z}\big)^3 \,
      \Gamma\big(1+\frac{\ap_1}{z}\big)
    }
    {
      \Gamma\big(\frac{1}{3}+\frac{p_2}{3z}\big)^3
    }
    \frac{
      \sin\big(\frac{\ap_1}{z}\pi\big)
    }
    {
      3\sin\big(\frac{\ap_1}{3z}\pi+\frac{1}{3}\pi\big)
    }
  \end{align*} 
  where $\ap_1=p_2-3p_1$, and the conclusion follows.  The value
  $\zeta(3)$ in the matrix for $\U$ comes from the expansion of the
  $\Gamma$-function:
  \[
  \Gamma(1+x)=\exp\(\textstyle -\gamma x + \frac{\pi^2}{12} x^2 -
  \frac{\zeta(3)}{3}x^3+O(x^4)\) 
  \]
  where $\gamma$ is Euler's constant.  The transformation $\U$ does
  not map $\cH_{\PP(1,1,1,3)}^-$ to $\cH_{\FF_3}^-$ because the matrix
  \eqref{eq:matrixofU} contains strictly positive powers of $z$.
\end{proof}

\begin{rem}
  The symplectic transformation $\U$ always has an ambiguity due to
  the monodromy action on $\cH$.  This corresponds to the choice of
  branch cuts in the process of analytic continuation.
\end{rem}

\begin{rem}
  A closely-related symplectic transformation (with $z=1$) occurs in
  work of Aganagic--Bouchard--Klemm \cite{ABK}*{equation~6.21}.  They
  studied a phase-transition from local $\PP^2$ (the total space of
  the canonical bundle $K_{\PP^2}$) to $\CC^3/\ZZ_3$.  Our example
  here is a global version of this but is not Calabi--Yau.
\end{rem} 

\subsection{The Proof of Conjecture~\ref{conj}}

In Theorem~\ref{thm:transformationofLagrangiancones} below, 
we prove Conjecture~\ref{conj:preciseCRC} when 
$\cX=\PP(1,1,1,3)$ and $Y=\FF_3$. 
Conjecture~\ref{conj} in the Introduction follows from 
this and the definition (\ref{eq:analyticLag}) of 
the analytic Givental's cone. 


\begin{thm}
\label{thm:transformationofLagrangiancones} 
Let $\E^{\FF_3}_\tau \subset \left.\cH_{\FF_3}\right|_{Q_i=1}$ and
$\E^{\PP(1,1,1,3)}_\tau \subset \left.\cH_{\PP(1,1,1,3)}\right|_{Q=1}$ be
the moving subspace realizations \eqref{eq:moving_A} 
of the analytic big A-model 
\VHSns s of $\FF_3$ and $\PP(1,1,1,3)$.  
Then there exists a map $\Upsilon$ from an open subset 
of $H_{\rm orb}^\bullet(\PP(1,1,1,3))$ 
to an open subset of $H^\bullet(\FF_3)$ 
such that, after analytic continuation if necessary,  
\[
\ov\U \big(\E^{\PP(1,1,1,3)}_\tau \big) 
= \E^{\FF_3}_{\Upsilon(\tau)}
\]
where $\U$ is the symplectic transformation from
Proposition~\ref{pro:symplectictransformationP1113}.
\end{thm}
\begin{proof} 
  We use Dubrovin's Reconstruction Theorem \cite{Dubrovin}.  This
  implies that the Frobenius manifold given by big quantum (orbifold)
  cohomology can be uniquely reconstructed from one semisimple fiber
  as an isomonodromic deformation of the differential equation
  \begin{equation}
  \label{eq:diffeq_z} 
  \left[\frac{1}{2} \Gr - \frac{1}{z}\nabla^z_E\right ] \psi = 
  \left[ z\parfrac{}{z} - \frac{1}{z} E\bullet_\tau 
  + \frac{1}{2}\Gr_0 \right]\psi = 0
  \end{equation} 
  where $\Gr,\Gr_0$ are grading operators.  
  We know from the mirror analysis
  that the small quantum cohomology algebras of $\FF_3$ and
  $\PP(1,1,1,3)$ are semisimple.  The big quantum cohomologies of
  $\FF_3$ and $\PP(1,1,1,3)$ are therefore determined \emph{as
    analytic Frobenius manifolds} by the small quantum cohomologies.
  
  Write $\cX = \PP(1,1,1,3)$ and $Y=\FF_3$.  Let $\mir_Y \colon
  \tcM_{\rm B} \to H^2(Y)$ and $\mir_\cX \colon \tcM_{\rm B}\supset
  \tcM_\cX \to H^2(\cX)$ be (analytic continuations of) the mirror
  maps.  They are given by the coefficients of $z^{-1}$ of the
  $I$-functions $I_{\FF_3}, I_{\PP(1,1,1,3)}$ in equations
  (\ref{eq:I-functionofF3}), (\ref{eq:I-functionofP1113}).  The mirror
  theorems discussed in Sections \ref{subsec:MS_F3},
  \ref{subsec:MS_P1113} imply that
  \begin{align*}
    \mapI_{Y}(\cE_y) = \E^{Y}_{\mir_Y(y)} && \text{and} &&
    \mapI_{\cX}(\cE_y) = \E^{\cX}_{\mir_\cX(y)}
  \end{align*}
  for $y\in \tcM_{\rm B}$ and $y\in \tcM_{\cX}$ respectively.  These
  equations hold \emph{a priori} in neighbourhoods of the cusps, but
  hold everywhere by analytic continuation.  By the definition of
  $\ov\U$, we have
  \begin{equation} 
    \label{eq:UmatchesE_small}
    \ov\U\big(\E^\cX_{\mir_\cX(y)}\big) = \E^Y_{\mir_Y(y)}
  \end{equation}
  for $y\in \tcM_\cX$. Take a semisimple point $y_0\in \tcM_\cX$ (in
  fact every point on $\tcM_\cX$ is semisimple) and a small open
  neighbourhood $U_0$ of $y_0$ in $\tcM_\cX$.  Since the B-model and
  A-model grading operators match (Proposition
  \ref{pro:matchingofgrading}) and the Euler vector field is tangent
  to $\tcM_\cX$, $\ov\U$ induces an isomorphism of graded \VHSns s:
  \begin{equation}
  \label{eq:VHSisom_restricted}
  \cE^{\cX, \rm big}_{\rm A}|_{\mir_\cX(U_0)} \cong 
  \cE^{Y, \rm big}_{\rm A}|_{\mir_Y(U_0)}
  \end{equation} 
  Take a homogeneous opposite subspace $\cH^-$ of $\cE$ at $y_0$.
  This gives rise to the opposite subspaces $\mapI_\cX(\cH^-)$ and
  $\mapI_Y(\cH^-)$ of $\cE^{\cX, \rm big}_{\rm A}$ and $\cE^{Y,\rm
    big}_{\rm A}$ --- these subspaces are opposite in neighbourhoods
  of $\sigma_0:=\mir_\cX(y_0)$ and $\tau_0:=\mir_Y(y_0)$ respectively
  --- and produces Frobenius manifold structures\footnote{These
    Frobenius manifold structures are not in general the quantum
    cohomology Frobenius manifold structures, because in general
    $\mapI_\cX(\cH^-) \ne \cHX^-$ and $\mapI_Y(\cH^-) \ne \cHY^-$.}
  on the analytic germs $\big(\HorbX,\tau_0\big)$ and
  $\big(\HY,\tau_0\big)$.  Since these two Frobenius manifolds are the
  unfolding of \emph{the same} differential equation
  \eqref{eq:diffeq_z} at $y_0$, by Dubrovin's Reconstruction Theorem
  we have a natural isomorphism of Frobenius manifolds $\Upsilon \co
  \big(\HorbX,\sigma_0\big) \cong \big(\HY,\tau_0\big)$.  Forgetting
  the opposite subspace, we conclude that there is an isomorphism of
  the underlying \VHSns s, i.e. that the isomorphism
  (\ref{eq:VHSisom_restricted}) extends to open neighbourhoods of
  $\sigma_0\in \HorbX$ and $\tau_0\in \HY$.  The moving subspace
  realizations of the two big A-model \VHSns s are therefore related
  by a constant $\Czzinv$-linear transformation.  Equation
  (\ref{eq:UmatchesE_small}) shows that this transformation is
  $\ov\U$.
\end{proof} 

\begin{rem}
  When reconstructing big quantum cohomology from small quantum
  cohomology, we could use quantum $H^2$-generation in place of
  Dubrovin's Reconstruction Theorem.  In fact, Dubrovin Reconstruction
  is a special case of quantum $H^2$-generation where the product
  $E\bullet_\tau$ of the Euler vector field generates the total
  cohomology.  In our case, orbifold cohomology $H^\bullet_{\rm
    orb}(\PP(1,1,1,3))$ is not generated by $H^2(\PP(1,1,1,3))$, but
  \emph{quantum} orbifold cohomology is generated by
  $H^2(\PP(1,1,1,3))$.  Reconstruction theorems of Hertling--Manin
  \cite{Hertling--Manin:unfolding}, Iritani \cite[Remark
  4.10]{Iritani:gen}, and Rose \cite{Rose} are also applicable here.
  These are generalizations of the First Reconstruction Theorem of
  Kontsevich-Manin \cite{Kontsevich--Manin:GW}, where  classical
  $H^2$-generation is assumed.
\end{rem}

\subsection{The Flat Co-ordinates Are Different}

We can see the difference between the big quantum cohomology Frobenius
manifolds for $\FF_3$ and $\PP(1,1,1,3)$ more explicitly as follows.
The vectors $I_{\FF_3}(y,-z)$ and $\U^{-1}I_{\FF_3}(y,-z)$
are on the cones $\cL_{\FF_3}$ and $\cL_{\PP(1,1,1,3)}$ respectively.
They expand as:
\begin{align}
  \label{eq:compareI}
  I_{\FF_3}(y,-z)&=-z+ \tau_1 p_1+ \tau_2 p_2 - 
  \left(-\frac{1}{3}\parfrac{F_0}{\tau_1} \ap_1^2 + 
    \parfrac{F_0}{\tau_2} p_1^2\right)
  \frac{1}{z}+O(z^{-2}) \\ \notag
  \U^{-1}\big(I_{\FF_3}(y,-z)\big) &=-z+ \atau_1 \fun_{2 \over 3}+ \atau_2 (3p) - 
  \left(3\parfrac{F_0^{\rm orb}}{\atau_1}\fun_{1 \over 3} + 
    \parfrac{F_0^{\rm orb}}{\atau_2} p^2\right)
  \frac{1}{z}+O(z^{-2}) 
\end{align}
where $(\tau_1,\tau_2)$ and $(\atau_1,\atau_2)$ are flat coordinates
on $\tcM_{\rm B}$ associated with the quantum cohomologies and $F_0,
F_0^{\rm orb}$ are the genus-zero Gromov--Witten potentials for
$\FF_3$ and $\PP(1,1,1,3)$.  Combining \eqref{eq:matrixofU} and
\eqref{eq:compareI}, we see that 
\begin{align}
  \label{eq:flatmatch1}
  \tau_1= -\frac{2\sqrt{3}\pi}{3\Gamma(\frac{2}{3})^3}\atau_1 + 
  \frac{2\sqrt{3}\pi}{\Gamma(\frac{1}{3})^3}\parfrac{F^{\rm orb}_0}{\atau_1}, 
  &&
  \tau_2+\frac{1}{3}\tau_1 =\atau_2.  
\end{align}
From \eqref{eq:IF3continued} we find 
\begin{align*}
  \atau_1= \sum_{n=0}^\infty (-1)^n 
  \frac{\prod_{k=0}^{n-1} (k+\frac{1}{3})^3}{(3n+1)!} \aq_1^{3n+1}, && 
  3\parfrac{F_0^{\rm orb}}{\atau_1} 
  =\sum_{n=0}^\infty (-1)^n 
  \frac{\prod_{k=0}^{n-1} (k+\frac{2}{3})^3}{(3n+2)!}
  \aq_1^{3n+2}  
\end{align*}
and thus:
\begin{equation}
  \label{eq:flatmatch2}
  3\parfrac{F_0^{\rm orb}}{\atau_1}= \frac{1}{2} \atau_1^2 
  -\frac{1}{3^2\cdot 5!} \atau_1^5 
  + \frac{1}{3\cdot 8!} \atau_1^8 - \frac{1093}{3^5\cdot 11!} \atau_1^{11}+ 
  \cdots 
\end{equation}
Combining \eqref{eq:flatmatch1} and \eqref{eq:flatmatch2} shows that
the flat co-ordinate systems $(\tau_1,\tau_2)$ and $(\atau_1,\atau_2)$
on $\cM_{\rm B}$ are different.

\subsection{The Proof of Theorem~\ref{thm:P1113}}
\label{sec:proofP1113}

We recommend that at this point the reader reviews the strategy
described in Section~\ref{sec:VHSmirrorsym}.   

The small quantum cohomology locus $\cM_{\PP(1,1,1,3)}$ of
$\PP(1,1,1,3)$ is given by
\[
\{\aq_1=0\} = \{\atau_1=0\} = \{\tau_1 = 0\}
\] 
Along this locus, the two flat co-ordinates coincide:
$\tau_2=\atau_2$.  As $q_1=e^{\tau_1}$, $q_2 = e^{\tau_2}$, and $q =
e^{3 \atau_2}$ we have
\begin{align*}
  q_1 = 1, && q_2 = \sqrt[3]{q}.
\end{align*}
We calculate the identification (along this locus) between the quantum 
cohomology algebras of $\PP(1,1,1,3)$ and $\FF_3$ by first finding
differential operators which represent our chosen basis for
$H^\bullet_{\text{\rm orb}}\big(\PP(1,1,1,3);\CC\big)$ through
derivatives of $I_{\PP(1,1,1,3)}$, then commuting these operators past the
symplectic transformation $\U$ in the equality
$\U\big(I_{\PP(1,1,1,3)}(y,-z)\big) = I_{\FF_3}(0,y^{1/3},-z)$, and
finally comparing the resulting derivatives of $I_{\FF_3}$ with our chosen
basis for $H^\bullet(\FF_3;\CC)$.  In detail, this goes as follows.

The
opposite subspace $\cH^-_{c_2}$ at the cusp $c_2$ for $\PP(1,1,1,3)$
determines a trivialization of $\cE$.  Define differential operators
$P_i(z \partial)$ by
\begin{align*}
&P_0 = 1, && P_1 =  z\partial, && P_2 = (z\partial)^2, \\
&P_3 = (z\partial)^3, && P_4 =  y^{-1/3} 3(z\partial)^4, 
&&P_5 = y^{-1/3} 3z\partial (y^{-1/3} 3(z\partial)^4)   
\end{align*}
where $\partial=y\parfrac{}{y}$.  The sections 
$e_i=P_i(\nabla^z_{\partial})[\exp(W_y)\,\omega_y]$ form a frame of
$\cE$ which is constant with respect to this trivialization. In fact,
from
\begin{equation}
\label{eq:expansionofIP1113}
I_{\PP(1,1,1,3)} =y^{p/z} 
\left(\fun_0+\frac{27y^{1/3}}{z^4} \fun_{1 \over 3}+ 
\frac{27 y^{2/3}}{16 z^5} \fun_{2 \over 3} +\frac{y}{6z^6} \fun_0
+O\big(z^{-7}\big)\right)
\end{equation}
one finds that:
\begin{align*}
& P_i I_{\PP(1,1,1,3)} = p^{i}+O\big(z^{-1}\big) && 0\le i\le 3 \\
& P_4 I_{\PP(1,1,1,3)} = \fun_{1 \over 3} +O\big(z^{-1}\big) \\
& P_5 I_{\PP(1,1,1,3)} =\fun_{2 \over 3} + O\big(z^{-1}\big) 
\end{align*}
and therefore that the differential operators $P_0$, \ldots, $P_5$
correspond to the basis $1,p,p^2,p^3,\fun_{1 \over 3},\fun_{2 \over
  3}$ for the quantum cohomology algebra of $\PP(1,1,1,3)$.  The
matrix of quantum multiplication $p \circ_q$ can be obtained as the
connection matrix of $\nabla^z_\partial$ with respect to the frame
$\{e_i\}$:
\begin{equation}
\label{eq:connectionp1113}
\nabla^z_\partial = z \, \partial +  
\begin{pmatrix}
0 & 0 & 0 & 0& 0& \frac{1}{3}y^{1/3} \\ 
1 & 0 & 0 & 0& 0&  0 \\
0 & 1 & 0 & 0& 0& 0 \\
0 & 0 & 1 & 0& 0& 0 \\
0 & 0 & 0 & \frac{1}{3}y^{1/3} & 0 & 0 \\
0 & 0 & 0 & 0 & \frac{1}{3}y^{1/3} & 0 
\end{pmatrix}
\end{equation}
(Recall that the mirror map for $\PP(1,1,1,3)$ is trivial, so $q=y$.)

On the other hand, the sections $e_i$ are not constant with respect
to the trivialization associated with $\cH^-_{c_1}$.  By using $P_i
I_{\FF_3}|_{\aq_1=0} = \Ubar P_i I_{\PP(1,1,1,3)}$ and the expansion
\eqref{eq:expansionofIP1113}, one finds that:
\begin{align*}
& P_i I_{\FF_3}|_{\aq_1=0} = \(\frac{p_2}{3}\)^{i} + 
O\big(z^{-1}\big) && 0\le i\le 2 \\
& P_3 I_{\FF_3}|_{\aq_1=0} = \frac{p_2^3}{27} - 
\sqrt{3}\beta_1 y^{1/3} \ap_1 + O\big(z^{-1}\big) \\
&P_4 I_{\FF_3}|_{\aq_1=0} = -\sqrt{3}\beta_1 \ap_1 z 
+ \frac{\pi}{3}\beta_1\ap_1^2 + O\big(z^{-1}\big) \\
&P_5 I_{\FF_3}|_{\aq_1=0} = \sqrt{3}\beta_2\ap_1 
+ \frac{1}{z} 
\(\frac{\pi}{3}\beta_2\ap_1^2 + y^{1/3}\) +O\big(z^{-2}\big)
\end{align*}
where $\beta_i = \frac{2\pi}{9\Gamma(i/3)^3}$.  From this, we see that
the frame of $\cE$ given by
\begin{gather*}
e_0, \quad e_1, \quad e_2, \quad e_3, \quad 
e_4+z\frac{\beta_1}{\beta_2} e_5, \quad 
e_5
\end{gather*}
corresponds to the non-constant basis 
\begin{gather*}
1, \quad \frac{p_2}{3}, \quad \frac{p_2^2}{9}, \quad 
\frac{p_2^3}{27} - \sqrt{3}\beta_1 y^{1/3}\ap_1, \quad 
\frac{\beta_1}{\beta_2} y^{1/3}+ 
\frac{2\pi}{3}\beta_1 \ap_1^2, 
\quad \sqrt{3}\beta_2\ap_1
\end{gather*}
for the quantum cohomology algebra of $\FF_3$.  The Dubrovin
connection for the quantum cohomology of $\FF_3$ can be obtained from
the connection \eqref{eq:connectionp1113} by the gauge transformation
$\Theta(y,z)\colon H^\bullet_{\text{\rm orb}}\big(\PP(1,1,1,3)\big) \rightarrow
H^\bullet(\FF_3)$ given by:
\begin{align*}
&\Theta(y,z)(p^i) = \(\frac{p_2}{3}\)^{i} &&0\le i\le 2 \\ 
&\Theta(y,z)(p^3) = \frac{p_2^3}{27} - 
\sqrt{3}\beta_1 y^{1/3}\ap_1 \\
&\Theta(y,z)(\fun_{1 \over 3})= \frac{\beta_1}{\beta_2} y^{1/3}+ 
\frac{2\pi}{3}\beta_1 \ap_1^2 - 
z \sqrt{3}\beta_1 \\
&\Theta(y,z) (\fun_{2 \over 3}) = \sqrt{3}\beta_2\ap_1
\end{align*} 
Therefore the quantum product by $p_2/3$ and by $p$, at $q \in
\cM_{\PP(1,1,1,3)}$, are related by conjugation by
$\Theta(q):=\Theta(q,0)$:
\[
\(\frac{p_2}{3}\circ_{(q_1,q_2)}\)\Big|_{(q_1,q_2) = (1,\sqrt[3]{q})} = \Theta(q) \(p\circ_q\) \Theta(q)^{-1}
\]
It is easy to check that $\Theta(q)$ preserves the (orbifold)
Poincar\'{e} pairing and grading.  Because $\Theta(q)$ preserves the
unit and $p\circ_q$ generates the small quantum cohomology algebra,
$\Theta(q)$ is an algebra isomorphism.  This proves
Theorem~\ref{thm:P1113}. \hfill \qedsymbol

\begin{rem}
  The symplectic transformation $\U$ does not induce an isomorphism
  between the Frobenius manifolds associated to the big quantum
  cohomologies of $\FF_3$ and $\PP(1,1,1,3)$ but it does induce an
  isomorphism between the corresponding $F$-manifolds.
\end{rem}

\begin{rem}
  The basis change operator $\Theta(q)$ becomes a ring isomorphism
  because it preserves the unit.  If we have a miniversal
  extended B-model moduli space $\hcM_{\rm B}$ corresponding to big
  quantum cohomology and a mirror $D$-module $\cE$ on it --- in fact
  we can reconstruct these from the small data --- then the basis
  change operator $\Theta(q,z)$ between two flat frames can be
  extended to $q\in \hcM_{\rm B}$.  The operator $\Theta(q)$ will not
  necessarily preserve the unit outside the original B-model moduli
  space $\cM_{\rm B} \subset \hcM_{\rm B}$, and so will not in general
  be a ring isomorphism there.  A ring isomorphism over the whole of
  $\hcM_{\rm B}$ is given by $v \longmapsto \Theta(q) (v \circ_q)
  \Theta(q)^{-1} \fun$, but outside of $\cM_{\rm B}$ this will not in
  general preserve the (orbifold) Poincar\'{e} pairing.
\end{rem} 

\section{Example: $\FF_2$ and $\PP(1,1,2)$}
\label{sec:P112}

We now consider the examples $\cX = \PP(1,1,2)$ and $Y = \FF_2$,
proving Theorem~\ref{thm:P112} and Conjecture~\ref{conj}.  The
argument is entirely parallel to that in Section~\ref{sec:P1113} and
so we omit many details.  The only significant difference is in the
conclusion which we draw.  Since the opposite subspaces associated to
the cusps for $\PP(1,1,2)$ and $\FF_2$ agree under parallel transport
--- or, more concretely, because the symplectic transformation
$\U:\cH_{\PP(1,1,2)} \to \cH_{\FF_2}$ in
Proposition~\ref{pro:symplectictransformationP112} maps
$\cH_{\PP(1,1,2)}^-$ to $\cH_{\FF_2}^-$ --- it follows that the flat
structures associated to $\PP(1,1,2)$ and $\FF_2$ agree under analytic
continuation.  This implies that the big quantum cohomology Frobenius
manifolds for $\PP(1,1,2)$ and $\FF_2$ become isomorphic after
analytic continuation, and hence that the original form of the
Bryan--Graber Conjecture holds.
See Appendix B for a more elementary 
proof of the crepant resolution conjecture in this surface case. 

\subsection{The Landau--Ginzburg Mirror to $\FF_2$}

The surface ${\FF_2}$ is a GIT quotient of $\CC^4$ by $(\Cstar)^2$
where $(\Cstar)^2$ acts via the inclusion:
\begin{align*}
(\Cstar)^2 \hookrightarrow (\Cstar)^4&&
(s,t) \mapsto (s,s,s^{-2}t,t)
\end{align*}
The mirror family $\pi:Z \to \cM_{\FF_2}$ is given by restricting
the dual of this inclusion
\begin{align*}
\pi\colon (\Cstar)^4 & \longrightarrow (\Cstar)^2 \\
(w_1,\dots,w_4) &\longmapsto (w_1w_2w_3^{-2},w_3w_4) 
\end{align*}
to the open subset $\cM_{\FF_2} \subset (\Cstar)^2$ defined by:
\[
\cM_{\FF_2} = \Big\{(y_1,y_2) \in (\Cstar)^2 : \textstyle y_1 \neq {1
  \over 4} \Big\}
\]  
The superpotential $W$ is
\[
W=w_1+w_2+w_3+w_4
\]
and the holomorphic volume form $\omega_y$ on the fiber
$Z_y=\pi^{-1}(y_1,y_2)$ is:
\[
\omega_y =\frac{d\log w_1\wedge \cdots \wedge d\log w_4}
{d\log y_1\wedge d\log y_2}
\] 
We deleted the locus $y_1 = {1 \over 4}$ from $\cM_{\FF_2}$ to ensure
that Assumptions~\ref{assumptions} hold.  The oscillating integrals
\eqref{eq:oscint} here satisfy the Picard--Fuchs equations:
\begin{equation}
  \label{eq:PFF2}
  \begin{split}
    D_2 (D_2 - 2 D_1) \,f &= y_2 \, f\\
    D_1^2 D_2 (D_2 - z) \, f&= y_1 y_2^2 \, f \\
    D_1^2 D_2 \, f &= y_1 y_2 (D_2 - 2 D_1) \, f \\
    D_1^2 \, f &= y_1 (D_2 - 2 D_1)(D_2 - 2 D_1 - z) \, f
  \end{split}
\end{equation}
where $D_1 = z y_1 {\partial \over \partial y_1}$ and $D_2 = z y_2
{\partial \over \partial y_2}$.  

\subsection{Mirror Symmetry for $\FF_2$}

Givental's mirror theorem \cite{Givental:toric}*{Theorem~0.1} implies
that the $I$\nobreakdash-function
\begin{equation}
  \label{eq:IF2}
  I_{\FF_2}(y_1,y_2,z) = z \, 
  \sum_{k,l \geq 0}
  { y_1^{k+p_1/z} y_2^{l+p_2/z}  \prod_{m=-\infty}^0 (p_2 - 2 p_1 + m z) 
    \over 
    \prod_{m=1}^k (p_1 + m z)^2
    \prod_{m=1}^l (p_2 + m z)
    \prod_{m=-\infty}^{l-2k} (p_2 - 2 p_1 + m z)}
\end{equation}
where $p_1$, $p_2$ is the basis of $H^2(\FF_2)$ defined in Section
\ref{sec:cohomology}, coincides with the small $J$-function
\eqref{eq:smallJF2} after a change of variables $(y_1,y_2) \mapsto
(q_1,q_2)$:
\[
J_{\FF_2}(q_1,q_2,z) = I_{\FF_2}(y_1,y_2,z)
\]
The components of $I_{\FF_2}(y_1,y_2,z)$ form another basis of
solutions to the Picard--Fuchs system \eqref{eq:PFF2}.  As we
explained in Section~\ref{sec:VHSmirrorsym}, this proves
Conjecture~\ref{conj:mirrorsym} for $\FF_2$.

As before, we can can read off the mirror map $(y_1,y_2) \mapsto
(q_1,q_2)$ by expanding the $I$-function as a Laurent series in
$z^{-1}$.  This gives:
\begin{equation}
  \label{eq:F2mirrormap}
  \begin{split}
    q_1 &= y_1 \exp\(2 \sum_{k\ge 1}  {(2k-1)! \over \(k!\)^2} y_1^k\)
    = {4y_1 \over \big(1 + \sqrt{1-4y_1}\big)^2}\\
    q_2 &= y_2 \exp\({-\sum_{k\ge 1} {(2k-1)! \over \(k!\)^2} y_1^k}\)
    = {y_2 \big(1+\sqrt{1-4y_1}\big) \over 2}
  \end{split}
\end{equation}

The mirror map identifies a neighbourhood of $y_1=y_2=0$ in
$\cM_{\FF_2}$ with the K\"{a}hler moduli space $\cM_{\rm A}$ of $\FF_2$;
here once again $(q_1,q_2)$ are co-ordinates on $\cM_{\rm A}$.  This
identification matches up the Jacobi ring $J(W_y)$ of $W_y$ with the
small quantum cohomology algebra of $\FF_2$ at $(q_1,q_2)$ and the
residue pairing with the Poincar\'e pairing.

\subsection{The Landau--Ginzburg Mirror to $\PP(1,1,2)$}

The mirror family $\pi$ is
\begin{align*}
\pi\colon Z=(\Cstar)^3 & \longrightarrow \cM_{\PP(1,1,2)}=(\Cstar) \\
(w_1,w_2,w_4) &\longmapsto w_1w_2w_4^2
\end{align*}
the superpotential $W$ is
\[
W=w_1+w_2+w_4
\]
and the holomorphic volume form $\omega_y$ on the fiber
$Z_y=\pi^{-1}(y)$ is:
\[
\omega_y =\frac{d\log w_1\wedge d\log w_2\wedge d\log w_4}
{d\log y}
\] 
The oscillating integrals \eqref{eq:oscint} satisfy the Picard--Fuchs
equation
\begin{equation}
  \label{eq:PFP112}
  D^2 (2D) (2D - z)\, f = y\, f.
\end{equation}
where $D = z y {\partial \over \partial y}$.

\subsection{Mirror Symmetry for $\PP(1,1,2)$}

A theorem of Coates--Corti--Lee--Tseng \cite{CCLT}*{Theorem~1.7} shows
that the $I$-function
\begin{equation}
  \label{eq:IP112}
  I_{\PP(1,1,2)}(y,z) = z \, y^{p/z} \sum_{\substack{d \colon 2 d \in \ZZ
      \\ d \geq 0}} { y^d \over \prod_{\substack{b \colon \fr{b} = \fr{d} \\ 0 < b \leq d}} (p + bz)^2
    \prod_{\substack{b \colon \fr{b} = 0 \\ 0 < b \leq 2d}} (2p + bz)}
  \fun_{\fr{d}}
\end{equation}
coincides with the small $J$-function \eqref{eq:smallJP112} after
the (trivial) change of variables $q=y$:
\[
J_{\PP(1,1,2)}(q,z) = I_{\PP(1,1,2)}(y,z)
\]
The components of $I_{\PP(1,1,2)}(y,z)$ give another basis of
solutions to the Picard--Fuchs equation \eqref{eq:PFP112}.  As before
this proves Conjecture~\ref{conj:mirrorsym} for $\PP(1,1,2)$.  The
moduli space $\cM_{\PP(1,1,2)}$ with co-ordinate $y$ is identified via
the map $q=y$ with the K\"ahler moduli space $\cM_{\rm A}$ of $\PP(1,1,2)$;
here $q$ is again a co-ordinate on $\cM_{\rm A}$.  This identification
matches the Jacobi ring $J(W_q)$ with the small quantum orbifold
cohomology algebra of $\PP(1,1,2)$ at $q$, and the residue pairing
with the Poincar\'e pairing.

\subsection{Constructing the B-model \VHSns}

As in Section~\ref{sec:BVSHSP1113}, we extend the Landau--Ginzburg
mirror of $\FF_2$ to a Landau--Ginzburg model with a larger base
defined in terms of the secondary fan for $\FF_2$ (Figure
\ref{fig:secondaryF_2}).

\begin{figure}[ht!]
  \begin{center}
    \begin{picture}(180,50)(40,0)
      \put(150,5){\vector(1,0){40}}
      \put(150,5){\vector(-2,1){80}}
      \put(150,5){\vector(0,1){40}}
      \multiput(70,5)(40,0){4}{\makebox(0,0){$\cdot$}}
      \multiput(70,45)(40,0){4}{\makebox(0,0){$\cdot$}}
      \put(195,0){$p_1$}
      \put(58,40){$\ap_1$}
      \put(150,50){\makebox(0,0){$\ap_2=p_2$}}
    \end{picture}
    \caption{The secondary fan for ${\FF_2}$.}
    \label{fig:secondaryF_2}
  \end{center}
\end{figure}
Take $w_1,w_4$ as co-ordinates on the fiber $Z_y$, so that:
\begin{equation}
  \label{eq:LGF2patch1}
  W_y=w_1+\frac{y_1y_2^2}{w_1w_4^2}+\frac{y_2}{w_4}+w_4
\end{equation} 
Let $\ccM$ be the toric orbifold associated to the secondary fan for
$\FF_2$.  There are two distinguished co-ordinate patches on $\ccM$,
one for each maximal cone in the secondary fan.  Let $(y_1, y_2)$ be the
co-ordinates dual to $(p_1, p_2)$ and $(\aq_1, \aq_2)$ be the co-ordinates
dual to $(\ap_1, \ap_2)$ (see Figure~\ref{fig:secondaryF_2}).  As
\begin{align*}
\aq_1=y_1^{-1/2}&& \aq_2= y_1^{1/2} y_2
\end{align*}
we see that $(\aq_1,\aq_2)$ is a uniformizing system near a $\ZZ/2\ZZ$
orbifold point at $(\aq_1,\aq_2)=0$.  In the co-ordinates
$(\aq_1,\aq_2)$ we have
\begin{equation}
  \label{eq:LGF2patch2}
  W_y=w_1+\frac{\aq_2^2}{w_1w_4^2}+\frac{\aq_1\aq_2}{w_4}+w_4 
\end{equation}
and so we can extend the family of tori $\pi$ and the superpotential
$W_y$ across the locus $\{\aq_1=0\}$.  Here we see
\begin{align*}
W=w_1+w_2+w_4, && w_1w_2w_4^2=\aq_2^2 
\end{align*}
which is the Landau--Ginzburg mirror to $\PP(1,1,2)$: the locus
$\aq_1=0$, $\aq_2\neq 0$ in $\ccM$ is identified with
$\cM_{\PP(1,1,2)}$ via the map $y = \aq_2^2$.

The \emph{B-model moduli space} $\cM_{\rm B}$ here, which is the base of
our extended Landau--Ginzburg model, is obtained from $\ccM$ by
deleting the closures of the loci $ \{y_1y_2=0\}$ and $\{y_1=1/4\}$.
Equations \eqref{eq:LGF2patch1} and \eqref{eq:LGF2patch2} together
define a Landau--Ginzburg model over $\cM_{\rm B}$ which contains the
Landau--Ginzburg mirrors for $\FF_2$ and for $\PP(1,1,2)$ as subsets.
The limit points $y_1=y_2=0$ and $\aq_1=\aq_2=0$ of $\cM_{\rm B}$ are called
the \emph{large radius limit points} or \emph{cusps} corresponding
respectively to $\FF_2$ and to $\PP(1,1,2)$.  Let $\cE$ denote the
B-model \VHS with base $\cM_{\rm B}$ defined by the Landau--Ginzburg model
just described and equipped  with the grading operator: 
\begin{align*}
\Gr \big[f(x,z)e^{W_y/z}\,\omega_y\big] = 
\left[\(2 z\partial_z + 2\sum_{i=1}^4 w_i \partial_{w_i}\) f(x,z)
e^{W_y/z} \, \omega_y\right]
\end{align*}
The Euler field and dimension here are
\begin{align*}
E=2y_2\parfrac{}{y_2}= 2\aq_2\parfrac{}{\aq_2}, && 
D=\dim \FF_2 =2.
\end{align*}
\begin{rem}
  The superpotentials \eqref{eq:LGF2patch1} and
  \eqref{eq:LGF2patch2} have isolated non-degenerate critical points;
  this implies that the small quantum cohomology algebra of $\FF_2$
  and the small quantum orbifold cohomology algebra of $\PP(1,1,2)$
  are semisimple.
\end{rem}

\subsection{An Opposite Subspace At Each Cusp}

As before, the $I$-functions \eqref{eq:IF2}, \eqref{eq:IP112}
and the isomorphism \eqref{eq:Dmodulegeneratedbyoscint} define
$D$-module homomorphisms:
\begin{equation}
  \label{eq:ImapsF2}
  \begin{split}
    & \mapI_{\FF_2} \colon \cE \longrightarrow H^\bullet(\FF_2)\otimes 
    \cOzzinv{\tcM_{\rm B}}, \\
    & \mapI_{\PP(1,1,2)} \colon \cE|_{\cM_{\PP(1,1,2)}}
    \longrightarrow H^\bullet_{\text{\rm orb}}\big(\PP(1,1,2)\big)\otimes 
    \cOzzinv{\tcM_{\PP(1,1,2)}}
  \end{split}
\end{equation}
by sending $P(y,z\partial,z)\in \cE \cong \cOz{\cM_{\rm B}}\langle z\cT_{\cM_{\rm
    B}}\rangle$ to $z^{-1} P I_{\FF_2}$ or to $z^{-1} P
I_{\PP(1,1,2)}$.  The maps \eqref{eq:ImapsF2} send $\nabla^z$-parallel
sections of $\cE$ to constant sections, and thus identify the space
$\cH$ of flat sections of $\cE$ with Givental's symplectic vector
space:
\begin{equation*}
  \begin{aligned}
    & \mapI_{\FF_2} \colon \cH \cong 
    H^\bullet(\FF_2)\otimes \Czzinv =\cH_{\FF_2}|_{Q_i=1} \\
    & \mapI_{\PP(1,1,2)} \colon \cH \cong H^\bullet_{\text{\rm orb}}\big(\PP(1,1,2)\big) 
    \otimes \Czzinv=\cH_{\PP(1,1,2)}|_{Q=1}
  \end{aligned}
\end{equation*}
Here $\mapI_{\FF_2}$, which is \emph{a priori} defined only in a small
neighbourhood of the cusp for $\FF_2$, is extended to the whole of
$\tcM_{\rm B}$ by analytic continuation.  

Let $c_1$ and $c_2$ denote the cusps of $\cM_{\rm B}$ corresponding
respectively to $\FF_2$ and $\PP(1,1,2)$. We define opposite subspaces
$\cH^-_{c_1}$ and $\cH^-_{c_2}$ of $\cH$ by:
\begin{align*}
  \mapI_{\FF_2}\big(\cH^-_{c_1}\big) =\cH_{\FF_2}^-|_{Q_i=1} && 
  \mapI_{\PP(1,1,2)}\big(\cH^-_{c_2}\big) =  \cH_{\PP(1,1,2)}^-|_{Q=1}
\end{align*}
These opposite subspaces are uniquely characterized by monodromy and
homogeneity properties, as in Proposition~\ref{thm:uniqueness}, but we
will not pursue this here.

\subsection{The Polarizations Match}
\label{sec:differenceofpolarization_F2}

Define $\Ubar \colon \cH_{\PP(1,1,2)}|_{Q=1} 
\rightarrow \cH_{\FF_2}|_{Q_i=1}$ to be
the composition 
\[
\begin{CD}
\cH_{\PP(1,1,2)}|_{Q=1} 
 @>{(\mapI_{\PP(1,1,2)})^{-1}}>> \cH @>{\mapI_{\FF_2}}>> 
\cH_{\FF_2}|_{Q_i=1}   
\end{CD}
\]
and let $\U: \cH_{\PP(1,1,2)}|_{Q=1} 
\rightarrow \cH_{\FF_2}|_{Q_i=1}$ be $\Ubar$
followed by changing the sign of $z$.  Arguing as in the proof of
Proposition~\ref{pro:symplectictransformationP1113} shows:

\begin{pro}
\label{pro:symplectictransformationP112}
The matrix of $\U$ with respect to the bases for $H^\bullet_{\rm
  orb}(\PP(1,1,2))$ and $H^\bullet(\FF_2)$ defined in
Section~\ref{sec:cohomology} is:
\begin{equation}
  \label{eq:matrixofUF2}
  \begin{pmatrix}
    1 & 0 & 0 & 0 \\
    \frac{ \pi \tti }{z} & 0 & 0 & \tti \\
    -\frac{\pi \tti }{2 z} & \frac{1}{2} & 0 & -\frac{\tti}{2} \\
    \frac{\pi ^2}{4 z^2} & 0 & \frac{1}{2} & \frac{\pi }{2 z}
  \end{pmatrix}
\end{equation} 
The linear transformation $\U$ preserves the grading, the symplectic
forms, and the standard opposite subspaces. \hfill \qedsymbol
\end{pro} 

Note that here $\U$ takes the form $\exp\Big({-\frac{\pi \tti
  (p_2-2p_1)}{2z}}\Big) \circ \U|_{z=\infty}$.

\subsection{The Proof of Conjecture~\ref{conj}}

By applying Dubrovin's Reconstruction Theorem, as in the proof of
Theorem~\ref{thm:transformationofLagrangiancones}, we deduce 
Conjecture~\ref{conj} and \ref{conj:preciseCRC} 
for $\cX = \PP(1,1,2)$ and $Y = \FF_2$: 

\begin{thm}
\label{thm:UP112}
Let $\E^{\FF_2}_\tau \subset \left.\cH_{\FF_2}\right|_{Q_i=1}$ and
$\E^{\PP(1,1,2)}_\tau \subset \left.\cH_{\PP(1,1,2)}\right|_{Q=1}$ 
be the moving subspace realizations \eqref{eq:moving_A} 
of the analytic 
big A-model \VHSns s of $\FF_2$ and $\PP(1,1,2)$. 
Then there exists a map $\Upsilon$ from an open subset 
of $H_{\rm orb}^\bullet(\PP(1,1,2))$ 
to an open subset of $H^\bullet(\FF_2)$ 
such that, after analytic continuation if necessary,  
\[
\ov\U \big(\E^{\PP(1,1,2)}_\tau \big) 
= \E^{\FF_2}_{\Upsilon(\tau)}
\]
where $\U$ is the symplectic transformation from
Proposition~\ref{pro:symplectictransformationP112}. \hfill \qedsymbol
\end{thm}


\subsection{The Proof of Theorem~\ref{thm:P112}}
\label{sec:P112thm}
It remains only to prove Theorem~\ref{thm:P112}.  But since the
transformation $\U$ maps:
\begin{itemize}
\item the big A-model \VHS for $\PP(1,1,2)$ to the big A-model
  \VHS for $\FF_2$
\item the standard opposite subspace $\cH_{\PP(1,1,2)}^-$ to the
  standard opposite subspace $\cH_{\FF_2}^-$
\item the dilaton shift $1 \in z\cH_{\PP(1,1,2)}^-/\cH_{\PP(1,1,2)}^-$
  to the dilaton shift $1 \in z\cH_{\FF_2}^-/\cH_{\FF_2}^-$
\end{itemize}
it follows immediately that $\U$ induces an isomorphism between the
big quantum cohomology Frobenius manifolds associated to $\cX =
\PP(1,1,2)$ and $Y = \FF_2$.  To compute this isomorphism explicitly,
consider the discussion before
Proposition~\ref{pro:VHStoFrobeniusmanifold}.  The underlying linear
isomorphism $\Theta:\HorbX \to \HY$ here is the isomorphism
$z\cHX^-/\cHX^- \cong z\cHY^-/\cHY^-$ induced by $\U$, so $\Theta =
\U|_{z = \infty}$:
\begin{equation}
  \label{eq:mapforP112}
  \begin{alignedat}{3}
    &\Theta(\fun_0) = 1 && \qquad \qquad \qquad&
    &\Theta(p) =  \textstyle {p_2 \over 2} \\
    &\Theta(p^2) = \textstyle \big({p_2 \over 2}\big)^2 & &&
    &\Theta(\fun_{1/2}) = \textstyle {-{\tti \over 2}}\big(p_2 - 2p_1\big)
  \end{alignedat}
\end{equation}
The map $\Theta$ evidently preserves the Poincar\'e pairings.  It gives an
isomorphism of algebras between
\begin{align*}
  \Big(\HorbX,\bullet_\tau\Big) && \text{and} && \Big(\HY,\bullet_{f(\tau)}\Big)
\end{align*}
where we can read off the affine-linear identification of flat
co-ordinates $\tau \mapsto f(\tau)$ from the big $J$-functions:
$\U\big(J_\cX(\tau,-z)\big) = J_Y\big(f(\tau),-z\big)$, and so
\begin{equation}
  \label{eq:flatshift}
  f(\tau) = \Theta(\tau) + \frac{\pi\tti}{2}(p_2 - 2 p_1).
\end{equation}
Here we are considering the analytic version of Givental's formalism,
with the Novikov variables $Q$, $Q_1$, and~$Q_2$ set to $1$.  Putting back
the Novikov variables using the Divisor Equation (see
Remark~\ref{rem:qandQ} above)  
one finds that one can absorb the shift of the origin in \eqref{eq:flatshift}
into the specialization of quantum parameters:
\begin{align*}
  Q_1 &= -1, && Q_2 = \tti Q^{1/2}.
\end{align*}
Theorem~\ref{thm:P112} is proved.\hfill \qedsymbol

\section{A Crepant Resolution Conjecture}
\label{sec:CRC}

In this final section we formulate our version of the Crepant
Resolution Conjecture.  This is a more precise version of
Conjecture~\ref{conj} from the Introduction.  We discuss its
relationship with theorems of Lupercio--Poddar and Yasuda, and show
that under a Hard Lefschetz condition it implies the original form of
the Bryan--Graber Conjecture.  We also indicate several aspects of the
story which remain to be explored.

\begin{conj}
\label{conj:preciseCRC}
Let $\cX$ be an orbifold with projective
coarse moduli space $X$ and let $\pi\colon Y \to X$ be a crepant
resolution.  Suppose that the big quantum products $\bullet_\tau$ in
\eqref{eq:bigQC} for $\cX$ and $Y$ are convergent as functions of
$\tau$ and $Q$, so that the analytic big A-model \VHS of $\cX$
(respectively of $Y$) with Novikov variables specialized to $1$ is
well-defined over an open subset of $\HorbX$ (respectively of $\HY$);
see Section \ref{subsubsec:analyticbig}.

Let $\E^\cX_\tau\subset \cH_\cX|_{Q_i=1}$, 
$\E^Y_\tau\subset \cH_Y|_{Q_i=1}$ be 
the moving subspace realizations \eqref{eq:moving_A} 
of the analytic big A-model \VHSns s of $\cX$ and $Y$ respectively. 
Define the limiting Hodge structure  
$\E^\cX_{{\rm lim},\tau}$ associated to $\tau\in \HorbX$ by:
  \[
  \E_{{\rm lim},\tau}^\cX 
  = \lim_{\sigma\to \text{l.r.l.}} e^{-\sigma/z} \, \E^\cX_{\tau+\sigma}  
  \]
Here $\sigma$ moves in $H^2(\cX;\CC)$ and $\lim_{\sigma\to
\text{l.r.l.}}$ denotes the large radius limit:
$\Re(\CHpair{\sigma}{d})\rightarrow -\infty$ for all $d\in
\Eff(\cX)$.

  There is a symplectic transformation 
  $\ov\U \colon \cHX|_{Q_i=1} 
  \rightarrow \cHY|_{Q_i=1}$ and a map 
  $\Upsilon$ from an open subset of $\HorbX$ to 
  an open subset of $\HY$ such that, 
  after analytic continuation if necessary, 
  \[
  \ov\U \big(\E_\tau^\cX \big) = \E^Y_{\Upsilon(\tau)}
  \]
  and that:
  \begin{equation}
    \label{eq:CRCconditions}
    \begin{minipage}{0.87\linewidth}
      \begin{itemize}
  \item[(a)] $\ov\U$ is degree-preserving and $\Czzinv$-linear
  \item[(b)] $\ov\U (\rho\cup) = (\pi^\star (\rho) \cup ) \ov\U$ for all
    non-twisted degree-two cohomology classes $\rho \in H^2(\cX;\CC)$;
    here the product on the left-hand side is the Chen--Ruan orbifold
    cup product and the product on the right-hand side is the usual
    cup product
  \item[(c)] There is a point $\tau_0\in H^2_{\rm orb}(\cX;\CC)$ such
    that the standard opposite subspaces $\cH^-_{\cX}|_{Q_i=1}$ 
    and $\cH^-_{Y}|_{Q_i=1}$ 
    are opposite to $\E_{{\rm lim}, \tau_0}^\cX$ and to 
    $\ov\U\big(\E_{{\rm lim},\tau_0}^\cX \big)$ 
    respectively. 
\end{itemize}
    \end{minipage}
  \end{equation}
\end{conj}

Let $\U \colon \cHX|_{Q_i=1}\to \cHY|_{Q_i=1}$ be $\ov\U$ followed by
changing the sign of $z$.  
Conjecture \ref{conj:preciseCRC} and the definition
(\ref{eq:analyticLag}) of the analytic version of Givental's cone
immediately imply Conjecture \ref{conj} from the Introduction:
\begin{equation} 
\label{eq:U_analyticLag} 
\U \big(\cLX^{\rm an} \big) = \cLY^{\rm an}. 
\end{equation} 
Note that $\U$ determines the map $\Upsilon$ 
uniquely (see \eqref{eq:J_intersection}):
\[
\ov\U\big(\E_\tau^\cX\big) \cap \big(1+\cHY^-|_{Q_i=1}\big) = 
\{1 + \Upsilon(\tau)/z + O(1/z^2)\}.  
\]
In terms of the Lagrangian cones, the base space of 
the big A-model \VHS arises as a space parametrizing 
tangent spaces to the cone $\cLX^{\rm an}$ or $\cLY^{\rm an}$. 
From this viewpoint, the map $\Upsilon$ 
can be interpreted as the map between the 
moduli spaces of tangent spaces to the cones 
induced from the isomorphism $\U \colon \cLX^{\rm an} \cong \cLY^{\rm an}$. 

\begin{rem}
  Even though the Gromov--Witten theories of $\cX$ and $Y$ are defined
  over $\QQ$, the transformation $\U$ may only be defined over $\CC$.
  This happened here for both $\cX = \PP(1,1,2)$ and $\cX= \PP(1,1,1,3)$: see
  Propositions~\ref{pro:symplectictransformationP112}
  and~\ref{pro:symplectictransformationP1113}.
\end{rem}

\begin{rem}
  \label{rem:ambiguity}
  The operator $\U$ will be far from unique because of various
  degree-preserving symmetries of the Lagrangian cones.  The ambiguity
  of $\U$ by scalar multiplication (dilation symmetry) can be fixed by
  the condition $\U(\fun)=\fun+O(z^{-1})$, but there will also be
  discrete symmetries coming from monodromy of the mirror \VHSns.  The
  Divisor Equation implies that 
  \begin{align}
    \label{eq:divisorequationandJfunction}
    J_{\cX}(\tau+\sigma,-z)|_{Q=1} = e^{-\sigma/z} 
    J_{\cX}(\tau,-z)|_{Q^d \mapsto e^{\CHpair{\sigma}{d}}} && \sigma \in H^2(\cX;\CC)
  \end{align}
  and in particular setting $\sigma=2\pi\tti\rho$, where $\rho\in
  H^2(X,\ZZ)$ is an integral degree-two class coming from the coarse
  moduli space $X$, shows that multiplication by
  $\exp\big(2\pi\tti\rho/z\big)$ preserves $\cL_\cX$. This symmetry comes
  from a monodromy around the large radius limit point. Part (b) can
  therefore be understood as a compatibility between the monodromy
  actions on $\cHX$ and $\cHY$.  There will also be discrete symmetries of
  other types; see the discussion in Section~\ref{sec:openquestions} below.
\end{rem} 

\begin{rem} \label{rem:Tlim} 
  When $\tau\in \HorbX$ is sufficiently close to the large radius limit, 
  the limiting Hodge structure $\E^{\cX}_{\lim,\tau}$ exists 
  and is calculated as: 
  \begin{equation}
  \label{eq:limitingHodgeingeneral}
  \E^\cX_{\lim,\tau}=\lim_{\sigma\to\text{l.r.l.}} 
  e^{-\sigma/z}\E^\cX_{\tau+\sigma} 
  = \Span_{\Cz}\left\langle \parfrac{J_\cX}{\tau^\alpha}(\tau,z)\Bigr|_{Q_i=0}
               \right\rangle_{1\le\alpha\le N}
  \end{equation}
  Part (c) implies that in a neighbourhood of the large radius limit
  ``$\lim_{\sigma\to \text{l.r.l.}} (\tau_0+\sigma)$", the Frobenius
  structures associated to both $\cX$ and $Y$ are well-defined.
\end{rem}

\begin{rem} \label{rem:Novikovrecovery} We can restore the Novikov
  variables in equation (\ref{eq:U_analyticLag}) as follows.  
  Given the analytic Givental cone $\cL^{\rm
    an}$, we can define a \emph{family} of analytic cones $\cL_Q^{\rm
    an}$ parametrized by $r$ \emph{complex numbers} $Q_1,\dots, Q_r
  \in \CC^\times$:
  \[
  \cL_Q^{\rm an} := \exp\(\sum_{i=1}^r p_i \log Q_i /z\) \cL^{\rm an}
  \] 
  Due to the discrete symmetries from Remark \ref{rem:ambiguity},
  $\cL_Q^{\rm an}$ depends only on suitable roots
  $Q_1^{1/m_1},\dots,Q_r^{1/m_r}$ of the Novikov variables.  From
  equation \eqref{eq:divisorequationandJfunction}, the original
  Givental cone over the Novikov ring $\Lambda$ can be interpreted as
  the completion of the family of cones $\{\cL_Q^{\rm an}\}_{Q\in
    (\C^\times)^r}$ at the origin $Q=0$.  By using (b), we have a
  family version of (\ref{eq:U_analyticLag}):
  \[
  \U \big( \cL^{\rm an}_{\cX,Q} \big) 
   = \cL^{\rm an}_{Y,\pi^* Q}
  \]
  where $\pi^*$ denotes 
  the map between Novikov variables induced by 
  the pull back $\pi^*\colon H^2(X)\to H^2(Y)$.  
\end{rem} 

\begin{rem}
  The method by which Conjecture~\ref{conj:preciseCRC} was proved here
  (see Sections~\ref{sec:P1113} and~\ref{sec:P112}) is applicable to a
  broad class of examples.  Every time we know a mirror for small
  quantum cohomology we should have a similar explanation for
  wall-crossing phenomena in genus-zero Gromov--Witten theory.  This
  has been emphasized in work of Coates \cite{Coates:wallcrossings2}
  and Iritani \cite{Iritani:wallcrossings3}, which provides evidence
  that something very like Conjecture~\ref{conj:preciseCRC} may also
  hold for more general crepant birational transformations.
\end{rem}

\begin{rem}
  Conjecture~\ref{conj:preciseCRC} is in keeping with ideas of Ruan
  \cite{Ruan:crepant2}: that the Lagrangian cones which encode
  genus-zero Gromov--Witten invariants of an orbifold $\cX$ and its
  crepant resolution $Y$ should coincide after a symplectic
  transformation, and that the total descendant potentials of $\cX$ and
  of $Y$ (which are generating functions encoding Gromov--Witten
  invariants of all genera) should be related by the quantization of
  this symplectic transformation.  These ideas were inspired by
  results of Givental \cite{Givental:symplectic}, who has found in a
  number of examples that operations in Gromov--Witten theory which in
  genus zero give a symplectic transformation of the Lagrangian cone
  act on higher-genus invariants by applying the quantization of that
  symplectic transformation (which is a differential operator) to the
  total descendant potential.
\end{rem}

\subsection{Consequences of Conjecture~\ref{conj:preciseCRC}}
\label{sec:consequences}
When Conjecture~\ref{conj:preciseCRC} holds we can distinguish two
cases:
\begin{itemize}
\item[(i)] $\U\big(\cHX^-\big) = \cHY^-$
\item[(ii)] $\U\big(\cHX^-\big) \ne \cHY^-$
\end{itemize}
In case (i), which occurs when matrix elements of $\U$ do not contain
positive powers of $z$, the big quantum cohomology Frobenius
structures of $\cX$ and $Y$ are related by analytic continuation ---
exactly as in Section~\ref{sec:P112thm}.  In this case $\U$ can be
thought of as the well-known ambiguity of fundamental solutions in the
theory of Frobenius manifolds.  In case (ii), which occurs when some
matrix elements of $\U$ contain strictly positive powers of $z$, $\U$
does not preserve the opposite subspaces and the Frobenius manifolds
associated to $\cX$ and $Y$ will in general be different.  In this
case $\HorbX$ and $\HY$ carry the same
$F$-manifold structure but have different flat co-ordinate systems.
Note that case (i) happened for $\cX=\PP(1,1,2)$ and that case (ii)
happened for $\cX=\PP(1,1,1,3)$.

In the next two sections we will show that if $\cX$ satisfies a Hard
Lefschetz condition then case (ii) cannot occur, and thus that our
Conjecture implies the Bryan--Graber Conjecture.

\subsection{Conjecture~\ref{conj:preciseCRC} and Theorems of
  Lupercio--Poddar and Yasuda}

\begin{thm}
  \label{thm:Yasuda}
  If Conjecture~\ref{conj:preciseCRC} holds then $\HorbX$ and $\HY$
  are isomorphic as graded vector spaces.
\end{thm}

\begin{rem}
   Lupercio--Poddar and Yasuda have shown that if $\cX$ and $Y$
  are $K$-equivalent orbifolds then $\HorbX$ and $\HY$ have the same
  Hodge numbers \citelist{\cite{Yasuda}\cite{Lupercio--Poddar}}.  Thus
  this consequence itself is not surprising.  We include a proof only
  because the isomorphism which we construct depends\footnote{In fact
    the isomorphism we construct depends only on the equivalence class
    of $\tau_0$ in $H^2_{\rm orb}(\cX;\CC)/H^2(\cX;\CC)$.} on the
  choice of $\tau_0$ in Conjecture~\ref{conj:preciseCRC}(c), and so
  from the point of view of Gromov--Witten theory there may be no
  \emph{distinguished} graded isomorphism between $\HorbX$ and $\HY$.
\end{rem}

\begin{proof}[Proof of Theorem~\ref{thm:Yasuda}]

  Let $\tau_0\in H^2_{\rm orb}(\cX)$ be as in
  Conjecture~\ref{conj:preciseCRC}(c).  
  Equation
  \eqref{eq:limitingHodgeingeneral} shows that the limiting Hodge
  structure $\E^\cX_{\lim, \tau_0}$ is spanned over $\CC\{z\}$ by
  homogeneous elements of $\cHX|_{Q_i=1}$, 
  and hence that $\E^\cX_{\lim,\tau_0}$ 
  is a homogeneous subspace of $\cHX|_{Q_i=1}$.  
  Because $\cHX^-|_{Q_i=1}$ is
  homogeneous and opposite to $\E^\cX_{\lim,\tau_0}$, there is a
  graded isomorphism:
  \[
  H^\bullet_{\rm orb}(\cX) \cong z \cH_-^\cX/\cH_-^\cX\Big|_{Q_i=1} \cong 
  \E^\cX_{\lim,\tau_0}/z \E^\cX_{\lim,\tau_0}
  \] 
  On the other hand, $\U\big(\E^\cX_{\lim,\tau_0}\big)$ is also homogeneous
  since $\U$ is degree-preserving.  Using
  Conjecture~\ref{conj:preciseCRC}(c) again, there is a graded
  isomorphism:
  \[
  H^\bullet(Y) \cong z \cH_-^Y/\cH_-^Y\Big|_{Q_i=1} \cong 
  \U\big(\E^\cX_{\lim,\tau_0}\big)/z
  \U\big(\E^\cX_{\lim,\tau_0} \big) 
  \] 
  The map $\U$ induces a graded isomorphism 
   $\E^\cX_{\lim,\tau_0}/z\E^\cX_{\lim,\tau_0}
    \cong \U\big(\E^\cX_{\lim,\tau_0}\big)/
    z\U\big(\E^\cX_{\lim,\tau_0}\big)$, and so the
  conclusion follows.
\end{proof}

\subsection{A Hard Lefschetz Condition and the Bryan--Graber
  Conjecture}

\begin{thm}
  \label{thm:hardLefschetz}
  Suppose that Conjecture~\ref{conj:preciseCRC} holds, that $\cX$ has
  complex dimension $n$, and that $\cX$ in addition satisfies the Hard
  Lefschetz condition
  \[
  \text{$\omega^i\cup \colon 
  H_{\rm orb}^{n-i} (\cX) \rightarrow H_{\rm orb}^{n+i}(\cX)$
  is an isomorphism for all $i\ge 0$} 
  \]
  where $\omega \in H^2(\cX;\CC)$ is a K\"{a}hler class and $\cup$ is
  the Chen--Ruan orbifold cup product.  Then
  $\U\big(\cHX^-\big) = \cHY^-$.
\end{thm}

In view of the discussion in Section~\ref{sec:consequences}, this
implies:

\begin{cor}
  Conjecture~\ref{conj:preciseCRC} implies the revised form of the
  Bryan--Graber Conjecture \cite{Bryan--Graber}.
\end{cor}

\begin{proof}[Proof of Theorem~\ref{thm:hardLefschetz}]
  We need to show that matrix elements of $\U$ do not contain strictly
  positive powers of $z$.  Because $\U$ is a symplectic operator, the
  inverse of $\U$ is given by the adjoint $\U^\dagger$ with the sign
  of $z$ flipped.  Thus it suffices to show that matrix elements of
  $\U^{-1}$ do not contain strictly positive powers of $z$.  By taking
  a Jordan normal form of the nilpotent operator
  $\pi^\star(\omega)\cup$ on $\HY$, we obtain a basis for $\HY$ of the
  form
  \begin{align*}
    \left\{\pi^\star(\omega)^i \phi_j: 1\le j\le l, \; 0\le i\le a_j\right\} &&
    a_1\ge a_2\ge \cdots \ge a_l  
  \end{align*}
  such that $\pi^\star(\omega)^{a_j+1}\phi_j=0$.  We can assume that
  $\phi_j$ is homogeneous of degree $n-a_j+\lambda_j$ for some
  $\lambda_j\in \ZZ$.  Since $\pi^\star(\omega)$ is conjugate to $\omega$
  over $\Czzinv$, the Jordan normal forms of $\omega$ and
  $\pi^\star(\omega)$ are the same.  The hard Lefschetz condition gives
  the Lefschetz decomposition of $H_{\rm orb}^\bullet(\cX)$:
\[
H_{\rm orb}^\bullet(\cX) = \bigoplus_{i=0}^n 
\bigoplus_{k=0}^i \omega^k PH^{n-i}_{\rm orb}(\cX)
\]
where $PH^{n-i}_{\rm orb}(\cX)$ is the primitive cohomology group:
\[
PH^{n-i}_{\rm orb}(\cX) = \left\{\phi \in H_{\rm orb}^{n-i}(\cX) :
\omega^{i+1} \phi =0 \right\}
\]
The numbers $a_j$ above are determined by the Lefschetz decomposition.
The variance $v_\cX$ of the spectrum of $H_{\rm orb}^\bullet(\cX)$ is:
\[
v_{\cX}= \sum_{i=0}^{2n} (i-n)^2  \dim H_{\rm orb}^i(\cX) 
= \sum_{j=1}^l \sum_{i=0}^{a_j}({-a_j}+2i)^2
\]
On the other hand, the variance $v_Y$ of $H^\bullet(Y)$ is:
\[
v_Y=\sum_{j=1}^l \sum_{i=0}^{a_j}\big({-a_j}+\lambda_j+2i\big)^2 
=v_{\cX}+ \sum_{j=1}^l (1+a_j) \lambda_j^2 
\]
Since there is a graded isomorphism $H_{\rm orb}^\bullet(\cX)\cong
H^\bullet(Y)$, we have $v_\cX=v_Y$ and so $\lambda_j=0$ for all $j$.
Thus $\deg \phi_j= n-a_j$.  Then $\U^{-1}(\phi_j)$ is in
$\Ker(\omega^{a_j+1})$ and also of degree $n-a_j$.  Using the
Lefschetz decomposition of $\HorbX$ again, we see that
$\U^{-1}(\phi_j)$ does not contain positive powers of $z$.  Thus
$\U^{-1}\big(\pi^\star(\omega)^i \phi_j\big)= 
\omega^i \, \U^{-1}(\phi_j)$ does not 
contain positive powers of $z$ either.
\end{proof} 

\begin{rem}
  Fernandez \cite{Fernandez} has shown that the Hard Lefschetz
  condition in Theorem~\ref{thm:hardLefschetz} is equivalent to the
  equality $\age \cX_i = \age I(\cX_i)$ for all components $\cX_i$ of
  the inertia stack $\cI\cX$.  This condition holds for $\PP(1,1,2)$
  and for any other two-dimensional Gorenstein orbifold but not for
  $\PP(1,1,1,3)$.
\end{rem} 

\begin{rem} 
Since the paper was written, 
a more general form of the Hard Lefschetz condition 
applicable to the case of partial resolutions 
or $K$-equivalence is studied by Iritani 
\cite[Section~3.7]{Iritani:RIMS}. 
\end{rem} 

\subsection{Conjecture~\ref{conj:preciseCRC} and the Ruan Conjecture}
\label{sec:RuanConjecture}
As we have seen in Theorems~\ref{thm:P112} and~\ref{thm:P1113}, by
proving Conjecture~\ref{conj:preciseCRC} we also proved the Ruan
Conjecture for $\cX = \PP(1,1,2)$ and $\cX = \PP(1,1,1,3)$.  This is
slightly misleading, however, as in general our Conjecture only
implies a modified version of the Ruan Conjecture.  This is explained
in detail in Coates--Ruan \cite{Coates--Ruan}*{Section~8}.  Coates
\cite{Coates:wallcrossings2} has proved our Conjecture, and hence the
modified Ruan Conjecture, in an example for which the modified Ruan
Conjecture and the original Ruan Conjecture differ: this example is
the canonical bundle to $\PP(1,1,3)$.  We expect that the original
version of the Ruan Conjecture is false in general.

\subsection{Open Questions}
\label{sec:openquestions}
We close by indicating several questions which deserve further study.
One such direction involves real and integral
structures\footnote{Since this paper was written, this question has
  been studied by Iritani \cite{Iritani:integral}.} on the \VHSns.
The B-model \VHS has a natural integral structure, coming from the
lattice of Morse cycles, but this is hard to see in the A-model.
The study of real structures should lead to $tt^*$-geometry of the
K\"{a}hler moduli space (see Cecotti--Vafa \cite{Cecotti--Vafa},
Dubrovin \cite{Dubrovin-fusion}, and Hertling \cite{Hertling}).  There
should also be a hidden real structure on Givental's symplectic space
$\cH$, and the symplectic transformation $\U$ from
Conjecture~\ref{conj:preciseCRC} should preserve such real structures.
The specialization of K\"{a}hler parameters to \emph{purely imaginary}
numbers (such as $\tau^1=2\pi\tti$,
$\tau^2={-\pi\tti}+\frac{1}{2}\tau$ in the case $\cX =
\PP(1,1,2)$, $Y=\FF_2$) might be explained from this viewpoint, as
they can be read off from $\U(\fun)$.

One should also consider higher-genus Gromov--Witten invariants.  In
Givental's quantization formalism \cite{Givental:symplectic}, the
total descendant potential corresponds to a quantization of the
(genus-zero) Lagrangian cones and lives in a Fock space produced from
$\cH$.  This suggests that the total descendant potentials of $\cX$
and $Y$ are related by
\begin{equation}
\label{eq:highergenusgeneralization}
\cD_{Y}\propto \widehat{\U}\cD_{\cX}. 
\end{equation}
We do not give rigorous meaning to this formula here; the idea is that
the two elements $\cD_{Y}$, $\cD_{\cX}$ of different Fock spaces will
be projectively identified by the change of polarization $\U$.  More
than a decade ago, Witten \cite{Witten} introduced a quantum
mechanical system on $H^3(X)$ for a Calabi--Yau threefold $X$ and
showed that the total potential of the B-model behaves like a wave
function of this quantum system.  Formula
\eqref{eq:highergenusgeneralization} fits with this picture.  It again
matches well with the ideas of Ruan disussed above, and also with
recent work of Aganagic, Bouchard, and Klemm \cite{ABK}.  They argue
that the fundamental group of the B-model moduli space should act as
``quantum symmetries'' of the total descendant potential.  The
monodromy around the large radius limit point is, as discussed in
Remark \ref{rem:ambiguity}, related to the Divisor Equation in
Gromov--Witten theory.  When we have a crepant resolution
$Y\rightarrow X$ of $\cX$, there should also be an ``extra'' monodromy
action on the Gromov--Witten theory of $Y$ coming from orbifold loops
around the large radius limit point for $\cX$; such monodromy will not
in general preserve the opposite subspace $\cHY^-$ for $Y$.  We hope
that these symmetries together with a hidden integral structure will
reveal a kind of \emph{quantum automorphic property} of the potential
$\cD_Y$.

\appendix
\section*{Appendix A} 
\label{sec:appendix}
In this appendix we give a brief account of the analytic continuation
of the $I$-function performed in \eqref{eq:IF3continued}.  We use an
integral representation of Barnes type, following Candelas \emph{et
  al.}  \cite{CDGP} and Horja \cite{Horja}.

Set $\ap_1=p_2-3p_1$, as in Figure \ref{fig:secondaryF_3}.  During the
analytic continuation we regard $p_1$ and $p_2$ as complex variables
and consider the $I$-function as an analytic function in $y_1$, $y_2$,
$p_1$, $p_2$, and $z$.  We obtain cohomology classes by, at the end of
the process, Taylor expanding in $p_1$ and $p_2$ and then regarding
$p_1$ and $p_2$ as cohomology classes.  We have:
\begin{multline}
 \label{eq:IF3}
 I_{\FF_3}(y_1,y_2,z) = z \, \textstyle
 \Gamma\big(1+\frac{p_1}{z}\big)^3 \, 
 \Gamma\big(1+\frac{p_2}{z}\big)
 \Gamma\big(1+\frac{\ap_1}{z}\big) \\
 \times \sum_{n,m\ge 0}
 \frac{
   y_1^{n+p_1/z}y_2^{m+p_2/z}z^{-2m}
 }
 {
   \Gamma\big(1+\frac{p_1}{z}+n\big)^3 \,
   \Gamma\big(1+\frac{p_2}{z}+m\big)
   \Gamma\big(1+\frac{\ap_1}{z}+m-3n\big)
 }
\end{multline}
Since $\Gamma(z)\Gamma(1-z)=\pi/\sin(\pi z)$ the coefficient of
$y_2^{m+p_2/z}$ in \eqref{eq:IF3} can be written, neglecting several
Gamma factors and powers of $z$, as:
\[
\sum_{n\ge 0} 
\frac{
  \Gamma\big({-\frac{\ap_1}{z}}+3n-m\big)
}
{
  \Gamma\big(1+\frac{p_1}{z}+n\big)^3
} 
\frac{
  \sin\big({-\frac{\ap_1}{z}}\pi+3n\pi-m\pi\big)
}
{
  \pi
} 
y_1^{n+p_1/z}
\]
This is the sum of residues:
\[
(-1)^m
\frac{
  \sin\big({-\frac{\ap_1}{z}}\pi\big)
}
{
  \pi
}
\sum_{n\ge 0} \Res_{s=n} 
\left[
\Gamma(s)\Gamma(1-s)  
\frac{
  \Gamma\big({-\frac{\ap_1}{z}}+3s-m\big)
}
{
  \Gamma\big(1+\frac{p_1}{z}+s\big)^3
} 
y_1^{s+p_1/z} \, ds
\right]
\]
For $|y_1|<{1 \over 27}$, it can be replaced by the integral along a
contour $C_m$ from $s = \tti\infty$ to $s = -\tti\infty$ which runs
along the imaginary axis for $|s|$ large and is such that $s=0,1,2,\dots$
are on the right hand side of $C_m$ and that $s=-1,-2,-3,\dots$ and
$s=\frac{m}{3}+\frac{\ap_1}{3z},
\frac{m}{3}+\frac{\ap_1}{3z}-\frac{1}{3},
\frac{m}{3}+\frac{\ap_1}{3z}-\frac{2}{3},\dots$ are on the left hand
side of $C_m$:
\[
(-1)^m
\frac{
  \sin\big({-\frac{\ap_1}{z}}\pi\big)
}
{
  \pi
} 
\frac{1}{2\pi\tti}
\int_{C_m} 
\frac{
  \Gamma\big({-\frac{\ap_1}{z}}+3s-m\big)
  \Gamma(s)
  \Gamma(1-s)
}
{
  \Gamma\big(1+\frac{p_1}{z}+s\big)^3
} 
y_1^{s+p_1/z}
\, ds 
\]
This integral converges on the region $| \arg(y_1)| < \pi$; see
e.g. Horja \cite[Lemma 3.3]{Horja}.  For $|y_1|>{1 \over 27}$ we can
close the contour to the left, finding:
\begin{multline*}
  (-1)^m
  \frac{
    \sin\big(\frac{\ap_1}{z}\pi\big)
  }
  {
    \pi
  } 
  \sum_{n\ge 0} \Res_{
    s=\frac{m}{3}+\frac{\ap_1}{3z}-\frac{n}{3}} 
  \left[
    \frac{
      \Gamma\big({-\frac{\ap_1}{z}}+3s-m\big)
      \Gamma(s)
      \Gamma(1-s)
    }
    {
      \Gamma\big(1+\frac{p_1}{z}+s\big)^3
    } 
    y_1^{s+p_1/z}\, ds
  \right]
  \\
  +   (-1)^m
  \frac{
    \sin\big(\frac{\ap_1}{z}\pi\big)
  }
  {
    \pi
  } \sum_{n\ge 0}
  \Res_{s=-1-n} 
  \left[
    \frac{
      \Gamma\big({-\frac{\ap_1}{z}}+3s-m\big)
      \Gamma(s)
      \Gamma(1-s)
    }
    {
      \Gamma\big(1+\frac{p_1}{z}+s\big)^3
    }
    y_1^{s+p_1/z} \, ds 
  \right]
\end{multline*}
The residues at $s=-1-n$ vanish in cohomology, as $p_1^3=0$, so
this is:
\[
\sum_{n\ge 0}
\frac{(-1)^{m+n} }{n!}
\frac{
  \sin\big(\frac{\ap_1}{z}\pi\big)
}
{
  3 \sin\big(\frac{\ap_1}{3z}\pi+\frac{m-n}{3}\pi\big)
} 
\frac{
  y_1^{(m-n)/3+p_2/(3z)}
}
{
  \Gamma\big(1+\frac{p_2}{3z}+\frac{m-n}{3}\big)^3
} 
\]
Changing variables from $(y_1,y_2)$ to $(\aq_1,\aq_2)$ yields
\eqref{eq:IF3continued}.

\appendix
\section*{Appendix B} 
Here we give an alternative elementary proof 
of the crepant resolution conjecture (Theorem \ref{thm:P112}) 
for $\PP(1,1,2)$ without using the notion of \VHSns. 
The method here, however, does not work 
beyond toric surfaces. 

Recall the Landau-Ginzburg mirrors of $\FF_2$ and 
$\PP(1,1,2)$ described in Section \ref{sec:P112}. 
The fact that the Landau-Ginzburg mirrors of 
$\FF_2$ and $\PP(1,1,2)$ are contained in the 
family $\{(Z_y, W_y)\}_{y\in \cM_{\rm B}}$ of 
Landau-Ginzburg models shows that 
quantum cohomology of $\FF_2$ and $\PP(1,1,2)$ 
are deformation equivalent,  
since the quantum cohomology is 
isomorphic to the family of Jacobi rings of $W_y$ 
under mirror symmetry.  
To show Theorem \ref{thm:P112}, 
we need to compare the  
flat co-ordinates and trivializations of 
$\FF_2$ and $\PP(1,1,2)$. 
We will do this by an elementary mirror analysis. 

Let
\[
\cJ = \bigcup_{y \in \cM_{\rm B}} J(W_y), \quad 
J(W_y) = \C[Z_y]/\langle dW_y \rangle  
\]
be the bundle of Jacobi rings and 
define a Kodaira--Spencer map 
$\KS\colon T\cM_{\rm B} \rightarrow \cJ$ by
\[
\KS\colon T\cM_{\rm B}  \ni v \longmapsto \tilde{v}(W) \in \cJ,
\]
where $\tilde{v}$ is a vector field on $Z$ such that
$d\pi(\tilde{v})=v$.  
This does not depend on the choice of lift $\tilde{v}$. 
The map $\KS$ is the Kodaira--Spencer map
for the B-model \VHS discussed after 
Proposition~\ref{pro:normalformofconn}. 
This corresponds under mirror symmetry to the inclusion 
$H^2({\FF_2})\hookrightarrow H^\bullet({\FF_2})$. 
By pulling back the residue metric on $\cJ$, we
obtain a non-degenerate symmetric $\C$-bilinear metric on 
$T\cM_{\rm  B}$.  
This turns out to be flat on $\cM_{\rm B}$.  
If one considers the frame of vector fields
\[
\varphi_1=\sqrt{1-4 y_1}\partial_1 + \frac{1}{2}(1-\sqrt{1-4 y_1})\partial_2, \quad
\varphi_2=\partial_2, 
\]
where $\partial_i=y_i\parfrac{}{y_i}$ --- this frame also appears in
\cite[Example 5.4]{Guest} --- then an easy calculation shows that
$\pairW{\KS(\varphi_i)}{\KS(\varphi_j)}$ is the non-degenerate
constant matrix
\[
\begin{pmatrix}
0 & 1 \\
1 & 2
\end{pmatrix},
\]
and that $[\varphi_1,\varphi_2]=0$.  Thus $\varphi_1, \varphi_2$ are
flat.  

One can find a co-ordinate system 
$(q_1,q_2)$ on a neighbourhood of $(y_1,y_2) = (0,0)$ 
in $\cM_{\rm B}$ such that $\varphi_i=q_i\parfrac{}{q_i}$:
\[
y_1=\frac{q_1}{(1+q_1)^2}, \quad y_2=q_2(1+q_1); \qquad 
\aq_1=\frac{1+q_1}{q_1^{1/2}},\quad \aq_2=q_1^{1/2}q_2.
\]
These flat co-ordinates $(q_1,q_2)$ are exactly 
the mirror map (\ref{eq:F2mirrormap}) for $\FF_2$
and the frame $\varphi_1$, $\varphi_2$ 
corresponds to the basis $p_1,p_2$ of $H^2({\FF_2})$.  
In the flat co-ordinates, the mirror locus 
$\cM_{\PP(1,1,2)}=\{\aq_1=0\}$ of $\PP(1,1,2)$ 
is given by the equation $q_1=-1$. 
Take a path from a neighbourhood of $y_1=y_2=0$ 
to that of $\aq_1=\aq_2=0$ and choose a branch of 
$\sqrt{q_1}$ as follows:
\[
\sqrt{q_1}= - \tti \lambda, 
\quad q_2=\text{small const}, \quad 
\lambda\in [0,1]. 
\]
Following this path, we see that 
$(H_{\rm orb}^\bullet(\PP(1,1,2)),\circ_q)$ is
obtained from $(H^\bullet({\FF_2}),\circ_q)$ 
by analytic continuation followed by the specialization
$q_1=-1, q_2=\tti \aq_2 = \tti \sqrt{q}$. 
When $\lambda \to 1$, 
the flat frame $\varphi_1$, $\varphi_2$ goes to
\[
\varphi_1\to -\tti \parfrac{}{\aq_1} + 
\frac{1}{2}\aq_2\parfrac{}{\aq_2}, \quad 
\varphi_2\to \aq_2\parfrac{}{\aq_2}. 
\]
But the elements
\[
\parfrac{}{\aq_1}W_y \Big|_{\aq_1=0}, \quad 
\frac{1}{2}\aq_2\parfrac{}{\aq_2} W_y\Big |_{\aq_1=0} 
\]
in the Jacobi ring $J(W_y)$ on the locus $\{\aq_1=0\}$ correspond
under mirror symmetry to the elements
\[
\fun_{1 \over 2},\quad p 
\]
of the small quantum orbifold cohomology algebra 
of $\PP(1,1,2)$. 
This explains the correspondence of the basis
\[
p_1 \mapsto -\tti  \fun_{1 \over 2}+ p, \quad p_2 \mapsto 2p. 
\] 
By construction, this correspondence preserves the algebra structure
and the Poincar\'{e} metric.

%
%
%
%
\bibliographystyle{amsxport}

\end{document}